\documentclass[11pt,final]{ucthesis}
\usepackage{amsmath}
\usepackage{amssymb}
\usepackage{amsfonts}
\usepackage{fancyheadings}
\usepackage{amscd}
\usepackage{latexsym}
\usepackage{theorem}
\usepackage{epsfig}
\usepackage{psfrag}
\usepackage{bbold}

\newtheorem{theorem}{Theorem}[chapter]
\newtheorem{lemma}[theorem]{Lemma}
\newtheorem{proposition}[theorem]{Proposition}
\newtheorem{corollary}[theorem]{Corollary}
\newtheorem{claim}[theorem]{Claim}

{\theorembodyfont{\rmfamily}
\theoremstyle{plain}
\newtheorem{definition}[theorem]{Definition}
\newtheorem{convention}[theorem]{Notational Convention}
\newtheorem{example}[theorem]{Example}
\newtheorem{problem}[theorem]{Problem}
\newtheorem{remark}[theorem]{Remark}
}

\input{xy}
\xyoption{all}

\newcommand{\Vol}{\operatorname{Vol}}

\renewcommand{\Im}{\operatorname{Im}}

\newcommand{\id}{\operatorname{id}}

\renewcommand{\ker}{\operatorname{Ker}}
\renewcommand{\dim}{\operatorname{dim}}
\newcommand{\Ann}{\operatorname{Ann}}
\newcommand{\Sym}{\operatorname{Sym}}
\newcommand{\Hom}{\operatorname{Hom}}
\newcommand{\tr}{\operatorname{tr}}

\newcommand{\Ld}{\displaystyle\lim}
\newcommand{\Od}{\displaystyle\bigoplus}
\newcommand{\Pd}{\displaystyle\prod}

\newcommand{\C}{{\mathbb{C}}}
\newcommand{\Z}{{\mathbb{Z}}}
\newcommand{\Q}{{\mathbb{Q}}}
\renewcommand{\H}{{\mathbb{H}}}
\newcommand{\R}{{\mathbb{R}}}
\newcommand{\Zt}{{\Z_2}}

\newcommand{\otn}{\{1,\ldots,n\}}
\newcommand{\half}{\frac{1}{2}}
\newcommand{\bd}{\partial}

\newcommand{\bigmid}{\hs\Big{|}\hs}
\newcommand{\subs}{\subseteq}
\newcommand{\sups}{\supseteq}
\newcommand{\hookto}{{\hookrightarrow}}
\newcommand{\onto}{{\twoheadrightarrow}}
\renewcommand{\iff}{\Leftrightarrow}
\newcommand{\impl}{\Rightarrow}
\newcommand{\becircled}{\mathaccent "7017}
\newcommand{\hs}{\hspace{3pt}}
\renewcommand{\ss}{\substack}

\newcommand{\D}{\Delta}
\renewcommand{\i}{\iota}
\renewcommand{\a}{\alpha}

\newcommand{\eps}{\varepsilon}
\renewcommand{\k}{\kappa}
\newcommand{\K}{\kappa}
\newcommand{\ga}{\gamma}
\renewcommand{\t}{\tau}

\renewcommand{\cot}{T^*\C^n}
\newcommand{\Cn}{\C^n}

\newcommand{\Mso}{\Ma^{\so}}
\newcommand{\Macs}{\Ma^{\gm}}
\newcommand{\Ma}{\M}
\newcommand{\Xa}{\X}
\newcommand{\Gap}{\Gamma_p}
\newcommand{\Gapt}{\widehat{\Gamma}_p}
\newcommand{\Gapto}{\Gapt\!\becircled{}}

\newcommand{\muc}{\mu_{\C}}
\newcommand{\mur}{\mu_{\R}}
\newcommand{\barmuc}{\bar{\mu}_{\C}}
\newcommand{\barmur}{\bar{\mu}_{\R}}
\newcommand{\mr}{\mur}
\newcommand{\mc}{\muc}
\newcommand{\mhk}{\mu_{\text{HK}}}

\newcommand{\ka}{K_{\a}}

\renewcommand{\mod}{{\!/\!\!/}}
\newcommand{\mmod}{{\!/\!\!/\!\!/\!\!/}}

\newcommand{\SL}{{SL(2,\C)}}
\newcommand{\sutd}{{\mathfrak{su}(2)}^*}
\newcommand{\sut}{{\mathfrak{su}(2)}}
\newcommand{\slt}{{\mathfrak{sl}(2,\C)}}
\newcommand{\so}{S^1}
\newcommand{\gm}{\C^{\!\times}}

\newcommand{\GC}{G_{\C}}
\newcommand{\gdc}{\gd_{\C}}
\newcommand{\gd}{\g^*}
\newcommand{\g}{\mathfrak{g}}
\newcommand{\Tk}{T^k}
\newcommand{\Tn}{T^n}
\newcommand{\Td}{T^d}

\newcommand{\tk}{\mathfrak{t}^k}
\newcommand{\tn}{\mathfrak{t}^n}
\newcommand{\td}{\mathfrak{t}^d}
\newcommand{\tkd}{(\tk)^*}
\newcommand{\tndz}{(\tn_{\Z})^*}
\newcommand{\tnz}{\tn_{\Z}}
\newcommand{\tdz}{\td_{\Z}}
\newcommand{\tnd}{(\tn)^*}
\newcommand{\tdd}{(\td)^*}
\newcommand{\tdu}{\mathfrak{t}^*}

\newcommand{\E}{\mathcal{E}}
\newcommand{\EA}{\E_A}

\newcommand{\Af}{\operatorname{Rep}(Q)}
\newcommand{\Diag}{\operatorname{Diag}}

\newcommand{\ktd}{\k_{\Td}}
\newcommand{\ktdso}{\k_{\Td\times\so}}
\newcommand{\kso}{\k_{\so}}

\newcommand{\Tdr}{T_{\R}^d}
\newcommand{\Tdso}{\Td\times S^1}

\newcommand{\hma}{H^*(\Ma)}
\newcommand{\htm}{H^*_{\Td}(\Ma)}
\newcommand{\htsm}{H^*_{\Td\times S^1}(\Ma)}
\newcommand{\hsm}{H^*_{S^1}(\Ma)}
\newcommand{\hr}{H^*_{\Tdr}(\M_{\R};\Zt)}
\newcommand{\hrs}{H^*_{\Tdr\times\Zt}(\M_{\R};\Zt)}

\newcommand{\bomp}{\C^d\setminus\cup_{i=1}^n H_i^{\C}}

\newcommand{\A}{\mathcal{A}}
\newcommand{\MA}{\mathcal{M}(\A)}
\newcommand{\hz}{H_{\Zt}^*(pt)}

\newcommand{\ctn}{{\C^{2n}}}
\newcommand{\ses}{st}

\newcommand{\polsa}{\X_S}
\newcommand{\pola}{\X}
\newcommand{\mcs}{\mathcal{S}}
\newcommand{\sa}{\mcs}
\newcommand{\spa}{\mcs'}
\newcommand{\core}{\mathfrak{L}}

\newcommand{\ms}{m_S}
\newcommand{\ns}{n_S}
\newcommand{\bi}{d_i}
\newcommand{\bk}{d_k}
\newcommand{\bj}{d_j}
\newcommand{\ba}{d_A}
\newcommand{\bns}{d_{\ns}}

\newcommand{\ck}{c_k}

\newcommand{\Sb}{\bar{S}}
\newcommand{\Lb}{\overline{S^c}}
\newcommand{\Ac}{A^c}
\newcommand{\asb}{A\cap\Sb}
\newcommand{\acsb}{\Ac\cap\Sb}
\newcommand{\pis}{\prod_{i\in\Sb}}
\newcommand{\pjl}{\prod_{j\in\Lb}}
\newcommand{\vs}{v_S}

\newcommand{\twn}{\{2,\ldots,n\}}

\newcommand{\us}{U_S}
\newcommand{\mtus}{X_T\cap\us}
\newcommand{\ts}{T\sups S}
\newcommand{\sz}{|S|-1}

\newcommand{\Kh}{\hat\k}
\newcommand{\Hn}{\H^n}
\newcommand{\rgt}{r_T^G}

\newcommand{\IM}{\int_M}
\newcommand{\IN}{\int_N}
\newcommand{\IG}{\int_G}
\newcommand{\IF}{\int_F}
\newcommand{\intxt}{\int_{X\mmod T}}
\newcommand{\intxg}{\int_{X\mmod G}}
\newcommand{\I}{\mathcal{I}}
\newcommand{\J}{\mathcal{J}}
\newcommand{\imkt}{\left(\Im\K_T\right)}
\newcommand{\hso}{H_{\so}^*}
\newcommand{\hhso}{\widehat{H}_{\so}^*}
\newcommand{\hmg}{H_{\so}^*(M\mmod G)}

\newcommand{\hxg}{H_{\so}^*(X\mmod G)}
\newcommand{\hxt}{H_{\so}^*(X\mmod T)}
\newcommand{\hhxg}{\widehat{H}_{\so}^*(X\mmod G)}
\newcommand{\hhxt}{\widehat{H}_{\so}^*(X\mmod T)}

\newcommand{\hhsg}{\widehat{H}_{\so\times G}^*(X)}
\newcommand{\hp}{H_{S^1}^*(pt)}
\newcommand{\hhp}{{\mathbb K}}
\newcommand{\hm}{\hso(M)}
\newcommand{\hhn}{\hhso(N)}
\newcommand{\hhm}{\hhso(M)}
\newcommand{\X}{{\mathfrak{X}}}
\newcommand{\M}{{\mathfrak{M}}}

\newcommand{\NN}{\mathfrak{N}}
\newcommand{\dso}{\Diag_*(1)}
\newcommand{\muh}{\mu_{G}}
\newcommand{\mubh}{\mu_{T}}
\newcommand{\bigmod}{{\Big/\!\!\!\!\Big/}}
\newcommand{\bigmmod}{{\Big/\!\!\!\!\Big/\!\!\!\!\Big/\!\!\!\!\Big/}}

\newcommand{\ra}{\rangle}

\newcommand{\PrA}{P^r_{\! A}}
\renewcommand{\Pr}{P^r}
\newcommand{\Dr}{\D^r}
\newcommand{\DrA}{\D^r_A}

\newcommand{\rt}{\tilde r}
\newcommand{\odra}{\mathbf{1}_{\DrA}}

\newcommand{\WA}{W_A}
\newcommand{\WAbd}{W_A^{bd}}
\newcommand{\Fbd}{\mathcal{F}^{bd}}
\newcommand{\F}{\mathcal{F}}

\newcommand{\<}{\left<}
\renewcommand{\>}{\right>}
\newcommand{\la}{\langle}
\renewcommand{\ra}{\rangle}
\renewcommand{\(}{\left(}
\renewcommand{\)}{\right)}

\newcommand{\qed}{\hfill \mbox{$\Box$}\medskip\newline}
\newenvironment{proof}{\noindent {\bf Proof:}}{\qed \par}

\newenvironment{proofcomb}{\noindent {\bf Proof of \ref{comb}:}}{\qed\par}
\newenvironment{proofint}{\noindent {\bf Proof of \ref{int}:}}{}
\newenvironment{proofintegration}{\noindent {\bf Proof of \ref{integration}:}}{\qed \par}
\newenvironment{proofordinary}{\noindent {\bf Proof of \ref{ordinary}:}}{\qed \par}
\newenvironment{proofhtsm}{\noindent {\bf Proof of \ref{htsm}:}}{\qed \par}
\newenvironment{proofp}{\noindent {\bf Proof of \ref{placement}:}}{\qed \par}
\newenvironment{proofc}{\noindent {\bf Proof of \ref{coorientation}:}}{\qed \par}

\newenvironment{proofmain}{\noindent {\bf Proof of \ref{main}:}}{\qed \par}
\newenvironment{prooftech}{\noindent {\bf Proof of \ref{tech}:}}{\qed \par}
\newenvironment{componentproof}{\noindent {\bf Proof of \ref{component}:}}{\qed \par}
\newenvironment{usproof}{\noindent {\bf Proof of \ref{us}:}}{\qed \par}

\newenvironment{eqcoreproof}{\noindent {\bf Proof of \ref{eqcore}:}}{\qed \par}

\begin{document}

\title{Hyperk\"ahler Analogues of K\"ahler Quotients}
\author{Nicholas James Proudfoot}
\degreesemester{Spring}
\degreeyear{2004}
\degree{Doctor of Philosophy}
\chair{Professor Allen Knutson}
\othermembers{Professor Robion Kirby \\
  Professor Robert Littlejohn}
\numberofmembers{3}
\prevdegrees{A.B. (Harvard University) 2000}
\field{Mathematics}
\campus{Berkeley}
 
\maketitle
\approvalpage
\copyrightpage
\pagestyle{fancyplain}

\begin{abstract}
Let $\X$ be a K\"ahler manifold that is presented as a K\"ahler
quotient of $\Cn$ by the linear action of a compact group $G$.
We define the {\em hyperk\"ahler analogue} $\M$ of $\X$ as a
hyperk\"ahler quotient of the cotangent bundle $\cot$ by the induced
$G$-action.  Special instances of this construction include
hypertoric varieties \cite{BD,K1,HS,HP1} and quiver varieties \cite{N1,N2,N3}.
One of our aims is to provide a unified treatment of these
two previously studied examples.

The hyperk\"ahler analogue $\M$ is noncompact, but this noncompactness
is often ``controlled'' by an action of $\gm$ descending from
the scalar action on the fibers of $\cot$.  Specifically, we are
interested in the case where the moment map for the action
of the circle $\so\subs\gm$
is proper.  In such cases, we
define the {\em core} of $\M$, a reducible, compact subvariety onto which
$\M$ admits a circle-equivariant deformation retraction.  One of the components
of the core is isomorphic to the original K\"ahler manifold $\X$.
When $\X$ is a moduli space of polygons in $\R^3$, we interpret
each of the other core components of $\M$ as related polygonal moduli spaces.

Using the circle action with proper moment map, we define an integration
theory on the circle-equivariant cohomology of $\M$, 
motivated by the well-known localization theorem 
of \cite{AB} and \cite{BV}.  This allows us to prove a hyperk\"ahler
analogue of Martin's theorem \cite{Ma}, which describes the cohomology
ring of an arbitrary K\"ahler quotient in terms of the 
cohomology of the quotient by a maximal torus.
This theorem, along with a direct analysis of the equivariant cohomology
ring of a hypertoric variety, gives us a method for computing
the equivariant cohomology ring of many hyperk\"ahler analogues,
including all quiver varieties.
\abstractsignature
\end{abstract}

\pagestyle{plain}
\begin{frontmatter}
\pagenumbering{roman}
 
\tableofcontents
\lhead[\fancyplain{}{}]{\fancyplain{}{}}
\rhead[\fancyplain{}{}]{\fancyplain{}{}}

\begin{acknowledgements}
During my years at Berkeley, I have had countless invaluable conversations
with other students, post-docs, professors, and visitors to Berkeley and MSRI.
Many of these conversations
have been relevant to the work described here, and many more
have greatly influenced my general mathematical perspective.
I would like to particularly acknowledge my coauthors, Megumi Harada
and Tam\'as Hausel, and my advisor, Allen Knutson.
I feel extremely lucky to be part of such a warm and lively community.
\end{acknowledgements}

\end{frontmatter}

\pagestyle{fancy}
\pagenumbering{arabic}
\setlength{\headrulewidth}{1pt}
\lhead[\rm\thepage]{\sl\leftmark}
\rhead[\sl\leftmark]{\rm\thepage}
\renewcommand{\chaptermark}[1]{\markboth{\chaptername\hs\thechapter.  #1}{ }}
\cfoot{}

\begin{chapter}{Introduction}
We begin with a quick overview of some of the structures that we will consider
in this thesis, and the types of questions that we will be asking.  
Detailed definitions will be deferred until the next chapter.

Let $G$ be a compact Lie group acting linearly on $\Cn$,
with moment map $\mu:\Cn\to\gd$, and suppose we are given 
a central (i.e. $G$-fixed) regular value $\a\in\gd$ of $\mu$.  
From this data, we may define the K\"ahler 
quotient $$\X := \Cn\mod G = \mu^{-1}(\a)/G,$$ which itself inherits the structure of a 
K\"ahler manifold. 
(We may also think of $\X$ as the geometric invariant theory
quotient of $\Cn$ by $G^{\C}$ in the sense of Mumford \cite[\S 8]{MFK};
see Proposition \ref{hkhs}.)
A {\em hyperk\"ahler} manifold is a riemannian manifold $(M,g)$ equipped
with three orthogonal complex structures $J_1, J_2, J_3$ and three compatible symplectic
forms $\omega_1, \omega_2, \omega_3$ such that $J_1J_2= -J_2J_1 = J_3$
for $i=1,2$, and $3$.
The cotangent bundle $\cot$ has a natural hyperk\"ahler structure,
and this structure is preserved by the induced action of $G$.
Furthermore, there exist maps $\mu_i:\cot\to\gd$ for $i=1,2,3$
such that $\mu_i$ is a moment map for the action of $G$ with respect
to the symplectic form $\omega_i$.
We define the {\em hyperk\"ahler analogue} of $\X$
to be the hyperk\"ahler quotient
$$\M := \cot\mmod G = \Big(\mu_1^{-1}(\a)\cap\mu_2^{-1}(0)
\cap\mu_3^{-1}(0)\Big)\Big/G.$$
(The set
$\mu_2^{-1}(0)\cap\mu_3^{-1}(0)\subs\cot$
is a complex subvariety with respect to $J_1$, and $\M$
may be thought of as the geometric invariant theory quotient
of this variety by $G^{\C}$.)
The quotient $\M$ is a complete hyperk\"ahler manifold \cite{HKLR},
containing $T^*\X$ as a dense open subset (see Proposition \ref{denseopen}).
The following is a description of some well-known classes of K\"ahler quotients,
along with their hyperk\"ahler analogues.

\vspace{.2cm}
\noindent{\bf Toric and hypertoric varieties.}
These examples comprise the case where $G$ is abelian.  
The geometry of toric varieties is deeply related to the combinatorics of polytopes;
for example, Stanley \cite{St} used the hard Lefschetz theorem for toric varieties
to prove certain inequalities for the $h$-numbers of a simplicial polytope.
Hypertoric varieties, introduced by Bielawski and Dancer \cite{BD}, interact in a similar
way with the combinatorics of real hyperplane arrangements.
Following Stanley's work, Hausel and Sturmfels \cite{HS} used the hard Lefschetz theorem
on a hypertoric variety to give a geometric interpretation of some previously-known
inequalities for the $h$-numbers of a rationally representable matroid.
In Chapter \ref{ht} we will explore further combinatorial properties
of the various equivariant cohomology rings of a hypertoric variety.

\vspace{.2cm}
\noindent{\bf Quiver varieties.}
For any directed graph,
Nakajima \cite{N1,N2,N3} defined a {\em quiver variety}
to be the hyperk\"ahler analogue of the moduli space of framed representations of that graph
(see Section \ref{quiver}).
Examples include the Hilbert scheme of $n$ points in $\C^2$ \cite{N4}, the moduli space
of instantons on an ALE space \cite{N1}, and Konno's hyperpolygon spaces \cite{K2,HP2},
which are the hyperk\"ahler analogues of the moduli spaces of $n$-sided polygons in $\R^3$ with
fixed edge lengths.
Quiver varieties have received much attention from representation theorists due to the actions
of various infinite-dimensional Lie algebras on the cohomology and
K-theory of a quiver variety (see, for example, \cite{N3,N5}).

\vspace{.2cm}
\noindent{\bf Moduli spaces of bundles and connections.}
Narasimhan and Seshadri \cite{NS} defined a notion of stability for a vector bundle
on a Riemann surface $\Sigma$, and proved that the moduli space of stable 
holomorphic bundles on $\Sigma$ may be identified with the moduli space of irreducible, 
flat, unitary connections.  
Atiyah and Bott presented this space as a K\"ahler quotient of the affine
space of all connections on a fixed bundle $E$ by the gauge group
of automorphisms of $E$.  This picture can be complexified by replacing holomorphic
bundles with Higgs bundles, and unitary connections with arbitrary ones \cite{Hi}.
The correspondence between Higgs bundles, flat connections, and representations
of the fundamental group is known as {\em nonabelian Hodge theory}, and has
been studied and generalized extensively by Simpson \cite{Si} in addition to
many other authors.
These constructions involve taking a quotient of an infinite dimensional affine
space by an infinite dimensional group, and therefore lie
it is beyond the scope of this work.
Many of our techniques, however, can be applied in this context.
See for example \cite{H1, HT1, HT2}.

\vspace{.2cm}
Consider the action of the multiplicative group
$\gm$ on $\cot$ given by scalar multiplication of the fibers.
The action of the compact subgroup $\so\subs\gm$ is hamiltonian
with respect to the first symplectic form $\omega_1$, and descends
to a circle action on $\M$ with moment map $\Phi:\M\to\R$, which is a Morse-Bott function.
The geometry and topology 
associated with this action will be our main object of study.
In Chapter \ref{analogues} we give a detailed discussion of the construction of $\M$,
along with the action of $\gm$.  In the case where $\Phi$ is proper,
we describe a reducible subvariety $\core\subs\M$ called the {\em core} of $\M$,
onto which $\M$ retracts $\so$-equivariantly.
The core $\core$ has $\X$ 
as one of its components, and if $\M$ is smooth, then $\core$ is equidimensional of dimension
$\dim\X = \half\dim\M$.  In particular, the fundamental classes of the components
of $\core$ provide a natural basis for the top degree cohomology of $\M$.  This fact is
exploited for hypertoric varieties in \cite{HS}, and for quiver varieties in various
papers of Nakajima.  Building on \cite{HP1}, this thesis is the first 
unified treatment of hyperk\"ahler analogues and their cores,
encompassing both hypertoric varieties and quiver varieties.

The geometry of the core of $\M$ will be one of two major themes that we consider.
In the case where $\M$ is a hypertoric variety, each of the components of the core
$\core$ is itself a toric variety (Lemma \ref{htcorecomp}), as first shown in \cite{BD}.
In section \ref{htcore}, we give an explicit description of the action of 
$\gm$ and the gradient flow of $\Phi$ on each piece.  The case of hyperpolygon spaces
is more interesting.  The ordinary polygon space $\X$ is the moduli space of $n$-sided polygons
in $\R^3$ with fixed edge lengths.  In Section \ref{modcore}, we show that the other
core components are smooth, and may themselves be interpreted as moduli spaces
of polygons in $\R^3$ satisfying certain conditions (Theorem \ref{us}).
Thus, for the special case of hyperpolygon spaces, 
we have solved the following general problem:

\begin{problem}
Given any moduli space $\X$ that can be constructed as a K\"ahler reduction
(or GIT quotient) of complex affine space, is it possible
to interpret the core of the hyperk\"ahler analogue
$\X$ as a union of moduli spaces corresponding to other, 
related moduli problems?
\end{problem}

Our second major theme will be the circle-equivariant cohomology ring of $\M$.
In Chapter \ref{ht} we compute the circle-equivariant
cohomology ring of a hypertoric variety, and as an application compute the
$\Zt = Gal(\C/\R)$ equivariant cohomology ring of the complement of a complex
hyperplane arrangement defined over $\R$.  The purpose of Chapter \ref{abelianization} is 
to extend to the hyperk\"ahler setting a theorem of Martin \cite{Ma}, which describes
how to compute the cohomology ring of a K\"ahler quotient $X\mod G$ in terms
of the cohomology ring of the abelian quotient $X\mod T$, where $T\subs G$ is a maximal
torus.  The main technical difficulty arises from the fact that 
Martin's theorem relies heavily on computing integrals, which is not possible
on the noncompact hyperk\"ahler analogues that we have defined.
Our approach is to make use of the localization theorem of \cite{AB,BV}, which
allows us to define an integration theory in the circle-equivariant cohomology
of $\so$-manifolds with compact, oriented fixed point set.  This is perhaps the single
most important reason for considering the circle action on a hyperk\"ahler analogue.
In Section \ref{hp} we combine the results of Chapters \ref{ht} and \ref{abelianization}
to compute the equivariant cohomology ring of a hyperpolygon space, and of each
of its core components.

Most of Chapter \ref{ht} (with the exception of Section \ref{cogs})
appeared first in \cite{HP1}, and Chapter \ref{abelianization}
is a reproduction of \cite{HP}.  Chapter \ref{hpspaces} 
is taken primarily from \cite{HP2}, with the exception of Section \ref{s1},
which comes from \cite{HP}.
\end{chapter}

\begin{chapter}{Hyperk\"ahler analogues}\label{analogues}
Our plan for this chapter is to provide a unified approach to the constructions of
hypertoric varieties and quiver varieties, which are the two major classes of examples
of hyperk\"ahler analogues of familiar K\"ahler varieties that appear in the literature.
In Section \ref{reduction} we give the basic construction of the hyperk\"ahler analogue
$\M$ of a K\"ahler quotient $\X=\Cn\mod G$, and show that $\M$ may be understood
as a partial compactification of the cotangent bundle to $\X$ 
(Proposition \ref{denseopen}).
In Section \ref{circleaction}, we define a natural action of the group $\gm$
on $\M$, which is holomorphic with respect to one of the complex structures.
This action will be our main tool for studying the geometry of $\M$
in future chapters.
Some of this material appeared first in \cite[\S 1]{HP1}.

\begin{section}{Hyperk\"ahler and holomorphic 
symplectic reduction}\label{reduction}
A {\em hyperk\"ahler manifold} is a Riemannian manifold $(M,g)$ along with
three orthogonal, parallel complex structures, $J_1,J_2,J_3$,
satisfying the usual quaternionic relations. 
These three complex structures allow us to define three symplectic forms
$$\omega_1(v,w) = g(J_{1}v, w), \hs\omega_2(v,w) = g(J_{2}v,w), 
\hs\omega_3(v,w) = g(J_{3}v,w),$$ so that $(g, J_i, \omega_i)$ is a K\"ahler
structure on $M$ for $i=1,2,3$.
The complex-valued two-form $\omega_2 + i
\omega_3$ is closed, nondegenerate, and holomorphic
with respect to the complex structure $J_1$.
Any hyperk\"ahler manifold can therefore be considered as a {\em holomorphic
symplectic} manifold with complex structure $J_1$, real symplectic
form $\omega_{\R} := \omega_1$, and holomorphic symplectic form
$\omega_{\C} := \omega_2 + i\omega_3$.
This is the point of view that we will adopt.

We will refer to an action of $G$ on a hyperk\"ahler manifold
$M$ as {\em hyperhamiltonian}
if it is hamiltonian with respect to $\omega_{\R}$ and holomorphic
hamiltonian with respect to $\omega_{\C}$, with $G$-equivariant
moment map
$$\mhk:=\mr\oplus\mc:M\to\gd\oplus\g_{\C}^*.$$
The following theorem describes the {\em hyperk\"ahler quotient}
construction, a quaternionic analogue of the K\"ahler quotient.

\begin{theorem}\label{hklr}{\em\cite{HKLR}}
Let $M$ be a hyperk\"ahler manifold equipped with a hyperhamiltonian action of a
compact Lie group $G$, with moment maps $\mu_1, \mu_2, \mu_3$.
Suppose $\xi = \xi_{\R}\oplus \xi_{\C}$ is a
central regular value of $\mhk$, with $G$ acting freely
on $\mhk^{-1}(\xi)$. Then there is a unique
hyperk\"ahler structure
on the hyperk\"ahler quotient
$\M = M \mmod_{\!\!\xi}G := \mhk^{-1}(\xi)/G$,
with associated symplectic and holomorphic symplectic forms
$\omega^{\xi}_{\R}$ and $\omega^{\xi}_{\C}$, such that
$\omega^{\xi}_{\R}$ and $\omega^{\xi}_{\C}$ pull back to the restrictions
of $\omega_{\R}$ and $\omega_{\C}$ to $\mhk^{-1}(\xi)$.
\end{theorem}

For a general regular value $\xi$, the action of $G$ on
$\mhk^{-1}(\xi)$ will not be free, but only locally free.
To deal with this situation we must introduce the notion
of a hyperk\"ahler orbifold.

An {\em orbifold} is a topological space $M$ locally modeled
on finite quotients of euclidean space.
More precisely, $M$ is a Hausdorff topological space,
equipped with an {\em atlas} $\mathcal U$ {\em of uniformizing charts}.
This consists of a collection of quadruples
$(\tilde U, \Gamma, U, \phi)$, where $\tilde U$ is an open subset
of euclidean space, $\Gamma$ is a finite group acting on $\tilde U$ and
fixing a set of codimension at least $2$, 
$U$ is an open subset of $M$, and $\phi$ is a homeomorphism from
$\tilde U/\Gamma$ to $U$.  The sets $U$ are required to cover $M$,
and the quadruples must satisfy a list of compatibility
conditions, as in \cite{LT}.

Given a point $p\in M$, the {\em orbifold group at $p$} is the isotropy
group $\Gamma_p\subs\Gamma$ of a point 
$\tilde p \in \phi^{-1}(p)\subs\tilde U$
for any quadruple $(\tilde U, \Gamma, U, \phi)$ such that $U$ contains $p$.
The orbifold tangent space $T_pM = T_{\tilde p}\tilde U_{\tilde p}$ 
should be thought of not as a vector space,
but rather as a representation of $\Gamma_p$ (see Proposition \ref{coreprops}).
A differential form on an orbifold may be thought of as a collection
of $\Gamma$-invariant differential forms on the open sets $\tilde U$,
subject to certain compatibility conditions.  We may define riemannian metrics,
complex structures, vector bundles, K\"ahler structures, and hyperk\"ahler
structures on orbifolds in a similar manner.

\begin{example}
Let $Z$ be a smooth manifold, and let $G$ be a compact Lie group
acting locally freely on $Z$.  Then $Z/G$ inherits the structure 
of an orbifold.  The orbifold group of an orbit of a point $z\in Z$ is simply
the stabilizer $G_z\subs G$.
Any $G$-invariant tensor on $Z$ descends to an orbifold
tensor on $Z/G$.  Any $G$-equivariant vector bundle on $Z$ descends
to a vector bundle on $Z/G$.
{\em All of the orbifolds that we consider will be of this form}
(and it is not known whether any other examples exist).
\end{example}

These definitions allow for a straightforward extension of 
Theorem \ref{hklr} to the case where $\xi$ is an arbitrary regular value
of $\mhk$.  This implies, by the moment map condition, that $G$
acts locally freely on $\mhk^{-1}(0)$, and that the quotient
$\mhk^{-1}(0)/G$ inherits the structure of a hyperk\"ahler orbifold.

Orbifolds are in many ways just as nice as manifolds; for example,
it is possible to adapt Morse theory to the orbifold case, as in \cite{LT},
which we will use in the next section.  When we say that a certain
K\"ahler or hyperk\"ahler quotient is an orbifold, we wish to express
the opinion that it is relatively well behaved, rather than
the opinion that it is nasty and singular.  For this reason, we will
use the positively connoted adjective $\Q$-{\em smooth}
to refer to orbifolds.

We now specialize to the case where $M = \cot$, and the action
of $G$ on $\cot$ is induced from a linear action of $G$ on $\Cn$
with moment map $\mu:\Cn \to \gd$.
Choose an identification of $\Hn$ with $\cot$
such that the complex structure $J_1$ on $\Hn$
given by right multiplication by $i$
corresponds to the natural complex structure on $\cot$.
Then $\cot$ inherits a hyperk\"ahler structure.
The real symplectic form $\omega_{\R}$ is given by
adding the pullbacks of the standard forms on $\Cn$ and $(\Cn)^*$,
and the holomorphic symplectic form $\omega_{\C} = d\eta$, where $\eta$ is the canonical
holomorphic 1-form on $\cot$.
 
We note that $G$ acts $\H$-linearly on $\cot\cong\Hn$
(where $n\times n$ matrices
act on the left on the space of column vectors $\Hn$, and scalar multiplication by $\H$
is on the right).
This action is hyperhamiltonian
with moment map $\mhk=\mc\oplus\mr$,
where
$$\mr(z,w) = \mu(z)-\mu(w)\hspace{10pt}\text{and}
\hspace{10pt}\mc(z,w)(v) = w(\hat{v}_z)$$
for $w\in T^*_z\Cn,\hs v\in\g_{\C}$, and $\hat{v}_z$ the element
of $T_z\Cn$ induced by $v$.
Given a central element $\a\in\gd$,
we refer to the hyperk\"ahler quotient $$\M = \cot \mmod_{\!\!(\a,0)} G$$
as the {\em hyperk\"ahler analogue} of the corresponding K\"ahler quotient
$$\X=\Cn\mod_{\!\!\a} G := \mu^{-1}(\a)/G.$$
In future sections we will often fix a parameter $\a$ and drop it from the notation.

At times it will be convenient to think of K\"ahler quotients
in terms of geometric invariant theory, as follows.
Let $G_{\C}$ be an algebraic group acting on an affine variety $V$,
and let $\chi:G_{\C}\to\gm$ be a character of $G_{\C}$.
This defines a lift of the action of $G_{\C}$ to the trivial line
bundle $V\times \C$ by the formula
$$g\cdot(v,z) = (g\cdot v, \chi(g)^{-1}z).$$
The {\em semistable locus} $V^{s}$ with respect to 
$\chi$ is defined to be the set
of points $v\in V$ such that for $z\neq 0$, the closure
of the orbit $G_{\C}(v,z)\subs V\times\C$ is disjoint from the
zero section $V\times\{0\}$ (see \cite{MFK} or \cite{N4}).
The {\em geometric invariant theory} (GIT) quotient $V\mod_{\!\!\chi}G_{\C}$
of $V$ by $G_{\C}$ at $\chi$ is an algebraic variety with underlying
space $V^{ss}\!/\!\!\sim$, where $v\sim w$ 
if and only if the closures $\overline{G_{\C}v}$
and $\overline{G_{\C}w}$ intersect in $V^{ss}$.
The {\em stable locus} $V^{st}$ with respect to $\chi$ is the set of points
$v\in V^{ss}$ such that the $G_{\C}$ orbit through $v$ is closed in
$V^{ss}$.  Clearly the geometric quotient $V^{st}\!/G_{\C}$ is an open set
inside of the categorical quotient $V\mod_{\!\!\chi} G_{\C}$.
The following theorem is due to Kirwan in the projective case \cite{Ki};
our formulation of it is taken from \cite[\S 3]{N4}
and \cite[\S 8]{MFK}.

\begin{theorem}\label{hkhs}
Let $G$ be a compact Lie group acting linearly on a complex vector
space $V$ with moment map $\mu:V\to\gd$.  Let $G_{\C}$ be the complexification
of $G$, with its induced action on $V$.  Let $\chi$ be a character of
$G$, and let $d\chi$ be the associated element of 5 center of $\gd$.
Then $v\in V^{ss}$ if and only if $G_{\C}v\cap\mu^{-1}(d\chi)\neq\emptyset$,
and the inclusion $\mu^{-1}(d\chi)\subs V^{ss}$ induces a homeomorphism
from $V\mod_{\! d\chi} G$ to $V\mod_{\!\chi} G_{\C}$.  
Furthermore, $d\chi$ is a regular
value of $\mu$ if and only if $V^{ss}=V^{st}$.
\end{theorem} 

Given a regular value $\a\in\gd$, Theorem \ref{hkhs} tells
us that we may interpret $V\mod_{\!\!\a} G$ as a GIT quotient
only in the case where $\a$ comes from a character of $G$.
We note, however, that the stability and semistability conditions 
are unchanged when $\chi$ is replaced by a high power of itself,
hence we may apply Theorem \ref{hkhs} whenever some multiple
of $\a$ comes from a character.  Furthermore, the GIT stability condition
is locally constant as a function of $\chi$.  Hence for any central 
regular value $\a\in\gd$,
we may perturb $\a$ to a ``rational'' point, thereby 
interpret $V\mod_{\!\!\a}G$ as a GIT quotient.
Accordingly, we will call an element of $V$ stable with respect
to $\a\in\gd$ if and only if it is stable with respect to $\chi$
for all $\chi$ such that $d\chi$ is close to a multiple of $\a$,
and we will write $V\mod_{\!\!\a}G_{\C} = V^{ss}\!/\!\!\sim$.

We may also use this theorem to formulate the hyperk\"ahler quotient
construction purely in terms of algebraic geometry.  
Theorem \ref{hkhs} says that, for $\a$ a regular value of $\mur$,
$$\cot\mod_{\!\!\a}G\,\cong\,\cot\mod_{\!\!\a}G_{\C}\,\cong\,\(\cot\)^{st}\!/G_{\C}.$$
Since $\muc:\cot\to\gd_{\C}$ is equivariant, we may take its vanishing locus
on both sides of the above equation, and we obtain the identity
$$\cot\mmod_{\!\!(\a,0)}G\,\cong\,\muc^{-1}(0)\mod_{\!\!\a}G_{\C}
\,\cong\,\muc^{-1}(0)^{st}/G_{\C}.$$

The following proposition is proven for the case
where $G$ is a torus in \cite[7.1]{BD}.
 
\begin{proposition}\label{denseopen}
Suppose that $\a$ and $(\a,0)$ are regular values for $\mu$ and $\mhk$,
respectively.
The cotangent bundle $T^*\X$ is isomorphic to an open subset of $\M$,
and is dense if it is nonempty.
\end{proposition}
 
\begin{proof}
Let $Y = \{(z,w)\in\mc^{-1}(0)^{st}\mid
z \in (\Cn)^{st}\}$,
where we ask $z$ to be semistable with respect to $\a$
for the action of
$G_{\C}$ on $\Cn$, so that $\X\cong (\Cn)^{st}/G_{\C}$.
Let $[z]$ denote the element of $\X$ represented by $z$.
The tangent space $T_{[z]}\X$ is equal to the quotient of $T_{z}\Cn$
by the tangent space to the $G_{\C}$ orbit through $z$,
hence
$$T^*_{[z]}\X\cong\{w\in T^*_{z}\Cn\mid w(\hat{v}_z)=0\text{ for all }
v\in\g_{\C}\} = \{w\in(\Cn)^*\mid\mc(z,w)=0\}.$$
Then $$T^*\X \cong \{(z,w)\mid z\in(\Cn)^{st}\text{ and }
\mc(z,w)=0\}/G_{\C} = Y/G_{\C}.$$
By the definition of semistability, $Y$ is an open subset of $\mc^{-1}(0)$,
and is dense if nonempty.  This completes the proof.
\end{proof}

\vspace{-\baselineskip}
\begin{remark}
We may significantly generalize the construction of hyperk\"ahler analogues
as follows.  Replace $\Cn$ by a smooth complex variety $X$, equipped with an
action of an algebraic group $G_{\C}$, an ample line bundle $L$, and a lift
of the action to $L$.  Then the cotangent bundle $T^*X$ is holomorphic symplectic,
and carries a natural holomorphic hamiltonian action of $G_{\C}$, along with a lift
of this action to the pullback of $L$.  We may then define the holomorphic symplectic
analogue of the GIT quotient $\X = X\mod G_{\C}$ to be the GIT quotient of the
zero level of the holomorphic moment map in $T^*X$, where the semistable
sets are defined by the action of $G_{\C}$ on $L$.  Theorem \ref{hkhs} tells us
that this agrees with our construction if $X=\Cn$, and Proposition \ref{denseopen}
generalizes to say that the holomorphic symplectic analogue of $\X$ is a partial
compactification of its cotangent bundle.

The reason for relegating this definition to a remark is that when $X$
is not equal to $\Cn$, its cotangent bundle $T^*X$ may not be the best holomorphic
symplectic manifold with which to replace it.  For example, if $X$ is itself a
K\"ahler quotient, then the holomorphic symplectic analogue of $X$ modulo the trivial
group, in the sense of the previous paragraph, would simply be the cotangent bundle
to $X$.  But this would (usually) not agree with the hyperk\"ahler
analogue of $X$.
\end{remark}

\end{section}
 
\begin{section}{The $\bf{\gm}$ action and the core}\label{circleaction}
Consider the action of $\gm$ on $T^*\Cn$ given by
scalar multiplication on the fibers, that is $\tau\cdot(z,w) = (z,\tau w)$.
The holomorphic moment map $\muc:\cot\to\gdc$ is $\gm$-equivariant
with respect to the scalar action on $\gdc$, 
hence $\gm$ acts on $\muc^{-1}(0)$.  Linearity of the action of $G$
on $\Cn$ implies that the actions of $\GC$ and $\gm$ on $\cot$ commute, 
therefore we obtain a $J_1$-holomorphic action of $\gm$ on 
$\Ma = \muc^{-1}(0)\mod\GC$.
Note that the $\gm$ action does {\em not} preserve the holomorphic symplectic form
or the hyperk\"ahler structure on $\Ma$;
rather we have $\t^*\omega_{\C}=\t\omega_{\C}$ for $\t\in\gm$.

If $\Ma$ is $\Q$-smooth, then the action of the compact subgroup
$\so\subs\gm$ is hamiltonian with respect to the real symplectic
structure $\omega_{\R}$, with moment map
$\Phi[z,w]_{\R} = \half |w|^2$.
This map is an orbifold Morse-Bott function\footnote{For a detailed discussion
of hamiltonian circle actions and Morse theory on orbifolds,
see \cite{LT}.} with image
contained in the non-negative real numbers, and
$\Phi^{-1}(0)=\Xa\subs\Ma$.

\begin{proposition}\label{proper}
If the original moment map $\mu:\Cn\to\gd$ is proper,
then so is $\Phi:\Ma\to\R$.
\end{proposition}

\begin{proof}
We must show that $\Phi^{-1}[0,R]$ is compact
for any $R\in\R$.  We have
$$\Phi^{-1}[0,R] = \{(z,w)\mid \mr(z,w)=\a,\hs
\mc(z,w)=0,\hs\Phi(z,w)\leq R\}\big/G$$ and $G$ is compact,
hence it is sufficient to show that
$\{(z,w)\mid \mr(z,w)=\a,\hs\Phi(z,w)\leq R\}$
is compact.  Since $\mr(z,w)=\mu(z)-\mu(w)$, this set is a closed subset
of $$\mu^{-1}\left\{\a+\mu(w) \,\Bigg{|}\, \half |w|^2\leq R\right\}
\times\left\{w \,\Bigg{|}\, \half |w|^2 \leq R\right\},$$
which is compact by the properness of $\mu$.
\end{proof}

\vspace{-\baselineskip}
\begin{remark}\label{tc}
In the case where $G$ is abelian and $\Xa$ is a nonempty toric
variety, properness of $\mu$ (and therefore of $\Phi$) is equivalent
to compactness of $\Xa$.
\end{remark}

Suppose that $\Ma$ is $\Q$-smooth and $\Phi$ is proper.
We define the {\em core} $\core\subs\Ma$ to be the union
of those $\gm$ orbits whose closures are compact.
Properness of $\Phi$ implies that $\Ld_{\t\to 0}\t\cdot p$ exists for all $p\in \M$,
hence we may write 
$$\core = \{p\in\Ma\mid\Ld_{\t\to\infty}\t\cdot p\text{  exists}\}.$$
For $F$ a connected component of $\Mso=\Macs$, 
let $U_F$ be the closure of the set
of points $p\in \M$ such that $\Ld_{\t\to\infty}\t\cdot p\in F$.
In Morse-theoretic language, 
$U(F)$ is the closure of the unstable orbifold of the critical set $F$.
We may then write $\core$ as a finite union of irreducible, compact
varieties as follows: $$\core = \bigcup_{F\subs\Macs}U_F.$$

\begin{proposition}\label{coreprops}
The core of $\Ma$ has the following properties:
\begin{enumerate}
\item $\core$ is an $\so$-equivariant deformation retract of $\Ma$
\item $U_F$ is isotropic with respect to $\omega_{\C}$
\item If $\Ma$ is smooth at $F$, then $\dim U_F = \half\dim\Ma$.
\end{enumerate}
\end{proposition}

\begin{proof}
Let $f:\Ma\to [0,1]$ be a smooth, $\so$-invariant function with $f^{-1}(0)=\core$.
For all $p\in \M$ and $t\geq 0$, let $\rho_t(p) = e^{f(p)t}\cdot p$.
This defines a smooth family of $\so$-equivariant maps $\rho_t:\Ma\to\Ma$, 
fixing $\core$, with $\rho_0 = \id$.
The limit $\rho_{\infty} = \Ld_{t\to\infty}\rho_t$ is a well-defined
smooth map from $\M$ to $\core$, hence (1) is proved.

Suppose that $\Ma$ is smooth at $F$ and consider a point $p\in F$.  
Since $p$ is a fixed point, $\so$ acts on $T_p\Ma$, 
and we may write $$T_p\Ma = \Od_{s\in\Z}H_s,$$
where $H_s$ is the $s$ weight space for the circle action.
Since $\t^*\omega_{\C}=\t\omega_{\C}$ and $\omega_{\C}$ is a nondegenerate
bilinear form on $T_p\Ma$, $\omega_{\C}$ restricts to a perfect pairing
$H_s\times H_{1-s}\to\C$ for all $s\in\Z$.
In particular, $$T_pU_F = \Od_{s\leq 0}H_s$$ is a maximal isotropic subspace
of $T_p\Ma$, thus proving (2) and (3).

Now suppose that $\Ma$ is only $\Q$-smooth at $F$, and let $\Gap$ be the
orbifold group at $p$.  A circle action on the orbifold tangent
space $T_p\M$ is an action of a group $\Gapt$,
where $\Gapt$ is an extension of $\so$ by $\Gap$.  
Let $\Gapto$ be the connected
component of the identity in $\Gapt$.  Then $\Gapto\,$ is itself isomorphic
to a circle, and maps to the original circle $\so$ with some degree $d\geq 1$.
We now decompose $T_p\Ma = \bigoplus H_s$ into $\Gapto$ weight spaces.
Again $\omega_{\C}$ is nondegenerate on $T_p\Ma$, but now
$\hat{\t}^*\omega_{\C} = \hat{\t}^d\omega_{\C}$ 
for $\hat{\t}\in\Gapto\cong\so$, hence $\omega_{\C}$ restricts to a perfect
pairing $H_s\times H_{d-s}\to\C$ for all $s\in\Z$.
It follows that $T_pU_F = \bigoplus_{s\leq 0}H_s$ is isotropic (though not necessarily
maximally isotropic\footnote{See Example \ref{orbi}.})
with respect to $\omega_{\C}$.  This completes the proof
of (2) in the orbifold case.
\end{proof}

\vspace{-\baselineskip}
\begin{remark}\label{cpt}
Proposition \ref{coreprops} provides a new way to understand Proposition \ref{denseopen}
in the case where $\M$ is $\Q$-smooth and $\Phi$ is proper.
The K\"ahler quotient $\Xa$ is an $\omega_{\C}$-lagrangian
suborbifold of $\Ma$, hence $\omega_{\C}$ identifies the normal
bundle to $\Xa$ in $\Ma$ with the cotangent bundle of $\Xa$.
The Proposition \ref{denseopen} follows from the fact that the normal bundle to
$\Xa$ in $\Ma$ can be identified with the dense open set of points in $\Ma$
that flow down to $\Xa=\Phi^{-1}(0)$ along the gradient of $\Phi$.
This also demonstrates that the $\gm$ action on $\M$ restricts to scalar multiplication
on the fibers of $T^*\Xa$.
\end{remark}

Given a space $M$ equipped with the action of a group $G$, we say that $M$
is {\em equivariantly formal} if the equivariant cohomology ring
$H^*_G(M)$ is a free module over $H^*_G(pt)$.
We end the section with the statement of a fundamental theorem
which we will use repeatedly throughout the paper.  Theorem \ref{formality}
is proven in the compact case in \cite{Ki}, and the proof goes through
in the noncompact case as long as $\Phi$ is proper 
(see, for example, \cite[\S 2.2]{H1} or \cite[4.2]{TW}).

\begin{proposition}\label{formality}
Let $M$ be a symplectic
orbifold with a hamiltonian circle action
such that the moment map $\Phi:M\to\R$ is proper and bounded below,
and has finitely many critical values.
Then $H^*_{\so}(M)$ is a free module over $H^*_{\so}(pt)$.
Moreover, if the action of $\so$ commutes with the action of another
torus $T$, then $H^*_{T\times\so}(M)$ is a free module
over $H^*_{\so}(pt)$.
\end{proposition}
\end{section}
\end{chapter}

\begin{chapter}{Hypertoric varieties}\label{ht}
In this chapter we give a detailed analysis of the construction described
in Chapter \ref{analogues} in the special case where the group $G$ is abelian.
In this case the K\"ahler quotient $\Xa$ is called a {\em toric variety}, and
its hyperk\"ahler analogue $\Ma$ is called a {\em hypertoric variety}.\footnote{Also
known as a {\em toric hyperk\"ahler variety}.}  
These latter spaces were first studied
systematically in \cite{BD}; other references include \cite{K1}, \cite{HS}, and 
\cite{HP1}.  

Just as there is a strong relationship between the geometry of toric
varieties and the combinatorics of polytopes (see, for example, \cite{St}), 
the geometry of hypertoric
varieties interacts with the combinatorics of real hyperplane arrangements.
Hausel and Sturmfels \cite{HS} gave an interpretation of the cohomology
ring of a hypertoric variety as the Stanley-Reisner ring of the matroid
associated to the corresponding arrangement of hyperplanes (Theorem \ref{htm}).  
In Section \ref{htcoh}
we interpret the $\so$-equivariant cohomology ring as an invariant of the oriented
matroid, a richer combinatorial structure that can be associated to a hyperplane
arrangement that is defined over the real numbers (Theorem \ref{htsm}
and Remark \ref{om}).
This result is applied in Section \ref{ossection} to obtain a combinatorial
presentation of the $\Zt$-equivariant cohomology ring of the complement
of an arrangment of complex hyperplanes defined over $\R$, thus enhancing
the classical result of Orlik and Solomon \cite{OS}.
In Section \ref{cogs}, we use the cogenerator approach of \cite{HS}
to explore an aspect of the relationship between the cohomology rings of toric and
hypertoric varieties.
Most of the material presented here, with the exception of 
the entirety of Section \ref{cogs}, has been taken from \cite{HP1} in a modified form.

\begin{section}{Hypertoric varieties and hyperplane arrangements}\label{arrangements}
Let $\tn$ and $\td$ be real vector spaces of dimensions $n$ and $d$, respectively,
with integer lattices $\tnz\subs\tn$ and $\tdz\subs\td$.
Let $\{x_1,\ldots,x_n\}$ be an integer basis for $\tnz$,
and let $\{\bd_1,\ldots,\bd_n\}$ be the dual basis for the dual lattice
$\tndz\subs\tnd$.
Suppose given $n$ nonzero integer vectors
$\{\a_1,\ldots,\a_n\}\subs\tdz$ that span $\td$ over the real numbers.\footnote{In 
each of \cite{BD,K1,HS,HP1}, some additional assumption is placed on the vectors
$\{a_i\}$.  Sometimes they are assumed to be primitive, and sometimes 
they are assumed to generate the lattice $\tdz$ over the integers.  Here
we make neither assumption.}
Define $\pi:\tn\to\td$ by $\pi(x_i)=a_i$, and let $\tk$ be the kernel of $\pi$,
so that we have an exact sequence
$$0 \longrightarrow \tk \stackrel{\i}{\longrightarrow} \tn
\stackrel{\pi}{\longrightarrow} \td\longrightarrow 0.$$
This sequence exponentiates to an exact sequence of abelian groups
$$0 \longrightarrow \Tk \stackrel{\i}{\longrightarrow} \Tn
\stackrel{\pi}{\longrightarrow} \Td\longrightarrow 0,$$
where $$\Tn=\tn/\tnz,\hs\hs\Td=\td/\tdz,\hs\hs\text{and }\Tk=\ker(\pi:\Tn\to\Td).$$
Thus $\Tk$ is a compact abelian group with Lie algebra $\tk$,
which is connected if and only if the vectors $\{a_i\}$ span the lattice
$\tdz$ over the integers.
It is clear that every closed subgroup of $\Tn$ arises in this way.

The restriction to $\Tk$ of the standard action of $\Tn$ on $\cot$
is hyperhamiltonian with hyperk\"ahler moment map
$$\mr\oplus\mc:\cot\to\tkd\oplus(\tk_{\C})^*,$$ where
$$\mur(z,w) = \i^* \left(\half \sum_{i=1}^n (|z_i|^2 - |w_i|^2) \bd_i
\right) \hspace{15pt}
\text{and}\hspace{15pt}\muc(z,w) = \i^*_{\C} \left(\sum_{i=1}^n (z_i w_i)\bd_i\right).$$
Given an element $\a\in\tkd$ with lift $r = (r_1,\ldots,r_n)\in\tnd$,
the K\"ahler quotient
$$\Xa= \Cn\mod\Tk = \mu^{-1}(\a)/\Tk$$
is called a {\em toric variety}, and its hyperk\"ahler
analogue $$\Ma=\cot\mmod\Tk = 
\Big(\mur^{-1}(\a)\cap\muc^{-1}(0)\Big)\Big/\Tk$$
is called a {\em hypertoric variety}.
Both of these spaces admit an effective residual action of the torus $\Td=\Tn/\Tk$
which is hamiltonian in the case of $\Xa$, and hyperhamiltonian in the case
of $\Ma$, with hyperk\"ahler moment map
\begin{eqnarray*}
\barmur[z,w]_{\R}\oplus\barmuc[z,w]_{\R}
&=& \half \sum_{i=1}^n (|z_i|^2 - |w_i|^2 - r_i)\,\bd_i 
\hs\oplus\hs\sum_{i=1}^n (z_i w_i)\,\bd_i\\
&\in&\ker(\i^*)\oplus\ker(\i^*_{\C}) = \tdd\oplus(\td_{\C})^*.
\end{eqnarray*}
In fact, this property may be used to give intrinsic definitions
of toric and hypertoric varieties in certain categories, 
as demonstrated by the following two theorems.

\begin{theorem}{\em\cite{De,LT}}
Any connected symplectic orbifold of real dimension $2d$ which admits
an effective, hamiltonian $\Td$ action with proper moment map is 
$\Td$-equivariantly symplectomorphic to a toric variety.
\end{theorem}

\begin{theorem}\label{bielawski}{\em\cite{Bi}}
Any complete, connected, hyperk\"ahler manifold of real dimension $4d$
which admits an effective, hyperhamiltonian $\Td$ action is $\Td$-equivariantly
diffeomorphic, and Taub-NUT deformation equivalent, to a hypertoric variety.
\end{theorem}

The data that were used to construct $\Xa$ and $\Ma$ consist of a collection of
nonzero vectors $a_i\in\tdz$ and an element $\a\in\tkd$.  It is convenient
to encode in terms of an arrangement of affine hyperplanes in $\tdd$ with some additional
structure.
A {\em rational, weighted, cooriented, affine hyperplane} 
$H\subs\tdd$ is an affine hyperplane
along with a choice of nonzero normal vector $a\in\tdz$.  The word rational
refers to integrality of $a$, and weighted means that $a$ is not required
to be primitive. 
Consider the rational, weighted, cooriented hyperplane
$$H_i =\{v\in\tdd \mid v\cdot a_i + r_i = 0\}$$
with normal vector $a_i\in\tdz$, along with the two half-spaces
\begin{equation}\label{fg}
F_i = \{v\in\tdd \mid v\cdot a_i + r_i \geq 0\}\hspace{10pt}\text{and}
\hspace{10pt}G_i = \{v\in\tdd \mid v\cdot a_i + r_i \leq 0\}.
\end{equation}
Let $$\D = \bigcap_{i=1}^n F_i=\{v \mid v\cdot a_i+r_i\geq 0
\text{  for all  }i\leq n\}$$
be the (possibly empty) weighted polyhedron in $\tdd$ defined by the 
weighted, cooriented, 
affine, hyperplane arrangement $\A = \{H_1,\ldots,H_n\}$.  
Choosing a different lift $r'$ of $\a$ corresponds combinatorially to translating
$\A$ inside of $\tdd$, and geometrically to shifting the K\"ahler and hyperk\"ahler 
moment maps for the residual $\Td$ actions by $r'-r \in\ker\i^*=\tdd$.
Our picture-drawing convention will be to encode the coorientations of the hyperplanes
by shading $\D$, as in Figure \ref{teepee}.  In every example that we consider,
all hyperplanes will have weight $1$; in other words we 
will choose the primitive integer
normal vector inducing the indicated coorientation.

\begin{figure}[h]
\centerline{\epsfig{figure=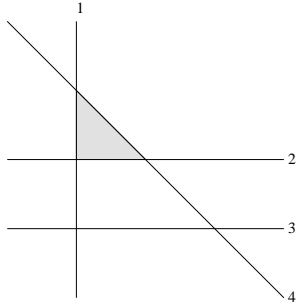, height=4cm}}
\caption{A cooriented arrangement representing a toric variety of complex dimension $2$,
or a hypertoric variety of complex dimension $4$, 
obtained from an action of $T^2$ on $\C^4$.}
\label{teepee}
\end{figure}

We call the arrangement $\A$ {\em simple} if every subset of $m$ hyperplanes
with nonempty intersection intersects in codimension $m$.
We call $\A$ {\em smooth} if every collection of $d$ linearly independent vectors
$\{a_{i_1},\ldots,a_{i_d}\}$ spans $\tdd$.
An element $r\in\tnd$ or $\a\in\tkd$ will be called simple if 
the corresponding arrangement $\A$ is simple.

\begin{theorem}\label{simple}{\em\cite[3.2,3.3]{BD}}
The hypertoric variety $\Ma$ is $\Q$-smooth if and only if $\A$ is simple,
and smooth if and only if $\A$ is smooth.
\end{theorem}


Let us pause to point out the different ways in which 
$\Xa$ and $\Ma$ depend on the
arrangement $\A$.  The toric variety $\Xa$ is in fact determined by the weighted
polyhedron $\D$ \cite{LT},
and is therefore oblivious to any hyperplane $H_i$ such that $\D$ is contained in the interior of $F_i$.
Thus the toric variety corresponding to Figure \ref{teepee} is 
$\C P^2$, the toric variety associated to a triangle.
This is not the case for $\Ma$; we will see, in fact, that the
hypertoric variety of Figure \ref{teepee} is topologically distinct
from the one that we would obtain by deleting the third hyperplane.
For this reason, it is slightly
misleading to call $\Ma$ the hyperk\"ahler analogue of $\Xa$;
more precisely, it is the hyperk\"ahler analogue of
{\em a given presentation} of $\Xa$ as a K\"ahler quotient of a 
complex vector space.

Just as the toric variety $\Xa$ fails to retain all of the data of
the arrangement $\A$, there is some data that goes unnoticed by the
hypertoric variety $\Ma$, as evidenced by the two following results.

\begin{lemma}\label{placement}
The hypertoric variety $\Ma$ is independent, up to $\Td$-equivariant
diffeomorphism,\footnote{Bielawski and Dancer \cite{BD} 
prove a weaker version of this statement,
involving the (nonequivariant) homeomorphism type of $\Ma$.} 
of the choice of a simple element $\a\in\tkd$.
\end{lemma}

\begin{lemma}\label{coorientation}
The hypertoric variety $\Ma$ is independent, up to $\Td$-equivariant isometry,
of the coorientations of the hyperplanes $\{H_i\}$.
\end{lemma}

\begin{proofp}
The set of nonregular values for $\mur\oplus\muc$
has codimension 3 inside of 
$\tkd \oplus (\tk_{\C})^*$ \cite{BD}, hence we may
choose a path connecting the two regular values $(\alpha,0)$
and $(\alpha',0)$ for any simple $\a,\a'\in\tkd$, 
and this path is unique up to homotopy.
Since the moment map $\mur\oplus\muc$
is not proper, we must take some care in showing that
two fibers are diffeomorphic. To this end, we note that the
norm-square function $\psi(z,w) = \|z\|^2 + \|w\|^2$ is $T^n$-invariant
and proper on $\cot$.
Let $(\cot)_{reg}$ denote the open submanifold of $\cot$ consisting of the
preimages of the regular values of $\mur\oplus\muc$.
By a direct computation, it is easy
to see that the kernels of $d\psi$ and $d\mur\oplus d\muc$ intersect
transversely at any point $p \in (\cot)_{reg}$.  Using the
$T^n$-invariant hyperk\"ahler metric on $\cot$, 
we define an Ehresmann connection
on $(\cot)_{reg}$ with respect to $\mur\oplus\muc$
such that the horizontal subspaces are contained in the kernel of $d\psi$.

This connection allows us to lift a path connecting the two
regular values to a horizontal vector field on its preimage in
$(\cot)_{reg}$. Since the horizontal subspaces are tangent to
the kernel of $d\psi$, the flow preserves level sets of $\psi$. Note that the
function $$\mur\oplus\muc\oplus\psi: \cot \to \tkd \oplus (\tk_{\C})^{*} \oplus\R$$
{\em is} proper. By a theorem of Ehresmann \cite[8.12]{BJ}, the
properness of this map implies that the flow of this vector field exists for all time, and
identifies the inverse image of $(\alpha,0)$ with that of
$(\alpha',0)$. Since the metric, $\psi$, and $\mur\oplus\muc$ are
all $T^n$-invariant, the Ehresmann connection is also $T^n$-invariant,
therefore the diffeomorphism identifying the fibers
is $T^n$-equivariant, and
the reduced spaces are $T^d$-equivariantly diffeomorphic.
\end{proofp}

\vspace{-\baselineskip}
\begin{proofc}
It suffices to consider the case when we change the orientation of a
single hyperplane within the arrangement.  Changing the coorientation
of a hyperplane $H_m$ is equivalent to defining a new map
$\pi':\tn\to\td$, with $\pi'(x_i)=a_i$ for $i\neq m$, 
and $\pi'(x_m) = -a_m$.
This map exponentiates to a map $\pi':\Tn\to\Td$, with $\ker(\pi')$
conjugate to $\ker(\pi)$ inside of $\operatorname{GL}_n(\H)$ 
(the group of quaternion-linear
automorphisms of $\cot\cong\Hn$) by the element
$(1,\ldots,1,j,1,\ldots,1)\in \operatorname{GL}_1(\H)^n\subs \operatorname{GL}_n(\H)$,
where the $j$ appears in the $m^{\text{th}}$ slot.
Hence the hyperk\"ahler quotient by $\ker(\pi')$ is isomorphic
to the hyperk\"ahler quotient by $\ker(\pi)$.
\end{proofc}

\vspace{-\baselineskip}
\begin{example}\label{same}
The three cooriented arrangements of Figure 2 all specify
the same hyperk\"ahler variety $\M$ up to equivariant
diffeomorphism.  The first has $\X\cong\widetilde{\C P^2}$ 
(the blow-up of $\C P^2$ at a point),
and the second and the third have $\X\cong\C P^2$.
Note that if we reversed the coorientation of $H_3$ in Figure 2(a)
or 2(c), then we would get a noncompact
$\X\cong\widetilde{\C^2}$.  If we reversed the coorientation of $H_3$ in
Figure 2(b), then $\X$ would be empty, 
but the topology of $\Ma$ would not change.

\begin{figure}[h]
\centerline{\epsfig{figure=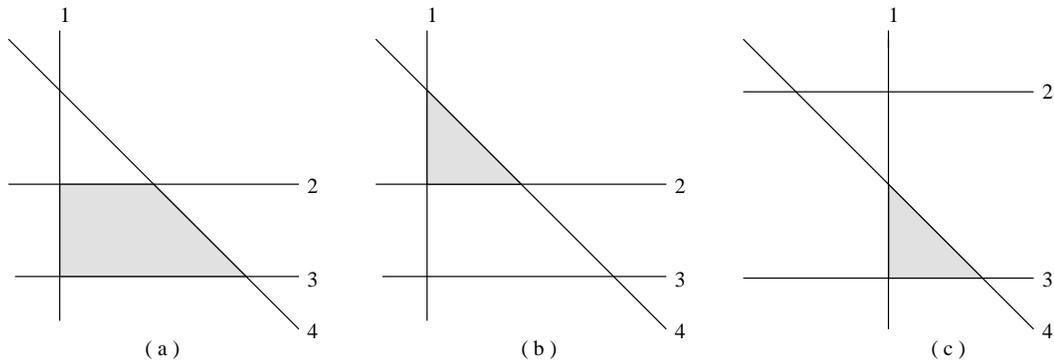}}
\caption{Three arrangements
related by reversing coorientations and translating hyperplanes.}
\label{triplet}
\end{figure}
\end{example}

Our purpose is to study not just the geometry
of $\Ma$, but the geometry of $\Ma$ along with the hamiltonian circle action
defined in Section \ref{circleaction}.
In order to define this action, we used the fact that we were reducing at a regular value
of the form $(\a,0)\in\tdd\oplus(\td_{\C})^*$.
Although the set of regular values of $\mur\oplus\muc$ is simply connected,
the set of regular values of the form $(\a,0)$ is disconnected, therefore the
diffeomorphism of Lemma \ref{placement} is not circle-equivariant.
Furthermore, left multiplication by the diagonal matrix 
$(1,\ldots,1,j,1,\ldots,1) \in 
\operatorname{GL}_n(\H)$ is not an $\so$-equivariant 
automorphism of $\cot\cong\Hn$, 
therefore the topological type of $\Ma$ along with a $\so$ action 
may depend nontrivially
on both the locations and the coorientations of the hyperplanes in $\A$.
Indeed it must, because we can recover $\Xa$ from $\Ma$ by taking the minimum
$\Phi^{-1}(0)$ of the moment map $\Phi:\Ma\to\R$, and we know that $\Xa$ depends in an essential way
on the combinatorial type of the polyhedron $\D$.
In this sense, the structure of a hypertoric variety $\Ma$
along with a hamiltonian circle action may be regarded as 
the universal geometric object associated to $\A$ from which
both $\Ma$ and $\Xa$ can be recovered.
\end{section}

\begin{section}{Geometry of the core}\label{htcore}
In this section we give a combinatorial description of the fixed point set 
$\M^{\gm} = \Ma^{\so}$
and the core $\core$ of a $\Q$-smooth hypertoric variety $\Ma$.  
We will assume that $\Phi$ is proper.  (If $\D$ is nonempty, this is equivalent to asking
that $\D$ be bounded, or that $\Xa$ be compact.)
First, we note that the holomorphic moment map $\barmuc:\Ma\to(\td_{\C})^*$
is $\gm$-equivariant with respect to the scalar action on $(\td_{\C})^*$,
hence both $\Ma^{\gm}$ and $\core$ will be contained in
$$\E=\barmuc^{-1}(0)=\Big\{[z,w]\in 
\Ma\bigmid z_iw_i=0\text{ for all } i\Big\},$$ 
which we call the {\em extended core} of $\Ma$.
It is clear from the defining equations that the restriction of
$\barmur$ from $\E$ to $\tdd$ is surjective.
The extended core naturally breaks into components
$$\EA = \Big\{[z,w]\in \Ma\bigmid w_i=0
\text{ for all $i\in A$ and $z_i=0$ for all $i\in A^c$}\Big\},$$
indexed by subsets $A\subs\otn$.
When $A=\emptyset$, $\EA = \X \subs \Ma$.  More generally,
the variety $\EA\subs \Ma$ is a $d$-dimensional K\"ahler subvariety 
of $\Ma$ with an effective hamiltonian $\Td$-action, and is therefore itself a toric variety.
(It is the K\"ahler quotient by $\Tk$ of an $n$-dimensional 
coordinate subspace of $\cot$,
contained in $\muc^{-1}(0)$.)
The hyperplanes $\{H_i\}$ divide $\tdd$ into a union of closed,
(possibly empty) convex polyhedra $$\Delta_A = 
\bigcap_{i\in A}F_i\,\,\cap\,\,\bigcap_{i\in A^c}G_i.$$

\begin{lemma}\label{sides}
If $w_i=0$, then $\barmur[z,w]_{\R}\in F_i$.
If $z_i=0$, then $\barmur[z,w]_{\R}\in G_i$.
\end{lemma}

\begin{proof}
We have
$$\barmur[z,w]_{\R}\cdot a_i + r_i = \mur(z,w)\cdot x_i
=\half\Big(|z_i|^2 - |w_i|^2\Big),$$
hence the statement follows from Equation~\eqref{fg}.
\end{proof}

\vspace{-\baselineskip}
\begin{lemma}\label{htcorecomp}{\em\cite{BD}}
The core component $\EA$ is isomorphic to the toric variety 
corresponding to the weighted polytope $\D_A$.
\end{lemma}

\begin{proof}
Lemma \ref{sides} tells us that $\barmur(\EA)\subs\Delta_A$,
and surjectivity of $\barmur|_{\E}$ says that this inclusion is an equality.
The lemma then follows from the classification theorems of \cite{De,LT}.
\end{proof}

Although $\gm$ does not act on $\Ma$ as a subtorus of $\Td_{\C}$, we show
below that when restricted to any single
component $\EA$ of the extended core, $\gm$
{\em does} act as a subtorus of $\Td_{\C}$, with the subtorus depending
combinatorially on $A$.  This will allow us to give a combinatorial
analysis of the gradient flow of $\Phi$ on the extended core.

Suppose that $[z,w]\in\EA$.  Then for $\tau\in\gm$,
$$\tau[z,w]=[z,\tau w] = 
[\tau_1 z_1,\ldots \tau_n z_n, \tau_1^{-1}w_1,\ldots \tau_n^{-1}w_n], \text{ where }
\tau_i =
\begin{cases}
\tau^{-1} & \text{if }i\in A,\\
1 & \text{if }i\notin A.
\end{cases}
$$
In other words, the $\gm$ action on $\EA$
is given by the one dimensional subtorus
$(\tau_1, \ldots, \tau_n)$
of the original torus $T^n_{\gm}$,
hence the moment map $\Phi|_{\EA}$ for the action of $\so\subs\gm$ is given by
$$\Phi[z,w] = \left<\mr[z,w],\,\,\sum_{i\in A}a_i\right>.$$
This formula allows us to compute the fixed points of the circle action.
For any subset $B\subs\{1,\ldots,n\}$, let $\EA^B$ be the toric subvariety
of $\EA$ defined by the conditions $z_i=w_i=0$ for all $i\in B$.
Geometrically, $\EA^B$ is defined by the (possibly empty) intersection of the hyperplanes
$\{H_i\mid i\in B\}$ with $\Delta_A$.

\begin{proposition}
The fixed point set of the action of $S^1$ on $\EA$
is the union of those toric subvarieties $\EA^B$ such that
$\sum_{i\in A}a_i\in \td_B:=\operatorname{Span}_{j\in B}a_j$.
\end{proposition}

\begin{proof}
The moment map $\Phi|_{\EA^B}$ will be constant if and only if
$\sum_{i\in A}a_i$ is perpendicular to
the face $\Phi(\EA^B)$,
i.e. if $\sum_{i\in A}a_i$ lies in the kernel of the projection
$\td\onto\td/\td_B$.
\end{proof}

\vspace{-\baselineskip}
\begin{corollary}
Every vertex $v\in\tdd$ of the polyhedral complex $|\A|$ given by our arrangement
is the image of an $S^1$-fixed point in $\M$.
Every component of $\M^{S^1}$ maps to a face of $|\A|$.
\end{corollary}

\begin{proposition}
The core $\core$ of $\Ma$ is equal to the union of those
subvarieties $\EA$ such that $\Delta_A$ is bounded.
\end{proposition}

\begin{proof}
Because $\gm$ acts on $\EA$ as a subtorus of the complex torus $\Td_{\C}$,
the set $$\{p\in\EA\mid\Ld_{\t\to\infty}\t\cdot p\text{ exists}\}$$
is a (possibly reducible) toric subvariety of $\EA$, i.e. a union of subvarieties
of the form $\EA^B$.  Fix a subset $B\subs\otn$.
The variety $\EA^B$ is stable under the $\gm$ action, hence
if $\EA^B$ is compact, then $\EA^B\subs\core$.  On the other hand if $\EA^B$ is noncompact,
then properness of $\Phi$ precludes it from being part of the core,
hence $$\core = \{p\in\E\mid\Phi(p)\text{ lies on a bounded face of }|\A|\}.$$
By \cite[6.7]{HS}, the bounded complex of $|\A|$ has pure dimension $d$,
and is therefore equal to the union of those polyhedra $\Delta_A$ that are bounded.
\end{proof}

\vspace{-\baselineskip}
\begin{corollary}\label{bij}
There is a natural injection from the set of bounded regions
$\{\Delta_A\mid A\in I\}$
to the set of connected components of $\M^{\gm}$.
If $\A$ is smooth, this map is a bijection.
\end{corollary}

\begin{proof}
To each $A\in I$, we associate the fixed subvariety $\EA^B$ corresponding
to the face of $\Delta_A$ on which the linear functional $\sum_{i\in A}a_i$
is minimized, so that $\EA=U(\EA^B)$.
Proposition \ref{coreprops}~(3) tells us that if $\A$ is smooth and 
$F$ is a component
of $\Ma^{\gm}$, then we will have $U(F)=\EA$ for some $A\subs\otn$.
\end{proof}

\vspace{-\baselineskip}
\begin{example}
In Figure 3, representing a reduction of $T^*\C^5$ by $T^3$,
we choose a metric on $(\mathfrak{t}^2)^*$ in order
to draw the linear functional $\sum_{i\in A}a_i$ as a vector
in each region $\Delta_A$.  We see that $\M^{S^1}$ has three components,
one of them $X\cong\widetilde{\C P}^2$, one of them a projective line with another 
$\widetilde{\C P}^2$ as its associated core component, and one of them a point with 
core component $\C P^2$.

\begin{figure}[h]
\centerline{\epsfig{figure=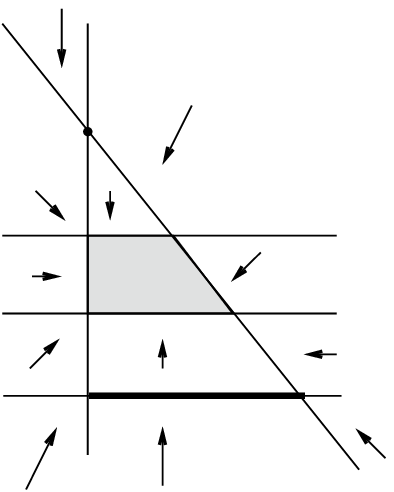}}
\caption{The gradient flow of $\Phi:\Ma\to\R$.}
\end{figure}
\end{example}

\begin{example}\label{orbi}
The hypertoric variety represented by Figure 4 has a fixed point set with four
connected components (three points and a $\C P^2$), but only three components in its core.
This phenomenon can be blamed on the orbifold point $p$ represented by the
intersection of $H_3$ and $H_4$, 
which has only a one-dimensional unstable orbifold (to its northwest).
In other words, this example illustrates the necessity of the smoothness
assumption to obtain a bijection in Corollary \ref{bij}.

\begin{figure}[h]
\centerline{\epsfig{figure=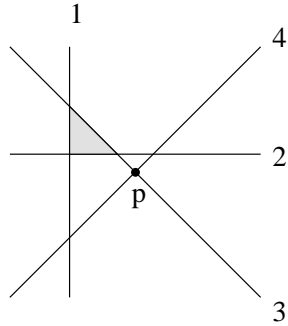}}
\caption{A singular example.}
\end{figure}
\end{example}
\end{section}

\begin{section}{Cohomology rings}\label{htcoh}
In this section we compute the $\so$ and $\Tdso$-equivariant cohomology
rings of of a $\Q$-smooth hypertoric variety $\Ma$, 
extending the computations of the ordinary and $\Td$-equivariant
rings given in \cite{K1} and \cite{HS}.
For the sake of simplicity, and with an eye toward the applications in Chapter
\ref{abelianization}, we will restrict our attention
to the case where $\Phi$ is proper (see Remark \ref{tc}).
This assumption will be necessary for the application of
Proposition \ref{formality} and the proof of Theorem \ref{tara}.

\begin{remark}\label{qorz}
Because we wish to treat the smooth and $\Q$-smooth cases simultaneously,
we will work with cohomology over the rational numbers.  We note, however,
that Konno proves Theorem \ref{htm} over the integers in the smooth case,
and therefore our Theorem \ref{htsm} holds over the integers when $\A$ is smooth.
This fact will be significant in Section \ref{ossection}, when we will need to reduce our
coefficients modulo 2.
\end{remark}

Consider the hyperk\"ahler Kirwan maps
$$\ktd:H^*_{\Tn}(\cot)\to\htm\hspace{15pt}\text{and}
\hspace{15pt}\k:H^*_{\Tk}(\cot)\to H^*(\Ma)$$
induced by the $\Tn$-equivariant inclusion of $\mur^{-1}(\a)\cap\muc^{-1}(0)$
into $\cot$.  Because the vector space $\cot$ is 
$\Tn$-equivariantly contractible, we have
$$H^*_{\Tn}(\cot)=\Sym\tnd\cong\Q[\bd_1,\ldots,\bd_n]$$ and
$$H^*_{\Tk}(\cot)=\Sym\tkd\cong\Q[\bd_1,\ldots,\bd_n]/\ker(\i^*).$$

\begin{theorem}\label{htm}\label{hm}{\em\cite{K1,HS}}
The Kirwan maps $\ktd$ and $\k$ are surjective, and the kernels of both are
generated by the elements
$$\prod_{i\in S}\bd_i\hspace{15pt}\text{  for all $S\subs\otn$ such that  }
\bigcap_{i\in S} H_i = \emptyset.$$
\end{theorem}

\begin{remark}
The kernel of $\ktd$ is precisely
the Stanley-Reisner ideal of
the matroid of linear dependencies among the vectors $\{a_i\}$ \cite{HS}
(see Remark \ref{om}).
\end{remark}

The inclusion of $\mur^{-1}(\a)\cap\muc^{-1}(0)$
into $\cot$ is also $\so$-equivariant, hence we may consider the 
analogous circle-equivariant Kirwan maps
$$\ktdso:H^*_{\Tn\times\so}(\cot)\to\htsm\hspace{15pt}\text{and}
\hspace{15pt}\k:H^*_{\Tk\times\so}(\cot)\to H^*_{\so}(\Ma),$$
where $$H^*_{\Tn\times\so}(\cot)\cong\Q[\bd_1,\ldots,\bd_n,x]$$ and
$$H^*_{\Tk}(\cot)\cong\Q[\bd_1,\ldots,\bd_n,x]/\ker(\i^*).$$
The remainder of this section will be devoted to proving the following theorem.

\begin{theorem}\label{htsm}\label{hsm}
The circle-equivariant Kirwan maps $\ktdso$ and $\kso$ are surjective,
and the kernels of both are generated by the elements
\begin{eqnarray*}
\prod_{i\in S_1}\bd_i\times\prod_{j\in S_2}(x-\bd_j)\hspace{15pt}
&\text{for}&\!\!\text{all disjoint pairs $S_1,S_2\subs\otn$}\\
&&\text{such that }
\bigcap_{i\in S_1} G_i\,\,\cap\,\,\bigcap_{j\in S_2} F_j = \emptyset.
\end{eqnarray*}
\end{theorem}

\begin{remark}\label{om}
For all $i\in\otn$, let $b_i = a_i\oplus 0\in\td\oplus\R$, 
and let $b_0=0\oplus 1$.  The
{\em pointed matroid} associated to $\A$ is a combinatorial object that 
tells us which subsets of $\{b_0,\ldots,b_n\}$ are linearly dependent.
By simplicity of $\A$, this is equivalent to knowing which subsets
$S\subs\otn$ have the property that $$\bigcap_{i\in S} H_i=\emptyset,$$ 
which is in turn
equivalent to knowing the dependence relations among the vectors $\{a_i\}$.
In particular, it does not depend on the relative positions of the hyperplanes,
encoded by the parameter $\a\in\tkd$.

The {\em pointed oriented matroid} associated to $\A$ encodes the data
of which subsets of $\{\pm b_0,\ldots,\pm b_n\}$ are linearly dependent
over the {\em positive} real numbers.  This is equivalent to knowing
which pairs of subsets $S_1,S_2\subs\otn$ have the property that
$$\bigcap_{i\in S_1} G_i\,\,\cap\,\,\bigcap_{j\in S_2} F_j = \emptyset,$$
which does indeed depend on $\a$.
Hence Theorem \ref{htm} shows that $\htm$ is an invariant
of the pointed matroid of $\A$, and Theorem \ref{htsm} demonstrates that $\htsm$
is an invariant of the pointed oriented matroid of $\A$.  
For more on this perspective, see \cite{H3} and \cite{Pr}.
\end{remark}

%
%
Consider the following commuting square of maps, where
$\phi$ and $\psi$ are each given by setting the image of 
$x\in\ H^*_{\so}(pt)\cong\Q[x]$ to zero.

$$\begin{CD}
H^*_{\Tn\times S^1}(\cot) @>\ktdso >> \htsm\\
@V\phi VV @ VV\psi V\\
H^*_{\Tn}(\cot) @>\ktd >> \htm\\
\end{CD}$$

\begin{lemma}\label{surjective}
The equivariant Kirwan map $\ktdso$ is surjective.
\end{lemma}

\begin{proof}
Suppose that $\gamma\in\htsm$ is a homogeneous class of minimal degree
that is {\em not} in the image of $\ktdso$.  By Theorem \ref{htm} $\ktd$ is surjective, 
hence we may choose a class $\eta\in\phi^{-1}\ktd^{-1}\psi(\gamma)$.
Theorem \ref{formality} tells us that the kernel of $\psi$ is generated by $x$,
hence $\ktdso(\eta)-\gamma = x\delta$ for some $\delta\in\htsm$.
Then $\delta$ is a class of lower degree that is not in the image of $\ktdso$.
\end{proof}

\vspace{-\baselineskip}
\begin{lemma}\label{enough}
If $\mathcal{I}\subs\ker\ktdso$ and $\phi(\mathcal{I}) = \ker\ktd$,
then $\mathcal{I}=\ker\ktdso$.
\end{lemma}

\begin{proof}
Suppose that $a\in\ker\ktdso\smallsetminus\mathcal{I}$ 
is a homogeneous class of minimal degree,
and choose $b\in\mathcal{I}$ such that $\phi(a-b)=0$.  Then $a-b = cx$ for some
$c\in H^*_{\Tn\times S^1}(\cot)$.  By Proposition \ref{formality},
$cx\in\ker\ktdso\impl c\in\ker\ktdso$, hence $c\in\ker\ktdso\smallsetminus\mathcal{I}$
is a class of lower degree than $a$.
\end{proof}

\vspace{-\baselineskip}
\begin{proofhtsm}
For any element $$h\in H^2_{T^n\times S^1}(\cot)\cong\Q\{\bd_1,\ldots,\bd_n,x\},$$
let $\tilde{L}_h = \cot\times\C_{h}$
be the $\Tn\times S^1$-equivariant line bundle on $\cot$
with equivariant Euler class $h$.
Let $$L_h = \tilde{L}_h\Big{|}_{\mur^{-1}(\a)\cap\muc^{-1}(0)}\Big/\Tk$$
be the quotient $\Td\times S^1$-equivariant line bundle on $\Ma$.
We will write $\tilde L_i = \tilde L_{\bd_i}$ and
$\tilde K = \tilde L_x$,
with quotients $L_i$ and $K$.
Since the $\Tdso$-equivariant Euler class $e(L_i)$
is the image of $\bd_i$ under the hyperk\"ahler Kirwan map
$H^*_{T^n\times S^1}(\cot)\to \htsm$, we will abuse notation and
denote it by $\bd_i$.  Similarly, we will denote $e(K)$ by $x$.
Lemma \ref{surjective} tells us that
$\htsm$ is generated by $\bd_1,\ldots,\bd_n,x$.

Consider the $\Tn\times S^1$-equivariant section $\tilde s_i$ of $\tilde{L}_i$
given by the function $\tilde s_i(z,w)=z_i$.
This descends to a $\Tdso$-equivariant section $s_i$ of $L_i$
with zero-set $$Z_i
:= \{[z,w]\in \Ma\mid z_i=0\}.$$
Similarly, the function $\tilde t_i(z,w)=w_i$ defines a $\Tdso$-equivariant
section of $L_i^*\otimes K$ with zero set $$W_i := \{[z,w]\in \Ma\mid w_i=0\}.$$
Thus the divisor $Z_i$ represents the cohomology class $\bd_i$,
and $W_i$ represents $x-\bd_i$.
Note, that by Lemma \ref{sides},
we have $\mr(Z_i) \subs G_i$ and $\mr(W_i) \subs F_i$ for all $i\in\otn$.

Let $S_1$ and $S_2$ be a pair of subsets of $\{1,\ldots n\}$ such that
$\big(\cap_{i\in S_1}G_i\big)\cap
\big(\cap_{j\in S_2}F_j\big) = \emptyset,$ and hence
$$\big(\cap_{i\in S_1}Z_i\big)\cap
\big(\cap_{j\in S_2}W_j\big)\subs\mr^{-1}\bigg(
\big(\cap_{i\in S_1}G_i\big)\cap
\big(\cap_{j\in S_2}F_j\big)\bigg) = \emptyset.$$
Now consider the vector bundle
$E=\Od_{i\in S_1}L_i\,\,\oplus\,\,\Od_{j\in S_2}L_j^*\otimes K$
with equivariant Euler class 
$$e(E)=\prod_{i\in S_1}\bd_i\times\prod_{j\in S_2}(x-\bd_j).$$
The section 
$\left(\oplus_{i\in S_1}s_i\right)\oplus\left(\oplus_{i\in S_2}t_i\right)$
is a nonvanishing equivariant
global section of $E$,
hence $e(E)$ is trivial in $\htsm$.
Theorem \ref{htm} and Lemma \ref{enough} 
tell us that we have found all of the relations.
The proofs of the analogous statements for $\hsm$ are identical.
\end{proofhtsm}

How sensitive are the invariants $\htsm$ and $\hsm$?
We can recover $\htm$ and $\hma$
by setting $x$ to zero, hence they are at least as fine
as the ordinary or $\Td$-equivariant cohomology rings.
The ring $\htsm$ does {\em not} depend on coorientations,
for if $\M'$ is related to $\M$ by flipping the coorientation
of the $l^{\text{th}}$ hyperplane $H_k$, then the map taking 
$\bd_l$ to $x-\bd_l$ is an isomorphism
between $\htsm$ and $H^*_{\Td\times S^1}(M')$.\footnote{The 
oriented matroid of a collection of nonzero vectors
in a real vector space does not change when one of the vectors is negated, hence
the independence of $\htsm$ on coorientations can be deduced from Remark \ref{om}.}
The ring {\em does}, however, depend on $\a$, as we see in Example \ref{ts}.

\begin{example}\label{ts}
We compute the equivariant cohomology ring $\htsm$
for the hypertoric varieties $\M_a$, $\M_b$, and $\M_c$
defined by the arrangements in Figure~\ref{triplet}(a), (b), and (c), 
respectively.
Note that each of these arrangements is smooth, hence Theorem \ref{htsm}
holds over the integers, as in Remark \ref{qorz}.
$$H^*_{\Td\times S^1}(\M_a)=\Z[\bd_1,\ldots,\bd_4,x]\big/
\< \bd_2\bd_3, \bd_1(x-\bd_2)\bd_4, \bd_1\bd_3\bd_4 \>,$$
$$H^*_{\Td\times S^1}(\M_b)=\Z[\bd_1,\ldots,\bd_4,x]\big/
\< (x-\bd_2)\bd_3, \bd_1\bd_2\bd_4, \bd_1\bd_3\bd_4 \>,$$
$$H^*_{\Td\times S^1}(\M_c)=\Z[\bd_1,\ldots,\bd_4,x]\big/
\< \bd_2\bd_3, (x-\bd_1)\bd_2(x-\bd_4), \bd_1\bd_3\bd_4 \>.$$
As we have already observed, the rings 
$H^*_{\Td\times S^1}(\M_a)$ and $H^*_{\Td\times S^1}(\M_b)$
are isomorphic by interchanging $\bd_2$ with $x-\bd_2$.
One can check that the annihilator of $\bd_2$ in $H^*_{\Td\times S^1}(\M_a)$
is the principal ideal generated by $\bd_3$, while the ring $H^*_{\Td\times S^1}(\M_c)$
has no degree $2$ element whose annihilator is generated by a single element of degree $2$.
Hence $H^*_{\Td\times S^1}(\M_c)$ is not isomorphic to the other two rings.
\end{example}

The ring $\hsm$, on the other hand, is sensitive to coorientations
as well as the parameter $\a$, as we see in Example \ref{s}.

\begin{example}\label{s}
We now compute the ring $\hsm$ for
$\M_a$, $\M_b$, and $\M_c$ of Figure 2.
Theorem \ref{hsm} tells us that we need only to quotient the ring
$\htsm$ by $\ker(\iota^*)$.
For $\M_a$, the kernel of $\iota_a^*$ is generated by
$\bd_1+\bd_2-\bd_3$ and $\bd_1-\bd_4$, hence we have
\begin{eqnarray*}
H^*_{S^1}(\M_a)&=&\Z[\bd_2,\bd_3,x]\big/
\< \bd_2\bd_3, (\bd_3-\bd_2)^2(x-\bd_2), (\bd_3-\bd_2)^2\bd_3 \>\\
&\cong & \Z[\bd_2,\bd_3,x]\big/\< \bd_2\bd_3, (\bd_3-\bd_2)^2(x-\bd_2), \bd_3^3 \>.\\
\end{eqnarray*}
Since the hyperplanes of 2(c) have the same coorientations as those
of 2(a), we have $\ker\iota_b^*=\ker\iota_a^*$, hence
\begin{eqnarray*}
H^*_{S^1}(\M_c) &=& \Z[\bd_2,\bd_3,x]\big/
\< \bd_2\bd_3, (x-\bd_3+\bd_2)^2\bd_2, (\bd_3-\bd_2)^2\bd_3 \>\\
&\cong & \Z[\bd_2,\bd_3,x]\big/\< \bd_2\bd_3, (x-\bd_3+\bd_2)^2\bd_2, \bd_3^3 \>.\\
\end{eqnarray*}
Finally, since Figure 2(b) is obtained from 2(a) by flipping the coorientation
of $H_2$, we find that $\ker(\iota_b^*)$ is generated by $\bd_1-\bd_2-\bd_3$
and $\bd_1-\bd_4$, therefore
$$H^*_{S^1}(\M_b)=\Z[\bd_2,\bd_3,x]\big/
\< (x-\bd_2)\bd_3, (\bd_2+\bd_3)^2\bd_2, (\bd_2+\bd_3)^2\bd_3 \>.$$
As in Example \ref{ts}, $H^*_{S^1}(\M_a)$ and $H^*_{S^1}(\M_c)$
can be distinguished by the fact that the annihilator
of $\bd_2\in H^*_{S^1}(\M_a)$ is generated by a single element of degree $2$,
and no element of $H^*_{S^1}(\M_c)$ has this property.
On the other hand, $H^*_{S^1}(\M_b)$ is distinguished from
$H^*_{S^1}(\M_a)$ and $H^*_{S^1}(\M_c)$ by the fact that neither $x-\bd_2$ nor $\bd_3$ cubes to zero.
\end{example}

\end{section}

\begin{section}{The equivariant Orlik-Solomon algebra}\label{ossection}
In this section we restrict our attention to smooth arrangements.
When $\A$ is smooth, all of the computations of Section \ref{htcoh}
hold over the integers (see Remark \ref{qorz}).  Since the rings in question
are torsion-free, the presentations are also valid when the coefficients
are taken in the field field $\Zt$.

Let $\M_{\R} \subs \M$ be the {\em real locus}
$\{[z,w]\in \M\mid z,w \text{ real}\}$ of $\M$ 
with respect to the complex structure $J_1$.
The full group $\Td\times S^1$ does not act on $\M_{\R} $,
but the subgroup $\Tdr\times\Zt$ does act,
where $\Tdr :=\Z_2^d\subs \Td$ is the fixed point set
of the involution of $\Td$ given by complex
conjugation.\footnote{It is interesting to note that
the real locus with respect to the complex structure $J_1$
is in fact a complex submanifold with respect to $J_3$,
on which $\Tdr$ acts holomorphically and $\Zt$ acts anti-holomorphically.}
In this section we will study the geometry of the real locus,
focusing in particular on the properties of the $\Zt$ action.
The following theorem is an extension of the results of \cite{BGH} and \cite{Sc}
to the noncompact case of hypertoric varieties.

\begin{theorem}\label{tara}{\em\cite{HH}}
Let $G = \Tdso$ or $\Td$,
and $G_{\R} = \Tdr\times\Zt$ or $\Tdr$, respectively.
Then we have $H^{*}_{G}(\M;\Zt)\cong
H^{*}_{G_{\R}}(\M_{\R} ;\Zt)$
by an isomorphism that halves the grading.\footnote{In particular,
$H^{*}_{G}(\M;\Zt)$ is concentrated in even degree.}
Furthermore, this isomorphism takes the cohomology classes represented
by the $G$-stable submanifolds $Z_i$ and $W_i$ to those represented by 
the $G_{\R}$-stable submanifolds $Z_i\cap \M_{\R}$ and $W_i\cap \M_{\R}$.
\end{theorem}

Consider the restriction $\Psi$ of the hyperk\"ahler moment map
$\mr\oplus\mc$ to $\M_{\R} $.  Since $z$ and $w$ are real
for every $[z,w]\in \M_{\R} $, the map $\mc$
takes values in $\td_{\R}\subs\td_{\C}$, which we will identify
with $i\R^d$, so that $\Psi$ takes values
in $\R^d\oplus i\R^d\cong\C^d$.
Note that $\Psi$ is $\Zt$-equivariant, with $\Zt$ acting on
$\Cn$ by complex conjugation.

\begin{lemma}
The map $\Psi:{\M_{\R}}\to\C^d$ is surjective, 
and the fibers are the orbits of $\Tdr$.
The stabilizer of a point $p\in {\M_{\R}}$ has order
$2^r$, where $r$ is the number of hyperplanes in the complexified
arrangement $\A_{\C}$ containing the point $\Psi(p)$.
\end{lemma}

\begin{proof}
For any point $a+bi\in\C^d$, choose a point $[z,w]\in \M$ such that
$\mr[z,w]=a$ and $\mc[z,w]=b$.  After moving $[z,w]$ by an element of $\Td$
we may assume that $z$ and $w$ are real, hence we may assume
that $[z,w]\in {\M_{\R}}$.
Then $$\Psi^{-1}(a+bi)=\mur^{-1}(a)\cap\muc^{-1}(b)\cap {\M_{\R}}
= \Td[z,w]\cap {\M_{\R}} = \Tdr[z,w].$$
The second statement follows easily from \cite[3.1]{BD}.
\end{proof}

Let $Y\subs {\M_{\R}}$ be the locus of points on which $\Tdr$
acts freely, i.e. the preimage under $\Psi$ of the space
$\MA:=\bomp$.
The inclusion map $Y\hookto {\M_{\R}}$
induces maps backward on cohomology, which we will denote
$$\phi:\hr\to H^*_{\Tdr}(Y;\Zt)\hspace{15pt}\text{and}\hspace{15pt}
\phi_2:\hrs\to H^*_{\Tdr\times\Zt}(Y;\Zt).$$
By Theorem \ref{tara}, we have
$$\hr\cong H^{*}_{\Td}(\M;\Zt)\hspace{15pt}\text{and}\hspace{15pt}
\hrs\cong H^{*}_{\Tdso}(\M;\Zt).$$
Furthermore, since $\Tdr$ acts freely on $Y$ with quotient $\MA$, we have
$$H^*_{\Tdr}(Y;\Zt)\cong H^*(\MA;\Zt)\hspace{15pt}\text{and}\hspace{15pt}
H^*_{\Tdr\times\Zt}(Y;\Zt)\cong H_{\Zt}^*(\MA;\Zt),$$ hence we may write
$$\phi:H^{*}_{\Td}(\M;\Zt)\to H^*(\MA;\Zt)\hspace{15pt}\text{and}\hspace{15pt}
\phi_2:H^{*}_{\Tdso}(\M;\Zt)\to H_{\Zt}^*(\MA;\Zt).$$
The ring $H^*(\MA;\Zt)$ is a classical invariant of the arrangement $\A$ known as
the {\em Orlik-Solomon algebra} (with coefficients in $\Zt$),
and is denoted by $A(\A;\Zt)$ \cite{OT}.
The ring $H^*_{\Zt}(\MA;\Zt)$ was introduced in \cite{HP1} and further studied
in \cite{Pr}.  We call it the {\em equivariant Orlik-Solomon algebra}
and denote it by $A_2(\MA;\Zt)$.

\begin{proposition}\label{ztf}{\em\cite[2.4]{Pr}}
The space $\MA$ is $\Zt$-equivariantly formal, i.e. $A_2(\A;\Zt)$ is a free module
over $\Zt[x]=\hz$.
\end{proposition}

\begin{theorem}\label{surj}
Both $\phi$ and $\phi_2$ are surjective, with kernels
$$\ker\phi = \Big\la\bd_i^2\,\bigmid\, i\in\otn\Big\ra\hspace{10pt}\text{and}
\hspace{10pt}\ker\phi_2 = \Big\la\bd_i\,(x-\bd_i)\,\bigmid\, i\in\otn\Big\ra.$$
\end{theorem}

\begin{proof}
Theorem \ref{tara} tells us that $\phi_2(\bd_i)$ is represented in 
$H^*_{\Tdr\times\Zt}(Y;\Zt)$ by the submanifold $Z_i\cap Y$, and likewise
$\phi_2(x-\bd_i)$ by the submanifold $W_i\cap Y$.  Since $\mur(Z_i\cap W_i)\subs H_i$,
we have $Z_i\cap W_i\cap Y=\emptyset$, hence $\bd_i(x-\bd_i)$ lies in the kernel of
$\phi_2$ (and therefore $\bd_i^2$ lies in the kernel of $\phi$).

By Proposition \ref{ztf} and a pair of formal arguments
identical to those of Lemmas \ref{surjective} and \ref{enough}, it is sufficient
to prove Theorem \ref{surj} only for $\phi$.  
Quotienting $Z_i\cap Y$ by $\Tdr$, we find that $\phi(\bd_i)$ is represented in $A(\A;\Zt)$
by the submanifold
$$\{v\in\MA\mid v\cdot a_i + r_i \in\R^-\}.$$
The standard presentation of $A(\A;\Zt)$ (see, for example, \cite{OT}) says that
these classes generate the ring, and that all relations between them are generated
by the monomials of Theorem \ref{hm} and $\bd_i^2$ for all $i$.
\end{proof}

\vspace{-\baselineskip}
\begin{remark}
Theorems \ref{htsm} and \ref{surj} combine to give a presentation of the equivariant
Orlik-Solomon algebra $A_2(\A;\Zt)$ in the case where $\A$ is rational, simple, and smooth.
This presentation first appeared in \cite{HP1}.  In \cite{Pr}, we generalize this
presentation to arbitrary real hyperplane arrangements, and in fact to arbitrary
pointed oriented matroids.
\end{remark}

\begin{remark}
The ring $A_2(\A;\Zt)$ is a deformation over the affine line $\operatorname{Spec}\Zt[x]$
from the ordinary Orlik-Solomon algebra $A(\A;\Zt)$ to the Varchenko-Gel$'$fand ring
$VG(\A;\Zt)$ of locally constant $\Zt$-valued functions on the real points of $\MA$ \cite{Pr}.
While the rings $A(\A;\Zt)$ and $VG(\A;\Zt)$ depend only on the matroid associated
to $\A$, the deformation $A_2(\A;\Zt)$ depends on the richer structure
of the pointed oriented matroid (see Remark \ref{om}).  
\end{remark}

\begin{example}\label{first}
Consider the arrangements $\A_a$ and $\A_c$ in Figure 2(a) and 2(c).
By Theorem \ref{surj} and Example \ref{ts} we have
$$H^*_{\Zt}(\mathcal{M}(\A_a);\Zt)\cong
\Zt [\bd_1,\ldots,\bd_4,x]\bigg/
\< \begin{array}{c}
\bd_1(x-\bd_1), \bd_2(x-\bd_2), \bd_3(x-\bd_3), \bd_4(x-\bd_4),\\
\bd_2\bd_3, \bd_1(x-\bd_2)\bd_4, \bd_1\bd_3\bd_4 
\end{array}\>$$
and
$$H^*_{\Zt}(\mathcal{M}(\A_c);\Zt)\cong
\Zt [\bd_1,\ldots,\bd_4,x]\bigg/
\< \begin{array}{c}
\bd_1(x-\bd_1), \bd_2(x-\bd_2), \bd_3(x-\bd_3), \bd_4(x-\bd_4),\\
\bd_2\bd_3, (x-\bd_1)\bd_2(x-\bd_4), \bd_1\bd_3\bd_4
\end{array}\>.$$
The map $f:H^*_{\Zt}(\mathcal{M}(\A_a);\Zt)\to H^*_{\Zt}(\mathcal{M}(\A_b);\Zt)$
given by
$$f(\bd_1) = \bd_1+\bd_2, f(\bd_2) = \bd_2+\bd_3+x, f(\bd_3) = \bd_3,
f(\bd_4) = \bd_2+\bd_4,\text{ and }f(x)=x$$
is an isomorphism of graded $\Zt[x]$-algebras, despite the fact that the
oriented matroids of the two arrangements differ.\footnote{We thank
Graham Denham for finding this isomorphism.}
\end{example}

\begin{example}\label{second}
Now consider the arrangements $\A_a'$ and $\A_c'$ obtained from $\A_a$ and $\A_c$
by adding a vertical line on the far left, as shown in Figure \ref{vertical}.
\begin{figure}[h]
\begin{center}
\psfrag{Ha'}{$\A_a'$}
\psfrag{Hc'}{$\A_c'$}
\includegraphics{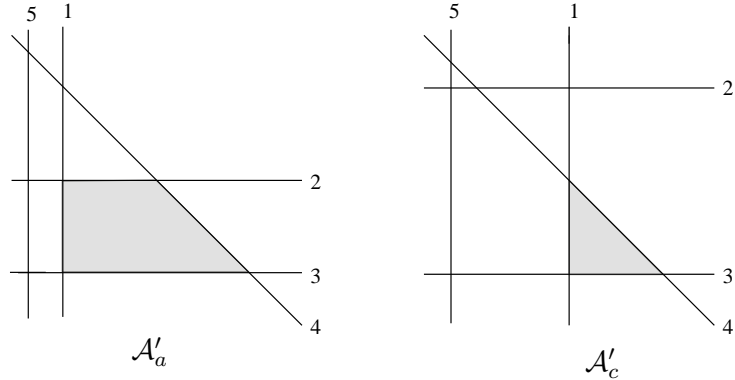}
\caption{Adding a vertical line to $\A_a$ and $\A_c$.}\label{vertical}
\end{center}
\end{figure}
Again by Theorem \ref{surj}, we have
$$H^*_{\Zt}(\mathcal{M}(\A_a');\Zt)\cong
\Zt [\bd_1,\ldots,\bd_4,x]\bigg/
\< \begin{array}{c}
\bd_1(x-\bd_1), \bd_2(x-\bd_2), \bd_3(x-\bd_3), \bd_4(x-\bd_4),\\ \bd_5(x-\bd_5), 
\bd_2\bd_3, (x-\bd_1)\bd_5, \bd_1(x-\bd_2)\bd_4,\\ \bd_1\bd_3\bd_4, (x-\bd_2)\bd_4\bd_5, \bd_3\bd_4\bd_5 
\end{array}\>$$
and
$$H^*_{\Zt}(\mathcal{M}(\A_c');\Zt)\cong
\Zt [\bd_1,\ldots,\bd_4,x]\bigg/
\< \begin{array}{c}
\bd_1(x-\bd_1), \bd_2(x-\bd_2), \bd_3(x-\bd_3), \bd_4(x-\bd_4),\\ \bd_5(x-\bd_5),
\bd_2\bd_3, (x-\bd_1)\bd_5, (x-\bd_1)\bd_2(x-\bd_4),\\ \bd_1\bd_3\bd_4, (x-\bd_2)\bd_4\bd_5, \bd_3\bd_4\bd_5 
\end{array}\>.$$
We have used Macaulay 2 \cite{M2} to check that the annihilator of the element
$\bd_2\in H^*_{\Zt}(\mathcal{M}(\A_a');\Zt)$ is generated by two linear elements
(namely $\bd_3$ and $x-\bd_2$) and nothing else, while there is no element
of $H^*_{\Zt}(\mathcal{M}(\A_c');\Zt)$ with this property.
Hence the two rings are not isomorphic.
\end{example}
\end{section}

\begin{section}{Cogenerators}\label{cogs}
Consider the K\"ahler Kirwan map
$$\ka:\Sym\tkd\cong H^*_{\Tk}(\C^n)\to H^*(\Xa_{\a})$$ 
induced by the $\Tk$-equivariant
inclusion of $\mu^{-1}(\a)$ into $\Cn$.
In this section we would like to consider simultaneously the Kirwan
maps $\ka$ for many different values of $\a$, so almost all of the
notation that we use will have a subscript or superscript indicating
the parameter $\a\in\tkd$ or a lift $r\in\tnd$.
An exception to this rule will be the 
hyperk\"ahler Kirwan map $$\k:\Sym\tkd\to H^*(\Ma_{\a}),$$
which, by Lemma \ref{placement} or Theorem \ref{hm}, 
is independent of our choice of simple $\a\in\tkd$.
The main result of this section is the following.

\begin{theorem}\label{int}
The kernel of the hyperk\"ahler Kirwan map $\k$ is equal to the intersection over all simple
$\a$ of the kernels of the K\"ahler Kirwan maps $\ka$.
\end{theorem}

\begin{remark}
Konno \cite[7.6]{K2} proves an analogous theorem about the kernels
of the Kirwan maps to the cohomology rings of polygon and
hyperpolygon spaces.  We may therefore conjecture a generalization
of Theorem  \ref{int} in which $\Tk$ is replaced by an arbitrary compact
group $G$.  Note that our proof of Theorem \ref{int}
depends strongly on the combinatorics associated to hypertoric varieties.
\end{remark}

We approach Theorem \ref{int} by describing the kernels of $\k$ and $\ka$ not
in terms of generators, but rather in terms of cogenerators.
Given an ideal $\mathcal{I}\subs\Sym\tkd$, a set of {\em cogenerators} for $\mathcal{I}$
is a collection of polynomials $\{f_i\}\subs\Sym\tk$ such that
$$\mathcal{I}=\{\bd\in\Sym\tkd\mid\bd\cdot f_i = 0\text{ for all }i\}.$$

The volume function $\Vol\Dr$ is locally polynomial in $r$.
More precisely, for every simple $r\in\tnd$, there exists a degree $d$ polynomial
$\Pr\in\Sym^d\tn$ such that for every simple $s\in\tnd$ lying in the
same connected component of the set of simple elements as $r$,
we have $$\Vol\D^{s}=\Pr\!(s).$$  We will refer to $\Pr$ as the
{\em volume polynomial} of $\Dr$.  The fact that the volume of a polytope
is translation invariant tells us that $\Pr$ lies in the image of the inclusion
$\i:\Sym^d\tk\hookto\Sym^d\tn$.

\begin{theorem}\label{toric}{\em\cite{GS,KP}}
Let $r\in\tnd$ be a simple element with $\i^*(r)=\a$.  Then
$$\ker\ka = \Ann\{\i^{-1}\Pr\} = 
\big\{\bd\in\Sym\tkd \mid \bd\cdot(\i^{-1}\Pr) = 0\big\}.$$
\end{theorem}

A similar description of the cohomology ring of a hypertoric variety
is given in \cite{HS}.  For any subset
$A\subs\otn$, consider the polyhedron $\DrA$ introduced in Section \ref{htcore}.
If $\DrA$ is nonempty, then it is bounded if and only if the vectors
$\{\eps_1(A)a_1,\ldots,\eps_n(A)a_n\}$ span $\td$ over the non-negative
real numbers,
where $\eps_i(A) = (-1)^{|A\cap\{i\}|}$.
We call such an $A$ {\em admissible}.  For all admissible $A$,
there exists a degree $d$ polynomial
$\PrA\in\Sym^d\tn$ such that for every simple $s\in\tnd$
lying in the same connected component of the set of simple elements as $r$,
we have $$\Vol\D^{s}_A=\PrA(s).$$
Once again, the translation invariance of volume implies that
$\PrA$ lies in the image of the inclusion
$\i:\Sym^d\tk\hookto\Sym^d\tn$.
Consider the linear span
$$U^r = \Q\big\{\PrA \mid A\text{ admissible}\big\}.$$

\begin{theorem}\label{hypertoric}{\em\cite{HS}}
Let $r\in\tnd$ be a simple element with $\i^*(r)=\a$.  Then
$$\ker\k = \Ann \,\i^{-1}U^r = \big\{\bd\in\Sym\tkd \mid 
\bd\cdot \i^{-1}P = 0\text{ for all }P\in U^r\big\}.$$
\end{theorem}

\begin{remark}\label{wow}
It is clear from the statement of Theorem \ref{hypertoric}
that $H^*(\Ma)$ does not depend on the coorientations of the hyperplanes
$\{H_1^r,\ldots,H_n^r\}$, as has been observed in Lemma \ref{coorientation}
and throughout Section \ref{htcoh}.  Indeed,
the polynomials $\PrA$ for $A$ admissible
are simply the volume polynomials of the maximal regions of $|\A|$.
What is not clear from this presentation is the independence of $H^*(\Ma)$
on $\a$.  In other words, it is a nontrivial
fact that the vector space $U^r$ is independent of the parameter $r\in\tnd$.
\end{remark}

\begin{proofint}
The statement of Theorem \ref{int} equates the kernel of $\k$, 
which is cogenerated by the vector
space $\i^{-1}U^r$, with the intersection of the kernels of the maps $\ka$,
each of which is cogenerated by the element $\i^{-1}P^r$.  Intersection of ideals
corresponds to linear span on the level of cogenerators, hence we have
$$\bigcap_{\a\text{ simple}}\ker\ka = \Ann \,\i^{-1}V,\hspace{15pt}\text{where}
\hspace{15pt}V = \Q\big\{\Pr\mid r\text{ simple}\}.$$
Our plan is to show that $V = U^r$ for any simple $r$.
Recall that the assignment of $\Pr$ to $r$ is locally constant on the set
of simple elements of $\tnd$, hence $V$ is finite-dimensional.
Since $\Pr\in U^r$ and $U^r$ does not depend on $r$ (see Remark \ref{wow}),
it is clear that $V\subs U^r$.  Thus to prove Theorem \ref{int}, it will suffice
to prove the opposite inclusion, as stated below.

\begin{proposition}\label{comb}
We have $\PrA\in V$ for every admissible $A\subs\otn$.
\end{proposition}

Let $\F$ be the infinite dimensional vector space consisting of all real-valued
functions on $\tdd$, and let $\Fbd$ be the subspace consisting of functions
with bounded support.  For all subsets $A\subs\otn$, let
$$\WA = \Q\big\{\odra\mid r\text{ simple}\big\}$$
be the subspace of $\F$ consisting of finite linear combinations
of characteristic functions of polyhedra $\DrA$, and let
$$\WAbd = \WA\cap\Fbd.$$
Note that $\WAbd=\WA$ if and only if $A$ is admissible.

\begin{lemma}\label{key}
For all $A,A'\subs\otn$, $\WAbd = W_{A'}^{bd}$.
\end{lemma}

\begin{proof}
We may immediately reduce to the case where $A' = A\cup\{j\}$.
Fix a simple $r\in\tnd$.
Let $\rt\in\tnd$ be another simple element obtained from $r$ by putting
$\rt_i=r_i$ for all $i\neq j$, and $\rt_j = N$ for some $N\gg 0$.
Then $\D_A^r\subs\D_A^{\rt}$, and 
\begin{eqnarray*}
\D_{A}^{\rt} \smallsetminus \D_{A}^{r} &=& 
\big\{v\in\tdd\mid \eps_i(A')(v\cdot a_i+r_i)\geq 0\text{  for all  }i\leq n
\text{ and } v\cdot a_j + N\geq 0\}\\
&=& \D_{A'}^{r} \cap G_j^{N}.
\end{eqnarray*}
Suppose that $f\in\Fbd$ can be written as a linear combination
of functions of the form $\mathbf{1}_{\D_{A'}^{r}}$.  Choosing $N$ large enough
that the support of $f$ is contained in $G_j^{N}$,
the above computation shows that $f$ can be written
as a linear combination of functions of the form $\odra$,
hence $W_{A'}^{bd}\subs W_{A}^{bd}$.  The reverse inclusion is obtained
by an identical argument.
\end{proof}

\vspace{-\baselineskip}
\begin{example}\label{pictures}
Suppose that we want to write the characteristic function for the
upper triangle $\Delta_{\{1,4\}}$ 
in Figure~\ref{triplet}(c) as linear combination
of characteristic functions of the shaded regions obtained by translating
the hyperplanes in any possible way.
Since $\{1,4\}$ has two elements, the procedure described
in Lemma \ref{key} must be iterated twice, and the result will have a total
of $2^2=4$ terms, as illustrated in Figure \ref{char}.
The first iteration exhibits $\mathbf{1}_{\Delta_{\{1,4\}}}$ as an element
of $W_{\{4\}}^{bd}$ by expressing it as the difference of the characteristic
functions of two (unbounded) regions.
With the second iteration, we attempt to express each of these two characteristic functions
as elements of $W_{\{1,4\}}^{bd} = W_{\{1,4\}}$.
This attempt must fail, because each of the two functions that we try
to express has unbounded support.  But the failures cancel out, and we succeed
in expressing the {\em difference} as an element of $W_{\{1,4\}}$.
\begin{figure}[h]
\begin{center}
\includegraphics[height=140mm]{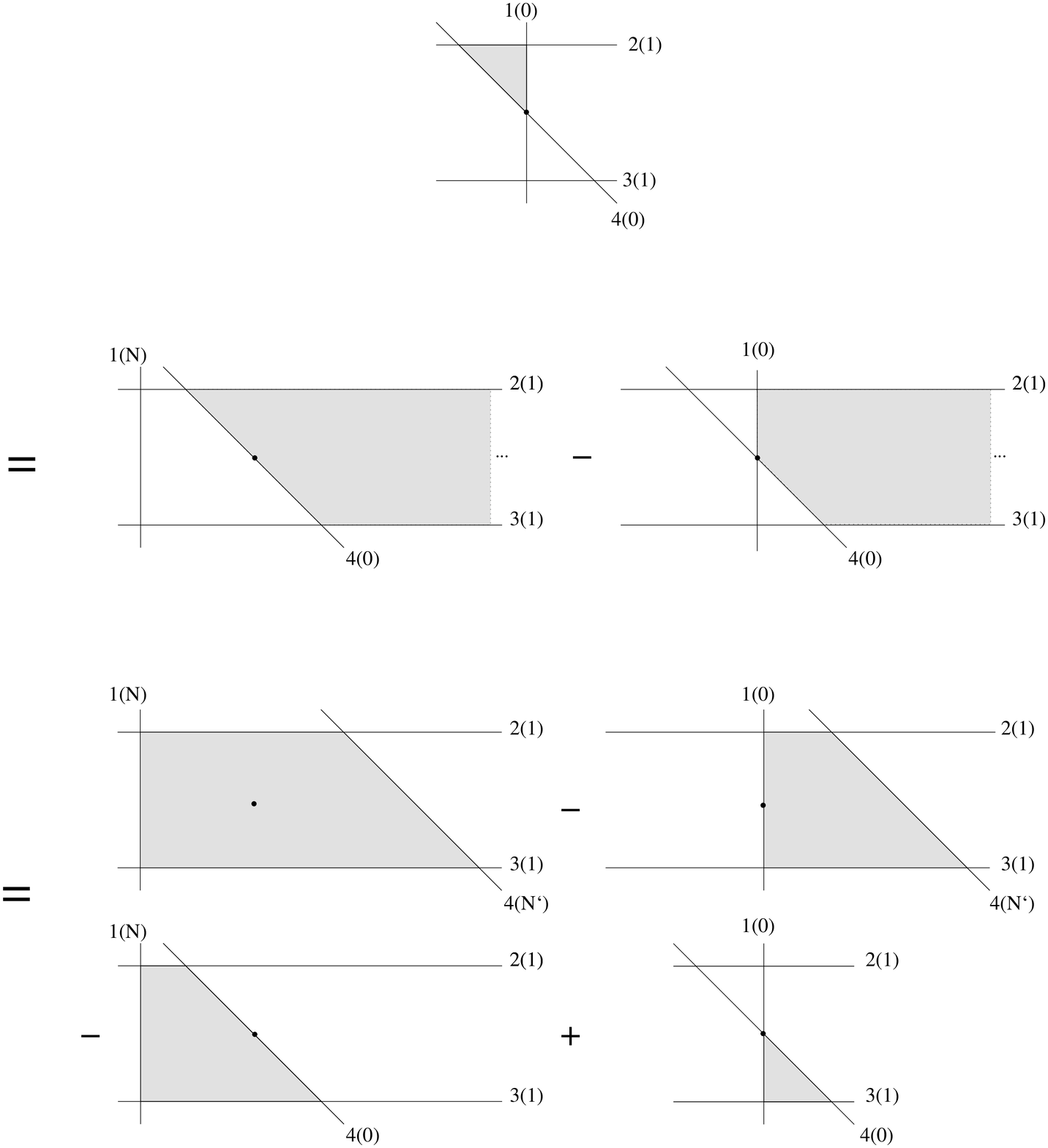}
\caption{An equation of characteristic functions.
We write two numbers next to each hyperplane:  the first is the 
index $i\in\{1,\ldots,4\}$, and the second is the parameter $r_i$
specifying the distance from the origin (denoted by a black dot) to $H_i$.
The two iterations of Lemma \ref{key} have produced two undetermined
large numbers, which we call $N$ and $N'$.}\label{char}
\end{center}\end{figure}
\end{example}

\begin{proofcomb}
By Lemma \ref{key}, we may write
$$\odra = \sum_{j=1}^m\eta_j \mathbf{1}_{\D^{r(j)}}$$
for any simple $r$ and admissible $A$, where $\eta_j\in\Z$ and
$r(j)$ is a simple element of $\tnd$ for all $j\leq m$.
Taking volumes of both sides of the equation, we have
\begin{equation}\label{polys}
\PrA(r) = \sum_{j=1}^m\eta_j P^{r(j)}\big(r(j)\big).
\end{equation}
Furthermore, we observe from the proof of Lemma \ref{key}
that for all $j\leq m$ and all $i\leq n$, the $i^{\text{th}}$
coordinate $r_i(j)$ of $r(j)$ is either equal to $r_i$, or to
some large number number $N\gg 0$.  The Equation~\eqref{polys}
still holds if we wiggle these large numbers a little bit,
hence the polynomial $P^{r(j)}$ must be independent of the variable
$r_i(j)$ whenever $r_i(j)\neq r_i$.  Thus we may substitute
$r$ for each $r(j)$ on the right-hand side, and we obtain the equation
\begin{equation*}
\PrA(r) = \sum_{j=1}^m\eta_j P^{r(j)}(r).
\end{equation*}
This equation clearly holds in a neighborhood of $r$,
hence we obtain an equation of polynomials
$$\PrA = \sum_{j=1}^m\eta_j P^{r(j)}.$$
This completes the proof of
Proposition \ref{comb}, and therefore also of Theorem \ref{int}.
\end{proofcomb}
\end{proofint}

\vspace{-\baselineskip}
\begin{example}
Let's see what happens when we take volume polynomials in the equation
of Figure \ref{char}.  The two polytopes on the top line have different
volumes, but the same volume polynomial, hence these two terms cancel.
We are left with the equation
$$P_{\{1,4\}}^{(0,1,1,0)} = P^{(0,1,1,0)} - P^{(N,1,1,0)},$$
which translates as
$$\half\(-x_1+x_2-x_4\)^2 = \half\(x_1+x_3+x_4\)^2 - 
\(x_2+x_3\)\(x_1+x_4+\half x_3 - \half x_2\).$$
\end{example}
\end{section}
\end{chapter}

\begin{chapter}{Abelianization}\label{abelianization}
Let $X$ be a symplectic manifold equipped with a hamiltonian
action of a compact Lie group $G$.  Let $T\subs G$
be a maximal torus, let $\Delta\subset\tdu$ 
be the set of roots\footnote{Not to be confused with the polyhedron $\Delta$
of Chapter \ref{ht}.} of $G$,
and let $W=N(T)/T$ be the Weyl group of $G$.
If $\mu:X\to\gd$ is a moment map for the action of $G$, then $pr\circ\mu:X\to\tdu$
is a moment map for the action of $T$, where $pr:\gd\to\tdu$ is the standard projection.
Suppose that $0\in\gd$ and $0\in\tdu$ are regular values for the two moment maps.
If the symplectic quotients
$$X\mod G = \mu^{-1}(0)/G\hspace{15pt}\text{and}
\hspace{15pt}X\mod T=(pr\circ\mu)^{-1}(0)/T$$ 
are both compact, then Martin's theorem \cite[Theorem A]{Ma}
relates the cohomology of $X\mod G$ to the cohomology of $X\mod T$.
Specifically, it says that $$H^*(X\mod G)\cong\frac{H^*(X\mod T)^W}{ann(e_0)},$$
where $$e_0 = \prod_{\a\in\Delta}\a\in\left(\operatorname{Sym}\tdu\right)^W
\cong H^*_T(pt)^W,$$ which acts naturally on
$H^*(X\mod T)^W \cong H^*_T(\mu_T^{-1}(0))^W$.
In the case where $X$ is a complex vector space and $G$ acts linearly
on $X$, a similar result was obtained by Ellingsrud and Str\o mme \cite{ES}
using different techniques.

Our goal is to state and prove an analogue of this theorem for
hyperk\"ahler quotients.  There are two main obstacles to this goal.
First, hyperk\"ahler quotients are rarely compact.
The assumption of compactness in Martin's theorem is crucial
because his proof involves integration.
Our answer to this problem is to work
with equivariant cohomology
of {\em circle compact} manifolds, by which we mean oriented manifolds
with an action of $\so$ such that the fixed point set is oriented and compact.
Using the localization theorem of Atiyah-Bott \cite{AB} and Berline-Vergne \cite{BV},
as motivation, we show that integration in rationalized $\so$-equivariant cohomology
of circle compact manifolds can be defined in terms of integration
on their fixed point sets.  Section~\ref{thepush} is devoted to
making this statement precise by defining a well-behaved
push forward in the rationalized $\so$-equivariant cohomology
of circle compact manifolds.

The second obstacle is that Martin's result uses surjectivity
of the K\"ahler Kirwan map from $H^*_G(X)$ to $H^*(X\mod G)$ \cite{Ki}.  
The analogous map for circle compact
hyperk\"ahler quotients is conjecturally surjective, but
only a few special cases are known (see Theorems \ref{hm} and \ref{hpsurj}, and 
Remarks \ref{generation} and \ref{horn}).
Our approach is to assume that the rationalized Kirwan map is
surjective, which is equivalent to saying that the cokernel
of the non-rationalized Kirwan map
$$\K_G:H^*_{\so\times G}(M)\to\hmg$$
is torsion as a module over $\hp$.  This is a weaker assumption
than surjectivity of $\K_G$; in particular, we show in Section
\ref{quiver} that this assumption holds for quiver varieties,
as a consequence of the work of Nakajima.

Under this assumption,
Theorem~\ref{main} computes the rationalized equivariant cohomology
of $M\mmod G$ in terms of that of $M\mmod T$.
We show that, at regular values of the hyperk\"ahler moment maps,
$$\hhso(M\mmod G)\cong\frac{\hhso(M\mmod T)^W}{ann(e)},$$
where $\hhso$ denotes rationalized equivariant cohomology
(see Definition \ref{rat}), and
$$e = \Pd_{\a\in\Delta}\a(x-\a) \in (\operatorname{Sym}\tdu)^W\otimes\Q[x]
\cong H^*_{\so\times T}(pt)^W\subs\widehat{H}^*_{\so\times T}(pt)^W.$$
Theorem~\ref{ordinary}
describes the image of the non-rationalized Kirwan map
in a similar way:
$$\hso(M\mmod G)\supseteq\Im(\K_G)\cong\frac{\imkt^W}{ann(e)},$$
where $\K_T:H^*_{\so\times T}(M)\to \hso(M\mmod T)$ is the Kirwan
map for the abelian quotient.
In many situations, such as when $M=\cot$,
$\K_T$ is known to be surjective (Theorem \ref{hm}).

This Chapter is a reproduction of \cite[\S 1-3]{HP}.

\begin{section}{Integration}\label{thepush}
The localization theorem of Atiyah-Bott \cite{AB} and Berline-Vergne \cite{BV} says that
given a manifold $M$ with a circle action, the restriction map from
the circle-equivariant cohomology of $M$ to the circle-equivariant
cohomology of the fixed point set $F$ is an isomorphism modulo torsion.
In particular, integrals on $M$ can be computed in terms of integrals on $F$.
If $F$ is compact, it is possible to use the Atiyah-Bott-Berline-Vergne formula
to {\em define} integrals on $M$.

We will work in the category of {\em circle compact}
manifolds, by which we mean oriented $\so$-manifolds with compact
and oriented fixed point sets.
Maps between circle compact
manifolds are required to be equivariant.

\begin{definition}\label{rat}
Let $\hhp = \Q(x)$, the rational function field of $\hp\cong \Q[x]$.
For a circle compact manifold $M$, let $\hhm = \hm\otimes\hhp$,
where the tensor product is taken over the ring $\hp$. We call $\hhm$
the {\em rationalized} $\so$-equivariant cohomology of $M$.
Note that because $\operatorname{deg}(x)=2$, $\hhm$ is supergraded,
and supercommutative with respect to this supergrading.
\end{definition}

An immediate consequence of \cite{AB} is that restriction gives an isomorphism
\begin{equation}\label{isomorphic}\hhm\cong \hhso(F)\cong H^*(F)\otimes_{\Q}\hhp,\end{equation}
where $F=M^{S^1}$ denotes
the compact fixed point set of $M$.  In particular
$\hhm$ is a finite dimensional vector space over $\hhp$,
and trivial if and only if $F$ is empty.

Let $i:N\hookto M$ be a closed embedding.
There is a standard notion of proper pushforward
$$i_*:\hso(N)\to\hm$$ given by the formula $i_*=r\circ\Phi$,
where $r:\hso(M,M\setminus N)\to\hso(M)$ is the restriction map,
and $\Phi:\hso(N)\to\hso(M,M\setminus N)$ is the Thom isomorphism.
We will also denote the induced map $\hhso(N)\to\hhm$ by $i_*$.
Geometrically, $i_*$ can be understood as the inclusion
of cycles in Borel-Moore homology.

This map satisfies two important formal properties
\cite{AB}:
\begin{equation}\label{functor}
\text{Functoriality:   } (i\circ j)_* = i_*\circ j_*
\end{equation}
\begin{equation}\label{module}
\text{Module homomorphism:   } i_*(\ga\cdot i^*\a)
=i_\ga\cdot\a
\text{   for all   }\a\in\hhm, \ga\in\hhn.
\end{equation}
We will denote the Euler
class $i^*i_*(1)\in\hhso(N)$ by $e(N)$.
If a class $\ga\in\hhn$ is in the image of $i^*$,
then property (\ref{module}) tells us that
$i^*i_*\ga = e(N)\ga$.  Since the pushforward construction is
local in a neighborhood
of $N$ in $M$, we may assume that $i^*$ is surjective, hence this identity
holds for all $\ga\in\hhn$.

Let $F= M^{\so}$ be the fixed point set of $M$.
Since $M$ and $F$ are each oriented, so is the normal bundle
to $F$ inside of $M$.  The following result is
standard, see e.g. \cite{Ki}.

\begin{lemma}
The Euler class $e(F)\in\hhso(F)$ of the normal bundle to $F$ in $M$
is invertible.
\end{lemma}

\begin{proof}
Let $\{F_1,\ldots,F_d\}$ be the connected components of $F$.
Since $\hhso(F)\cong\bigoplus\hhso(F_i)$ and $e(F)=\oplus e(F_i)$,
our statement is equivalent to showing that $e(F_i)$ is invertible for all $i$.
Since $\so$ acts trivially on $F_i$, $\hhso(F_i)\cong H^*(F_i)\otimes_{\Q}\hhp$.
We have $e(F_i) = 1\otimes ax^k + nil$, where
$k= \operatorname{codim}(F_i)$, $a$ is the product of the weights of the $S^1$
action on any fiber of the normal bundle, and $nil$ consists of terms
of positive degree in $H^*(F_i)$.  Since $F_i$ is a component of the fixed
point set, $\so$ acts freely on the complement of the zero section
of the normal bundle, therefore $a\neq 0$.  Since $ax^k$ is invertible and
$nil$ is nilpotent, we are done.
\end{proof}

\vspace{-\baselineskip}
\begin{definition}
For $\a\in\hhm$, let $$\IM\a=\IF\frac{\a|_F}{e(F)}\in\hhp.$$
\end{definition}

Note that this definition does not depend on our choice
of orientation of $F$.  Indeed, reversing the orientation of $F$
has the effect of negating $e(F)$, {\em and} introducing a second factor
of $-1$ coming from the change in fundamental class.  These two effects
cancel.

For this definition to be satisfactory, we must be able to prove
the following lemma, which is standard in the setting of
ordinary cohomology of compact manifolds.

\begin{lemma}\label{Npush}
Let $i:N\hookto M$ be a closed immersion.
Then for any $\a\in\hhm,\ga\in\hhso(N)$, we have
$\IM\a\cdot i_*\ga = \IN i^*\a\cdot\ga$.
\end{lemma}

\begin{proof}
Let $G=N^{\so}$, let $j:G\to F$ denote the restriction of $i$ to $G$,
and let $\phi:F\to M$ and $\psi:G\to N$ denote the inclusions of $F$ and $G$
into $M$ and $N$, respectively.
$$\begin{CD}
N        @>i>> M\\
@A\psi AA @AA\phi A\\
G        @>j>> F
\end{CD}$$
Then $$\IM\a\cdot i_*\ga = \IF\frac{\phi^*\a\cdot\phi^*i_*\ga}{e(F)},$$
and $$\IN i^*\a\cdot\ga = \IG\frac{\psi^*i^*\a\cdot\psi^*\ga}{e(G)}
=\IG\frac{j^*\phi^*\a\cdot\psi^*\ga}{e(G)}
=\IF\phi^*\a\cdot j_*\left(\frac{\psi^*\ga}{e(G)}\right),$$
where the last equality is simply the integration formula applied
to the map $j:G\to F$ of compact manifolds \cite{AB}.
Hence it will be sufficient to prove that
$$\phi^*i_*\ga=e(F)\cdot j_*\left(\frac{\psi^*\ga}{e(G)}\right)\in\hhso(F).$$
To do this, we will
show that the difference of the two classes lies in
the kernel of $\phi_*$, which we know is trivial because the composition
$\phi^*\phi_*$ is given by multiplication by the invertible
class $e(F)\in\hhso(F)$.
On the left hand side we get
$$\phi_*\phi^*i_*\ga = \phi_*(1)\cdot i_*\ga\hspace{15pt}
\text{by    }(\ref{module}),
$$
and on the right hand side we get
\begin{eqnarray*}
\phi_*\left(e(F)\cdot j_*\left(\frac{\psi^*\ga}{e(G)}\right)\right)
&=& \phi_*\left(\phi^*\phi_*(1)
\cdot j_*\left(\frac{\psi^*\ga}{e(G)}\right)\right)\\
&=& \phi_*(1)\cdot \phi_*j_*\left(\frac{\psi^*\ga}{e(G)}\right)
\hspace{1in}\text{    by    }(\ref{module})\\
&=& \phi_*(1)\cdot i_*\psi_*\left(\frac{\psi^*\ga}{e(G)}\right)
\hspace{1in}\text{    by    }(\ref{functor}).
\end{eqnarray*}
It thus remains only to show that
$\ga = \psi_*\left(\frac{\psi^*\ga}{e(G)}\right)$.
This is seen by applying $\psi^*$ to both sides,
which is an isomorphism (working over the field $\hhp$) by \cite{AB}.
\end{proof}

For $\a_1,\a_2\in\hhm$, consider the
symmetric, bilinear, $\hhp$-valued pairing
$$\la\a_1,\a_2\ra_M=\IM\a_1\a_2.$$

\begin{lemma}[Poincar\'e Duality]\label{perfect}
This pairing is nondegenerate.
\end{lemma}
 
\begin{proof}
Suppose that $\a\in\hhm$ is nonzero, and therefore $\phi^*\a\neq 0$.
Since $F$ is
compact, there must exist a class $\ga\in\hhso(F)$ such
that $0\neq\IF\phi^*\a\cdot\ga=\IM\a\cdot \phi_*\ga=\la\a,\phi_*\ga\ra_M$.
\end{proof}

\vspace{-\baselineskip}
\begin{definition}\label{push}
For an arbitrary equivariant map $f:N\to M$, we may now define
the pushforward $$f_*:\hhn\to\hhm$$ to
be the adjoint of $f^*$ with respect to the pairings
$\la\cdot,\cdot\ra_N$ and $\la\cdot,\cdot\ra_M$.
This is well defined because, according to (\ref{isomorphic}),
$\hhm$ and $\hhn$ are finite dimensional
vector spaces over the field $\hhp$. Lemma~\ref{Npush} tells us that this definition generalizes
the definition for closed immersions.  Furthermore,
properties (\ref{functor}) and (\ref{module}) for pushforwards along
arbitrary maps
are immediate corollaries of the definition.
If $f$ is a projection, then $f_*$ will be given by
integration along the fibers.
Using the fact that
every map factors through its graph as a closed immersion
and a projection, we always have a geometric interpretation of the pushforward.
\end{definition}

As an application, let us consider the manifold $M\times M$,
along with the two projections $\pi_1$ and $\pi_2$,
and the diagonal map $\Diag:M\to M\times M$.
Suppose that we can write $$\Diag_*(1)=\sum\pi_1^*a_i\cdot\pi_2^*b_i$$
for a finite collection of classes $a_i,b_i\in\hhm$.
The following Proposition will be used in Section~\ref{quiver}.

\begin{proposition}\label{basis}
The set $\{b_i\}$ is an additive spanning set for $\hhm$.
\end{proposition}

\begin{proof}
For any $\a\in\hhm$, we have
\begin{eqnarray*}
\a &=& \id_*\id^*\a\\
&=& (\pi_2\circ\Diag)_*(\pi_1\circ\Diag)^*\a\\
&=& \pi_{2*}\big(\Diag_*\left(1\cdot\Diag^*\pi_1^*\a\right)\big)\\
&=& \pi_{2*}\big(\pi_1^*\a\cdot\Diag_*(1)\big)\\
&=& \pi_{2*}\left(\sum\pi_1^*(a_i\a)\cdot\pi_2^*b_i\right)\\
&=& \sum\pi_{2*}\pi_1^*(a_i\a) \cdot b_i\\
&=& \sum\la a_i,\a\ra\cdot b_i,
\end{eqnarray*}
hence $\a$ is in the span of $\{b_i\}$.
\end{proof}
\end{section}

\begin{section}{Hyperk\"ahler abelianization}
Let $M$ be a hyperk\"ahler manifold with a circle action, and
suppose that a compact Lie group $G$ acts hyperhamiltonianly on $M$
with hyperk\"ahler moment map
$$\muh=\mur\oplus\muc:M\to\gd\oplus\gd_{\C},$$
where $\muc$ is holomorphic with respect to the distinguished
complex structure $I$.
We require that the action of $G$ commute with the action of $\so$,
that $\mur$ is $\so$-invariant, and that $\muc$ is $\so$-equivariant
with respect to the action of $\so$ on $\gd_{\C}$ by complex multiplication.

Let $T\subs G$ be a maximal torus, and let $pr:\gd\to\tdu$
be the natural projection.  Then $T$ acts on $M$ with hyperk\"ahler
moment map
$$\mubh=pr\!\circ\!\mur\,\oplus\, pr_{\C}\!\circ\!\muc:
M\to\tdu\oplus\tdu_{\C}.$$
Let $\xi\in\gd$ be a central element such that $(\xi,0)$ is a regular
value of $\muh$ and $(pr(\xi),0)$ is a regular value of $\mubh$.
Assume further that $G$ acts freely on $\muh^{-1}(\xi,0)$,
and $T$ acts freely on $\mubh^{-1}(pr(\xi),0)$.\footnote{We make this
simplifying assumption in order to talk about manifolds, rather than
orbifolds, which makes the integration formulae easier to state.
In fact, Theorems \ref{main} and \ref{ordinary}
generalize easily to the orbifold case, as in \cite[\S 6]{Ma}.}
Let $$M\mmod G = \muh^{-1}(\xi,0)/G\hs\hs\hs\text{   and   }\hs\hs\hs
M\mmod T = \mubh^{-1}(pr(\xi),0)/T$$
be the hyperk\"ahler quotients of $M$ by $G$ and $T$, respectively.
Because $\muh$ and $\mubh$ are circle-equivariant,
the action of $S^1$ on $M$ descends to actions on the hyperk\"ahler quotients.
Note that $M\mmod T$ also inherits an action of the Weyl group $W$ of $G$.

\begin{example}\label{standard}
The main example to keep in mind is $M=\cot$, where $\so$
acts on $M$ by scalar multiplication on the fibers and the action
of $G$ on $M$ is induced by a linear action of $G$ on $\Cn$,
as in Chapter \ref{analogues}.
\end{example}

Consider the Kirwan maps
$$\k_G:H_{\so\times G}^*(M)\to \hso(M\mmod G)\hs\hs\hs\text{   and   }\hs\hs\hs
\k_T:H_{\so\times T}^*(M)\to \hso(M\mmod T),$$
induced by the inclusions of $\muh^{-1}(\xi,0)$ and $\mubh^{-1}(pr(\xi),0)$ into $M$,
along with their rationalizations
$$\hat\k_G:\widehat H_{\so\times G}^*(M)\to \hhso(M\mmod G)\hs\hs\hs\text{   and   }\hs\hs\hs
\hat\k_T:\widehat H_{\so\times T}^*(M)\to \hhso(M\mmod T).$$
Let $$\rgt:\widehat{H}^*_{\so\times G}(M)\to\widehat{H}^*_{\so\times T}(M)^W$$
be the standard isomorphism.

Let $\Delta=\Delta^+\sqcup\Delta^-\subset\tdu$ be the set of roots of $G$.
Let $$e = \Pd_{\a\in\Delta}\a(x-\a) \in (\operatorname{Sym}\tdu)^W\otimes\Q[x]
\cong H_{\so\times G}(pt)\subs\widehat{H}_{\so\times G}(pt),$$
and $$e' = \Pd_{\a\in\Delta^-}\a\cdot\Pd_{\a\in\Delta}(x-\a)
\in\operatorname{Sym}\tdu\otimes\Q[x]\cong H_{\so\times T}(pt)
\subs\widehat{H}_{\so\times T}(pt).$$
The following two theorems are analogues of
Theorems B and A of \cite{Ma}, adapted to circle compact
hyperk\"ahler quotients.  Our proofs follow closely those of Martin.

\begin{theorem}\label{integration}
Suppose that $M\mmod G$ and $M\mmod T$ are both circle compact.
If $\ga\in\hhsg$, then
$$\intxg \Kh_G(\ga) = \frac{1}{|W|}\intxt \Kh_T\circ\rgt(\ga)\cdot e.$$
\end{theorem}

\begin{theorem}\label{main}
Suppose that $M\mmod G$ and $M\mmod T$ are both circle compact,
and that the rationalized Kirwan map $\Kh_G$
is surjective.  Then
$$\hhso(M\mmod G)\cong\frac{\hhso(M\mmod T)^W}{ann(e)}
\cong\left(\frac{\hhso(M\mmod T)}{ann(e')}\right)^W.$$
\end{theorem}

\begin{proofintegration}
Consider the following pair of maps:
$$\begin{CD}
\muh^{-1}(\xi,0)/T  @>i>> \mubh^{-1}(pr(\xi),0)/T\cong M\mmod T\\
@V\pi VV \\
\muh^{-1}(\xi,0)/G\cong M\mmod G.
\end{CD}$$
Each of these spaces is a complex $S^1$-manifold with a compact, complex
fixed point set, and therefore satisfies the hypotheses of Section
\ref{thepush}.
Let $$b = \Pd_{\a\in\Delta^+}\a\in H^*_{\so\times T}(pt)$$ be the product of the
positive roots of $G$, which we will think of as an
element of $\hhxt$.  Martin shows that $\pi_*i^*b=|W|$,
and that $i^*\circ\Kh_T\circ\rgt = \pi^*\Kh_G$ \cite{Ma}, hence
we have
\begin{eqnarray*}
\intxg\Kh_G(\ga) &=& \frac{1}{|W|}\intxg\Kh_G(\ga)\cdot\pi_*i^*b\\
&=& \frac{1}{|W|}\int_{\muh^{-1}(\xi,0)/T}\pi^*\Kh_G(\ga)\cdot i^*b\hspace{20pt}
\text{by Definition~\ref{push}}\\
&=& \frac{1}{|W|}\int_{\muh^{-1}(\xi,0)/T}i^*\circ\Kh_T\circ\rgt(\ga)\cdot i^*b\\
&=& \frac{1}{|W|}\intxt \Kh_T\circ\rgt(\ga)\cdot b \cdot i_*(1)\hspace{20pt}
\text{by Lemma~\ref{Npush}}.
\end{eqnarray*}
It remains to compute $i_*(1)\in\hhxt$.
For $\a\in\Delta$, let $$L_{\a} = \mubh^{-1}((pr(\xi),0)\times_T\C_{\a}$$
be the line bundle on $M\mmod T$ with $S^1$-equivariant
Euler class $\a$.
Similarly, let $L_x$ be the (topologically trivial)
line bundle with $S^1$-equivariant
Euler class $x$.
Following the idea of \cite[1.2.1]{Ma}, we observe
that the restriction of $\muh - (\xi,0)$ to
$\mubh^{-1}(pr(\xi),0)$ defines
an $\so\times T$-equivariant map $$s:\mubh^{-1}(pr(\xi),0)
\to V\oplus V_{\C},$$ where $V = pr^{-1}(0)$
and $V_{\C}=pr_{\C}^{-1}(0)$.
This descends to an $S^1$-equivariant section of the
bundle $$E = \mubh^{-1}(pr(\xi),0)\times_T \left(V\oplus V_{\C}\right)$$
with zero locus $\muh^{-1}(\xi,0)/T$.  The fact that $(\xi,0)$
is a regular value implies that this section is generic,
hence the equivariant Euler class $e(E)\in\hhxt$ is equal to $i_*(1)$.

The vector space $V$ is isomorphic as a $T$-representation to
$\bigoplus_{\a\in\Delta^-}\C_{\a}$,
with $S^1$ acting trivially.
Similarly, $V_{\C}$ is isomorphic to $V\otimes\C\cong V\oplus V^*$,
with $S^1$ acting diagonally by scalars.
Hence
\begin{eqnarray*}
E &\cong& \bigoplus_{\a\in\Delta^-}L_{\a}
\oplus\bigoplus_{\a\in\Delta^-}\left(L_x\otimes L_{\a}\right)
\oplus \left(L_x\otimes L_{-\a}\right)\\
&\cong& \bigoplus_{\a\in\Delta^-}L_{\a}\oplus
\bigoplus_{\a\in\Delta}L_x\otimes L_{-\a},
\end{eqnarray*}
and therefore $$i_*(1)=e(E) = 
\prod_{\a\in\Delta^-}\a\cdot\prod_{\a\in\Delta}(x-\a)=e'.$$
Multiplying by $b$ we obtain $e$, and the theorem is proved.
\end{proofintegration}

%

\begin{proofmain}
Observe that the restriction
of $\pi^*$ to the Weyl-invariant part
$\hhso\left(\muh^{-1}(\xi,0)/T\right)^W$ is given by the composition of
isomorphisms
$$\hhso\left(\muh^{-1}(\xi,0)/T\right)^W
\cong\widehat{H}^*_{\so\times T}\left(\muh^{-1}(\xi,0)\right)^W
\cong\widehat{H}^*_{\so\times G}\left(\muh^{-1}(\xi,0)\right)
\cong\hhxg,$$
hence we may define
$$i^*_W:= (\pi^*)^{-1}\circ i^*:
\hhxt^W\to\hhso\left(\muh^{-1}(\xi,0)/T\right)^W.$$
Furthermore, we have
$\Kh_G = i^*_W\circ\Kh_T\circ\rgt$, hence
$i^*_W$ is surjective.
As in \cite[\S 3]{Ma},
\begin{eqnarray*}
i^*_W(a)=0 &\iff& \forall c\in\hhxt^W, \intxg i^*_W(c)\cdot i^*_W(a)=0\hspace{15pt}
\text{by \ref{perfect} and surjectivity of $i^*_W$}\\
&\iff& \forall c\in\hhxt^W, \intxt c\cdot a\cdot e = 0\hspace{15pt}
\text{by Theorem~\ref{integration}}\\
&\iff& \forall d\in\hhxt, \intxt d\cdot a\cdot e=0\hspace{15pt}
\text{by using $W$ to average $d$}\\
&\iff& a\cdot e=0\hspace{15pt}\text{by Lemma~\ref{perfect}},
\end{eqnarray*}
hence $\ker i^*_W = ann(e)$.
By surjectivity of $i^*_W$,
\begin{eqnarray*}
\hhso(M\mmod G) &\cong&
\frac{\hhxt^W}{\ker i^*_W}\cong\frac{\hhxt^W}{ann(e)}.
\end{eqnarray*}
By a second application of Lemma~\ref{perfect}, for any $a\in\hhxt$, we have
\begin{eqnarray*}
i^*(a) = 0 &\impl& \forall f\in\hhso(\muh^{-1}(\xi,0)/T),
\int_{\muh^{-1}(\xi,0)/T}f\cdot i^*(a) = 0\\
&\impl& \forall c\in\hhxt,
\int_{\muh^{-1}(\xi,0)/T} i^*(c)\cdot i^*(a)=0\\
&\impl& \forall c\in\hhxt,
\intxt c\cdot a\cdot i_*(1)=0\hspace{15pt}\text{by Lemma~\ref{Npush}}\\
&\impl& a\cdot e'= a\cdot i_*(1) = 0\hspace{15pt}\text{by Lemma~\ref{perfect}},
\end{eqnarray*}
hence $\ker i^* \subs ann(e')$.
This gives us a natural surjection
$$\frac{\hhxt^W}{ann(e)}=\frac{\hhxt^W}{\ker i^*_W}
\cong \left(\frac{\hhxt}{\ker i^*}\right)^W
\to \left(\frac{\hhxt}{ann(e')}\right)^W,$$
which is also injective because $e'$ divides $e$.
This completes the proof of Theorem \ref{main}.
\end{proofmain}

For the non-rationalized version of Theorem~\ref{main}, we make the additional
assumption that $M\mmod G$ and $M\mmod T$ are equivariantly formal
$\so$-manifolds, i.e. that $\hxg$ and $\hxt$ are free modules
over $\hp$.  Proposition \ref{formality} tells us that this 
is the case whenever the circle action
is hamiltonian and its moment map is proper and bounded below.

\begin{theorem}\label{ordinary}
Suppose that $M\mmod G$ and $M\mmod T$ are equivariantly formal, circle compact, and
that the rationalized Kirwan map $\Kh_G$ is surjective.
Then
$$\hso(M\mmod G)\supseteq\Im(\K_G)\cong\frac{(\Im\K_T)^W}{ann(e)}
\cong\left(\frac{\Im\K_T}{ann(e')}\right)^W.$$
\end{theorem}

\begin{remark}\label{notsobad}
In the context of Example~\ref{standard} with $pr\circ\mu$ proper,
$M\mmod G$ and $M\mmod T$ are both circle compact and equivariantly
formal (Proposition \ref{formality}) and $\K_T$ is always surjective 
(Theorem \ref{hm}).
\end{remark}

\begin{proofordinary}
Consider the following exact commutative diagram
$$
\xymatrix{
0 \ar[r] & A \ar[r]\ar[d] & \hxt^W \ar[r]^{i_W^*} \ar[d]
& \hxg \ar[d]\\
0 \ar[r] & \widehat{A} \ar[r] & \hhxt^W \ar[r]^{i_W^*} & \hhxg.}
$$
Equivariant formality implies that the downward maps in the above
diagram are inclusions, hence the map on top labeled $i^*_W$
is simply the restriction of the map on the bottom
to the subring $\hxt\subs\hhxt$.
We therefore have $$A = \widehat{A}\cap\hxt^W = ann(e).$$
Just as in the rationalized case, we have $\K_G = i_W^*\circ\K_T\circ\rgt$,
hence
$$\Im(\K_G)\cong i_W^*\left(\Im\K_T\circ\rgt\right)
\cong \frac{(\Im\K_T)^W}{ann(e)}.$$

Now consider the analogous diagram
$$
\xymatrix{
0 \ar[r] & B \ar[r]\ar[d] & \hxt \ar[r]^{\!\!\!\!\!\!\!\!\! i^*} \ar[d]
& \hso(\muh^{-1}(\xi,0)/T) \ar[d]\\
0 \ar[r] & \widehat{B} \ar[r] & \hhxt \ar[r]^{\!\!\!\!\!\!\!\!\! i^*} 
& \hso(\muh^{-1}(\xi,0)/T).}
$$
Since we have not assumed that $\muh^{-1}(\xi,0)/T$ is equivariantly formal,
we only know that the first two downward arrows are inclusions, and hence
can only conclude that $B$ is contained in the annihilator of $e'$.
Since $e'$ divides $e$, we have a series of natural surjections
$$\frac{\imkt^W}{ann(e)}\cong
\frac{\imkt^W}{A}\cong
\left(\frac{\Im\K_T}{B}\right)^W\to
\left(\frac{\Im\K_T}{ann(e')}\right)^W\to
\left(\frac{\Im\K_T}{ann(e)}\right)^W.$$
The composition of these maps is an isomorphism, hence so is each one.
\end{proofordinary}
\end{section}
\end{chapter}

\begin{chapter}{Hyperpolygon spaces}\label{hpspaces}
A hyperpolygon space is the hyperk\"ahler analogue of a polygon space,
which parameterizes $n$-sided polygons in $\R^3$ with fixed edge lengths.
It is also an example of a quiver variety, introduced by Nakajima \cite{N1,N2},
and since studied by many authors.  In Section \ref{quiver} we give the basic
constructions of quiver varieties and hyperpolygon spaces, and show that they
satisfy all of the hypotheses of Chapter \ref{abelianization}.
Section \ref{modcore} is devoted to understanding the components of the core
of a hyperpolygon space; in particular, we show that they are all smooth
(Theorem \ref{component}), and interpret them as moduli spaces of spatial polygons
with certain properties (Theorem \ref{us}).
Sections \ref{s1} and \ref{last} contain computations of the $\so$-equivariant
cohomology rings of the hyperpolygon space as well as its core components,
making use of the abelianization technique of Chapter \ref{abelianization}.

This chapter is taken from \cite{HP2} and \cite{HP}.  The reader is warned
that our notation differs significantly from that of \cite{HP2}; most glaring
is the fact that our spaces $\X$ and $\M$ correspond to the spaces 
$M$ and $X$ (respectively) in \cite{HP2}.  This abrupt switch is necessary to
conform with the conventions of Chapters \ref{analogues}-\ref{abelianization}.

\begin{section}{Quiver varieties}\label{quiver}
Let $Q$ be a quiver with vertex set $I$ and edge set $E\subs I\times I$,
where $(i,j)\in E$ means that $Q$ has an arrow pointing from $i$ to $j$.
We assume that
$Q$ is connected and has no oriented cycles.
Suppose given two collections of vector spaces $\{V_i\}$ and $\{W_i\}$,
each indexed by $I$, and consider the affine space
$$\Af = \bigoplus_{(i,j)\in E}\Hom(V_i,V_j)
\,\oplus\,\bigoplus_{i\in I}\Hom(V_i,W_i).$$
The group $U(V) = \prod_{i\in I}U(V_i)$ acts on $\Af$ by conjugation, and
this action is hamiltonian.  Given an element
$$(B,J) = \bigoplus_{(i,j)\in E}B_{ij}\oplus\bigoplus_{i\in I}J_i$$
of $\Af$,
the $\mathfrak{u}(V_i)^*$ component of the moment map
is $$\mu_i(B,J) = J_i^{\dagger}J_i+\sum_{(i,j)\in E}B_{ij}^{\dagger}B_{ij},$$
where $\dagger$ denotes adjoint, and $\mathfrak{u}(V_i)^*$
is identified with with the set of hermitian matrices via
the trace form.
Given a generic central element $\xi\in\mathfrak{u}(V)^*$,
the K\"ahler quotient $\Af\,\mod_{\!\xi}U(V)$ parameterizes 
isomorphism classes of $\xi$-stable,
framed representations
of $Q$ of fixed dimension \cite{N2}.  If $W_i=0$ for all $i$, then
the diagonal circle $U(1)$ in the center of $U(V)$ acts trivially,
and we instead quotient by
$PU(V) = U(V)/U(1)$.

Consider the hyperk\"ahler quotient $$\M = T^*\Af\,\mmod_{\!(\xi,0)}U(V),$$
or, if $W_i=0$ for all $i$,
$$\M = T^*\Af\,\mmod_{\!(\xi,0)}PU(V).$$
As in Section~\ref{circleaction}, $\M$ has a natural action 
of $\gm$ induced from scalar
multiplication on the fibers of $T^*\Af$.
We now show that $M=T^*\Af$ satisfies
the hypotheses of Theorems~\ref{main} and \ref{ordinary}.

\begin{proposition}
Let $T(V)\subs U(V)$ be a maximal torus, and let $pr:\mathfrak{u}(V)^*\to\mathfrak{t}(V)^*$
be the natural projection.
The moment maps $\mu=\Od_{i\in I}\mu_i:\Af\to\mathfrak{u}(V)^*$
and $pr\circ\mu:\Af\to\mathfrak{t}(V)^*$ are each proper.
\end{proposition}

\begin{proof}
To show that $\mu$ and $pr\circ\mu$ is proper, it suffices to find an element
$t\in T(V)\subs U(V)$ such that the weights of the
action of $t$ on $\Af$ are all strictly positive.
Let $\lambda=\{\lambda_i\mid i\in I\}$ be a collection of integers,
and let $t\in T(V)$ be the central element of $U(V)$ that acts
on $V_i$ with weight $\lambda_i$ for all $i$.
Then $t$ acts on $\Hom(V_i,V_j)$ with weight $\lambda_j-\lambda_i$,
and on $\Hom(V_i,W_i)$ with weight $-\lambda_i$.
Hence we have reduced the problem to showing that it is possible
to choose $\lambda$ such that $\lambda_i<0$ for all $i\in I$ and
$\lambda_i<\lambda_j$ for all $(i,j)\in E$.

We proceed by induction on the order of $I$.
Since $Q$ has no oriented cycles, there must exist a source $i\in I$;
a vertex such that for all $j\in I$, $(j,i)\notin E$.
Deleting $i$ gives a smaller (possibly disconnected) quiver with
no oriented cycles, and therefore we may choose
$\left\{\lambda_j\mid j\in I\smallsetminus\{i\}
\right\}$ such that $\lambda_j<0$ for all $j\in I\smallsetminus\{i\}$
and $\lambda_j<\lambda_k$ for all $(j,k)\in E$.
We then choose $\lambda_i<\min\left\{\lambda_j\mid j\in I\smallsetminus\{i\}\right\}$,
and we are done.
\end{proof}

\vspace{-\baselineskip}
\begin{proposition}
The rationalized Kirwan map
$\Kh_{U(V)}:\widehat{H}^*_{\so\times U(V)}\Big(T^*\Af\Big)\to\hhso(\M)$
is surjective.
\end{proposition}

\begin{proof}
Nakajima \cite[\S 7.3]{N2} shows that there exist cohomology classes
$a_i, b_i$ in the image of $\Kh_{U(V)}$
such that $\dso = \sum\pi_1^*a_i\cdot\pi_2^*b_i$.
(Nakajima uses a slightly different circle action, but
his proof is easily adapted to the circle action that we have defined.)
It follows from Proposition~\ref{basis} that
the classes $\{b_i\}$ generate $\hhso(\M)$.
\end{proof}

\vspace{-\baselineskip}
\begin{remark}\label{generation}
This Proposition shows that the assumptions of 
Theorems~\ref{integration},~\ref{main},
and~\ref{ordinary} are satisfied for Nakajima's quiver varieties.
Thus integration in equivariant cohomology yields a
description of the rationalized $\so$-equivariant cohomology, and also of
the image of the non-rationalized Kirwan map $\K_G$. Therefore
if we know that $\K_G$ is surjective for a particular quiver variety,
then we have a concrete description of the ($\so$-equivariant)
cohomology ring of that quiver variety.
It is known that $\K_G$ is surjective
for Hilbert schemes of $n$ points on an ALE space,
so our theory applies and gives a description
of the cohomology ring of these quiver varieties.
More examples of quiver varieties with surjective Kirwan map,
including hyperpolygon spaces,
are discussed in Remark \ref{horn}.
\end{remark}

\begin{remark}
Another interesting application of Proposition~\ref{basis} is
to the moduli space ${\mathcal M}$ of stable rank $n$ and degree $1$
Higgs bundles  on a genus $g>1$ smooth projective algebraic curve $C$ (see \cite{H2}).
It is an easy exercise to write down the cohomology class 
of the diagonal in ${\mathcal M}\times {\mathcal M}$
as an expression in a certain set of tautological classes.
Proposition~\ref{basis} implies that
the rationalized $\so$-equivariant cohomology ring $\hhso({\mathcal M})$ 
is generated by these classes.
In fact the same result follows from the argument of \cite{HT1}. 
There ${\mathcal M}$ was embedded into a circle compact manifold
${\mathcal M}_\infty$, whose cohomology is the free algebra on the 
tautological classes. The argument in \cite{HT1} then goes by
showing that the embedding of the $\so$-fixed point set of ${\mathcal M}$ 
in that of ${\mathcal M}_\infty$ induces a surjection
on cohomology. This already implies that $\hhso({\mathcal M}_\infty)$ 
surjects onto $\hhso(\mathcal M)$. In \cite{HT1} it is shown
that in the rank $2$ case this embedding also implies the surjection on 
non-rationalized cohomology, and then
a companion paper \cite{HT2} describes the cohomology ring of ${\mathcal M}$ explicitly. 
However for higher rank
this part of the argument of \cite{HT1} breaks down. 
Later Markman \cite{Mk} used similar diagonal arguments on
certain compactifications of $\mathcal M$ and
Hironaka's celebrated theorem on desingularization of algebraic 
varieties to deduce that the cohomology ring of ${\mathcal M}$
is generated by tautological classes for all $n$.
\end{remark}

A {\em hyperpolygon space}, introduced in \cite{K2},
is a quiver variety associated to the
following quiver (Figure \ref{fig:quiver}),
with $V_0 = \C^2$, $V_i=\C^1$ for $i\in\{1,\ldots,n\}$, and $W_i=0$ for all $i$.
\begin{figure}[h]
\centerline{\epsfig{figure=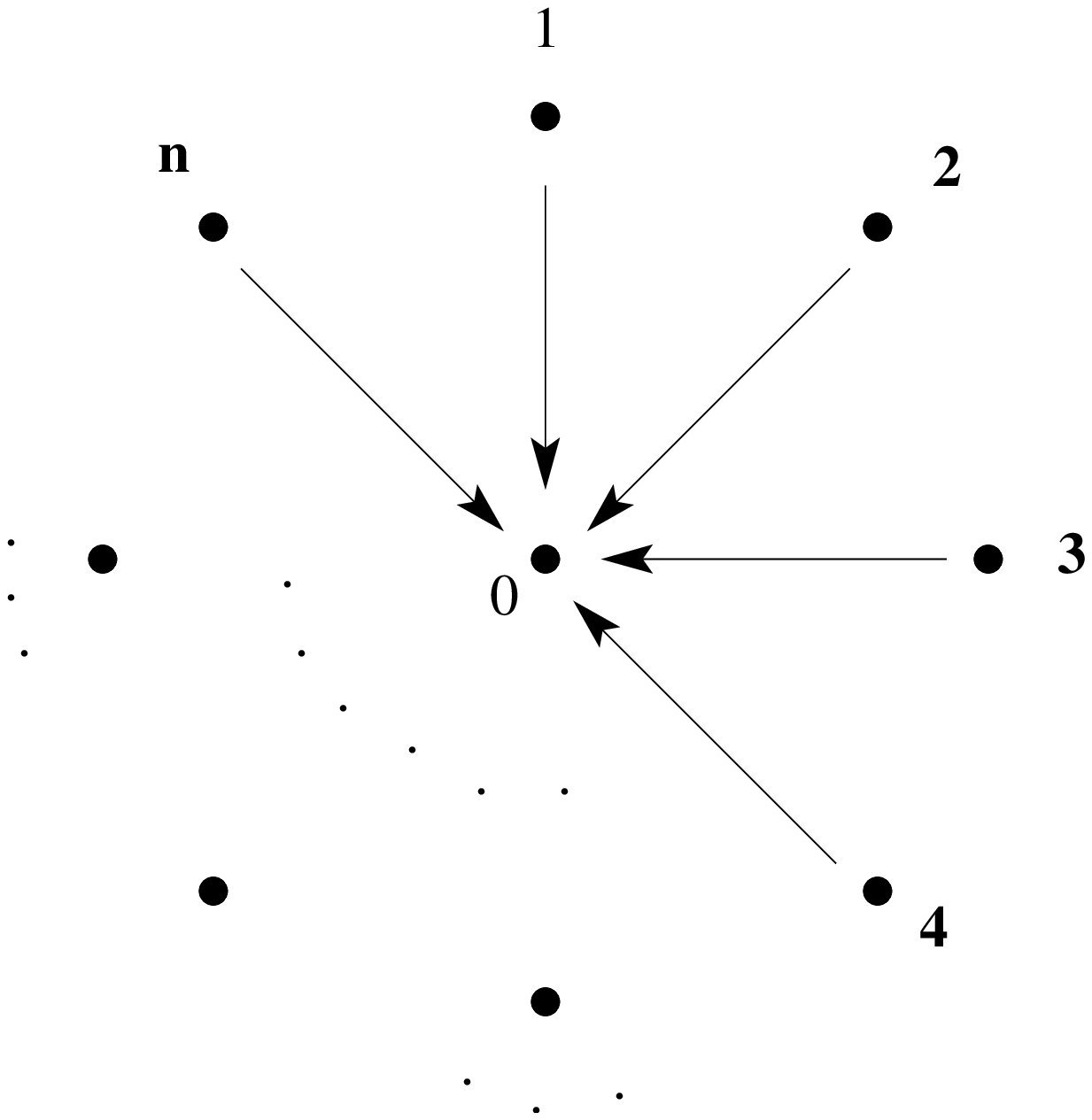,height=4cm}}
\caption{The quiver for a hyperpolygon space.}\label{fig:quiver}
\end{figure}

Let $$G := PU(V) = \Big(SU(2)\times U(1)^n\Big)\Big/\Z_2,$$
and $$\GC := PGL(V) = \Big(\SL\times (\gm)^n\Big)\Big/\Z_2,$$
where $\Zt$ acts by multiplying each factor by $-1$.
We represent an element of $\Af\cong\ctn$ by an $n$-tuple
of column vectors $$q = (q_1,\ldots,q_n).$$
Following the conventions in \cite{K2}, we consider the {\em right}
action of $G\subs\GC$ on $\Af$ given explicitly by
$$q[\Theta;e_1,\ldots,e_n] = (\Theta^{-1}q_1e_1,\ldots,\Theta^{-1}q_ne_n).$$
The compact group $G$ acts with moment map $\mu:\ctn\to\sutd\oplus\tnd$
given by the equation $$\mu(q) = \sum_{i=1}^n (q_iq_i^*)_0 \oplus 
\left(\half|q_1|^2,\ldots,\half|q_n|^2\right),$$
where $q_i^*$ denotes the conjugate transpose of $q_i$, $(q_iq_i^*)_0$ 
denotes the
traceless part of $q_iq_i^*$, and $\sutd$ is identified with 
$i\cdot\sut$ via the trace form.
Given an $n$-tuple of real numbers $(\a_1,\ldots,\a_n)$, 
we define the {\em polygon space} $$\X(\a) := \ctn\mod_{\!\!\a}G,$$
where $\a = 0\oplus (\a_1,\ldots,\a_n)\in\sutd\oplus\tnd$.
If we break the reduction into two steps, reducing first by $U(1)^n$
and then by $SU(2)$, we find that
\begin{equation}\label{polygon}
\X(\a)\cong\left\{(v_1,\ldots,v_n)\in (\R^3)^n\,\bigg|\, 
\|v_i\|=\a_i\text{ and }\sum v_i=0\right\}\bigg/SO(3)
\end{equation}
(see Remark \ref{interp} and the proof of Theorem \ref{us}).
Here $\sutd$ is being identified with $\R^3$, 
and the coadjoint action of $SU(2)$
on $\sutd$ is being replaced by the standard action of $SO(3)$ on $\R^3$ \cite{HK}.
This space, therefore, may be thought of as the moduli space of $n$-sided polygons in $\R^3$,
with fixed edge lengths $(\a_1,\ldots,\a_n)$, up to rotation.  
In particular, $\X(\a)$ is empty unless $\a_i\geq 0$ for all $i$.

We call $\a$ {\em generic} if there does not exist a
subset $S\subs\{1,\ldots,n\}$
such that $\sum_{i\in S}\a_i=\sum_{j\in S^c}\a_j$.
Geometrically, this means that there is no element of $\X(\a)$ represented by a polygon that is
contained in a single line in $\R^3$.  If $\a$ is generic, then $\X(\a)$ is smooth \cite{HK}.
Throughout this paper we will assume that $\a$ is generic, and that $\a_i> 0$ for all $i$.

To define the hyperk\"ahler analogue of $\X(\a)$, we consider the induced action of $G$
on $\cot$.
Explicitly, we write an element of $T^*\Af$ as $(p,q)$,
where $q = (q_1,\ldots,q_n)$ is an $n$-tuple
of column vectors
and $p = (p_1,\ldots,p_n)$ an $n$-tuple
of row vectors, and we put
$$(p,q)[\Theta;e_1,\ldots,e_n] = 
\big((e_1^{-1}p_1\Theta,\ldots,e_n^{-1}p_n\Theta),(\Theta^{-1}q_1e_1,\ldots,\Theta^{-1}q_ne_n)\big).$$
The action of $G$ on $T^*\Af$ is hyperhamiltonian
with hyperk\"ahler moment map
$$\mur\oplus\muc:\cot\to\Big(\sutd\oplus\tnd\Big)\oplus\Big(\slt^*\oplus
(\mathfrak{u}(1)^n_{\C})^*\Big)$$
given by the equations
\begin{equation*}\label{mur}
\mu_{\R}(p,q)  =  \frac{\sqrt{-1}}{2}\sum_{i=1}^n \left(q_i q_i^* - p_i^* p_i\right)_0
\oplus\left(\half\left(|q_1|^2 - |p_1|^2\right), \ldots, \half\left(|q_n|^2 - |p_n|^2\right)\right) 
\end{equation*}
and
\begin{equation*}\label{muc}
\mu_{\C}\left(p,q\right)  =  - \sum_{i=1}^n \left(q_i p_i\right)_0\oplus \left(\sqrt{-1} p_1 q_1,
\ldots, \sqrt{-1} p_n q_n\right).
\end{equation*}
We then define the hyperpolygon space to be the hyperk\"ahler quotient
$$\M(\a) := \cot\mmod_{(\a,0)}G = 
\Big(\mur^{-1}(\a)\cap\muc^{-1}(0)\Big)\Big/G,$$
a smooth, noncompact hyperk\"ahler manifold of complex dimension $2(n-3)$.
Recall that we also have
$$\X(\a)\cong\left(\ctn\right)^{\ses}\!\big/G_{\C}
\hspace{.7cm}\text{ and }\hspace{.7cm}
\M(\a)\cong\muc^{-1}(0)^{\ses}\big/G_{\C},$$
where $\ses$ means stable with respect to the weight $\a$ in the sense of
geometric invariant theory (see Section \ref{reduction}).\footnote{Recall
from Theorem \ref{hkhs} that the notions of stability and semistability
agree for generic $\a$.}
Nakajima gives a stability criterion for general quiver varieties \cite{N1,N2},
which Konno interprets in the special case of hyperpolygon spaces.
Call a subset $S\subs\{1,\ldots,n\}$ {\em short} if $\sum_{i\in S}\a_i < \sum_{j\in S^c}\a_j$,
and call it {\em long} if its complement is short.  (Assuming that $\a$ is generic is equivalent to assuming
that every subset is either short or long.)
Given a point $(p,q)\in\cot$ and a subset $S\subs\{1,\ldots,n\}$,
we will say that $S$ is {\em straight} 
in $(p,q)$ if $q_i$ is proportional to $q_j$ for every $i,j\in S$.
The terminology comes from K\"ahler polygon spaces, in which this condition
is equivalent to asking that the vectors $v_i$ and $v_j$ be proportional
over $\R_+$, or that the edges of lengths $\a_i$ and $\a_j$ (if they happen to be adjacent)
line up to make a single edge of length $\a_i+\a_j$,
as in Figure \ref{straight}.
\begin{figure}[h]
\begin{center}
\psfrag{a1}{$\a_1$}
\psfrag{a2}{$\a_2$}
\psfrag{a3}{$\a_3$}
\psfrag{a4}{$\a_4$}
\psfrag{a5}{$\a_5$}
\psfrag{a6}{$\a_6$}
\psfrag{a7}{$\a_7$}
\psfrag{a8}{$\a_8$}
\includegraphics[height=40mm]{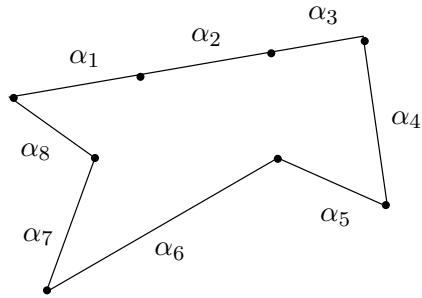}
\caption{The subset $\{1,2,3\}$ is straight.}\label{straight}
\end{center}\end{figure}

\begin{theorem}\label{stability}{\em\cite[4.2]{K2}}
Suppose that $\a$ is generic, and $\a_i> 0$ for all $i$.
Then a point $(p,q)\in\cot$ is stable with respect to $\a$
if and only if the following two conditions are satisfied:
\begin{eqnarray*}
&1)& q_i\neq 0 \text{ for all }i,\text{ and}\\
&2)& \text{if $S$ is straight and $p_j=0$ for all $j\in S^c$, then $S$ is short.}
\end{eqnarray*}
\end{theorem}
As in Chapter \ref{analogues}, we will use the notation $[p,q]$ to denote
the $\GC$-equivalence class of a point $(p,q)\in\muc^{-1}(0)^{\ses}$,
and $[p,q]_{\R}$ to denote
the $G$-equivalence class of a point $(p,q)\in\mur^{-1}(\a)\cap\muc^{-1}(0)$.
Recall that $\X$ sits inside of $\M$ 
as the locus of points $[p,q]$ with $p=0$.
This observation, along with Theorem \ref{stability}, allows us to recover
the stability condition for the action of $G$ on $\ctn$.
A point $q\in\ctn$ is stable if and only if $q_i\neq 0$ for all $i$,
and no long subset is straight, as first shown in \cite{Kl}.
The polygonally-minded reader is warned that in the hyperpolygon space $\M(\a)$,
long subsets {\em can} be straight.
\end{section}

\begin{section}{Moduli theoretic interpretation of the core}\label{modcore}
For the rest of the section we fix a generic 
$\a = 0\oplus (\a_1,\ldots,\a_n)\in\sutd\oplus\tnd$,
with $\a_i>0$ for all $i$, and write $\M=\M(\a)$, $\X=\X(\a)$.
Following Konno, we define
$$\mathcal{S} = \big\{S\subs\{1,\ldots,n\}\hs\hs\big{|}\hs\hs S\text{ is short}\big\}$$
and $$\spa = \big\{S\in\mathcal{S}\hs\hs\big{|}\hs\hs |S|\geq 2\big\}.$$

\begin{theorem}\label{fixed}{\em\cite{K2}}
The fixed point set 
$\displaystyle{\M^{\gm} = \M^{S^1} = \pola\cup\bigcup_{S\in\spa}\M_S,}$
where
$$\M_S = \big\{[p,q]\hs\hs\big{|}\hs\hs S\text{ and }S^c\text{ are each straight, and }
p_j = 0\text{ for all }j\in S^c\big\}.$$  Furthermore, $\M_S$
is diffeomorphic to $\C P^{|S|-2}$.
\end{theorem}

For all $S\in\spa$, let
$U_S = U_{\polsa}$ be the piece of the core $\core\subs\M$
defined in Section \ref{circleaction}.
{\em A priori} we know only that $U_S$ is an irreducible, isotropic
subvariety of dimension at most $n-3$ (Proposition \ref{coreprops}).

\begin{theorem}\label{component}
The core component $U_S$ is smooth of complex
dimension $n-3$, and we have
$$U_S = \big\{[p,q]\bigmid S\text{ is straight, and }
p_j = 0\text{ for all }j\in S^c\big\}.$$
\end{theorem}

Before proving Theorem \ref{component}, 
we describe the way in which the various components
of the core fit together.
For all $S\in\spa$, let $$\polsa = U_S\cap \pola =
\big\{[0,q]\hs\hs\big{|}\hs\hs S\text{ is straight}\big\}.$$
We call this space the {\em polygon subspace} of $\pola$ 
corresponding to the short subset $S$.
Note that $\polsa$ is itself a polygon space with $n-|S|+1$ edges, of lengths
$\{\a_j\mid j\in S^c\}\cup\{\sum_S\a_i\}$.  In particular, it is smooth.
Now suppose given two short subsets $S,T\in\spa$, 
and consider the intersection $U_S\cap U_T$.
\begin{itemize}
\item If $S\cap T=\emptyset$, then $U_S\cap U_T = \polsa\cap X_T$, a polygon subspace
both of $\polsa$ and of $X_T$.
\item If $S\cap T\neq\emptyset$
and $S\cup T$ is long, then $U_S\cap U_T=\emptyset$.
\item If $S\cap T\neq\emptyset$ and 
$S\cup T$ is short,
then $$U_S\cap U_T = \big\{[p,q]\bigmid S\cup T\text{ is straight, and }p_j=0\text{ for all }
j\in (S\cap T)^c\big\}.$$
This is a subvariety of $U_{S\cup T}$ given by taking the closure inside of $U_{S\cup T}$
of a certain subbundle of the conormal bundle to $X_{S\cup T}\subs\pola$,
defined by setting $p_j=0$ for all $j\in (S\cap T)^c\supseteq (S\cup T)^c$.
\end{itemize}
Each of these descriptions generalizes to higher intersections
without modification.

Finally, we compute the fixed point set $U_S^{\gm}$.
If $[p,q]\in U_S^{\gm}$, then either $p=0$ and $[p,q]\in\polsa$,
or $[p,q]\in \X_T$ for some $T\in\spa$.  If $[p,q]\in \X_T$ then Theorem \ref{fixed}
tells us that $T$ and $T^c$ are each straight, hence $S\subs T$ or $S\subs T^c$.
Since $p\neq 0$, we must have $S\subs T$.  Indeed, $U_S\cap \X_T$ is the linear subspace
of $\X_T\cong \C P^{|T|-2}$ given by the condition $p_j=0$ for all $j\in T\smallsetminus S$.
In particular, $U_S\cap \X_T$ is isomorphic to $\C P^{|S|-2}$ for any $T\supseteq S$.

\begin{example}\label{ooh}
Let $n=5$, $\a_1=\a_2=1$, and $\a_3=\a_4=\a_5=3$, and consider the short subset
$S=\{1,2\}$.  The fixed point set of $U_S$ consists of $\polsa\cong\C P^1$,
and four points $\X_S$, $U_S\cap \X_{T_3}$, 
$U_S\cap \X_{T_4}$, and $U_S\cap \X_{T_5}$,
where $T_j = \{1,2,j\}$ for $j=3,4,5.$  
For each $j$, $U_S\cap U_{T_j}$ is isomorphic to $\C P^1$,
and touches $\polsa$ at the point $X_{T_j}$.
In the following picture, an ellipse represents a copy of $\C P^1$ flowing between two
fixed points, where the numbers or pairs of numbers
indicate subsets that are straight on this $\C P^1$.  (For example, 
$12,45$ means that $1$ and $2$ are straight,
as are $4$ and $5$.)
We will revisit this example at the end of Section \ref{last}.
\begin{figure}[h]
\begin{center}
\psfrag{Xs}{$\M_S$}
\psfrag{Xt3}{$\M_{T_3}$}
\psfrag{Xt4}{$\M_{T_4}$}
\psfrag{Xt5}{$\M_{T_5}$}
\psfrag{Ms(a)}{$\X_{S}$}
\psfrag{Mt3(a)}{$\X_{T_3}$}
\psfrag{Mt4(a)}{$\X_{T_4}$}
\psfrag{Mt5(a)}{$\X_{T_5}$}
\includegraphics[height=40mm]{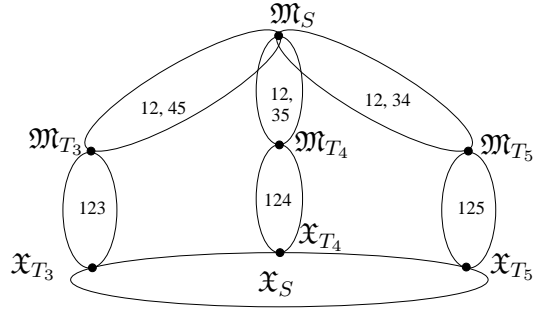}
\caption{$U_S$, with $S=\{1,2\}$}
\label{surface}
\end{center}
\end{figure}
\end{example}

\begin{componentproof}
%
Consider a point $[p,q]\in \M$ with $S$ straight, and $p_j=0$ for all $j\in S^c$.
By applying an element of $G$, we may assume that $q_i = \binom{1}{0}$ for all $i\in S$.
Suppose further that there exists an $i\in S$ with $p_i\neq 0$, and that
no strict superset of $S$ is straight.  In other words, if 
$q_j = \binom{a_j}{b_j}$ for $j\in S^c$,
suppose that $b_j\neq 0$.
For $t\in\gm$, let $\Theta(t) = 
\footnotesize{\Big(\begin{array}{cc}
t&0\\
0&t^{-1}
\end{array}
\Big)}$, 
let $e_i(t)=t$ for all $i\in S$,
and let $e_j(t)=t^{-1}$ for all $j\in S^c$.
Then for $i\in S$, we have $e_i(t)^{-1}p_i \Theta(t) = t^{-2} p_i$ and 
$\Theta(t)^{-1}q_i e_i = q_i$.
For $j\in S^c$, we have $\Theta(t)^{-1}q_j e_j = \binom{t^{-2}a_j}{b_j}$.
Hence
\begin{eqnarray*}
\lim_{t\to\infty}t\cdot[p,q] &=& \lim_{t\to\infty}t^2\cdot[p,q]\\
&=& \lim_{t\to\infty}[t^2 p,q]\\
&=& \lim_{t\to\infty}[t^2 e(t)^{-1}p \Theta(t), \Theta(t)^{-1}q e(t)]\\
&=& [p,q'],
\end{eqnarray*}
where $q'_i = q_i$ for $i\in S$, and $q'_j = \binom{0}{b_j}$ for $j\in S^c$.
Since we have assumed that $b_j\neq 0$ for all $j\in S^c$ and that $p_i\neq 0$
for some $i\in S$, $(p,q')$ is stable, and hence defines an element of $\M_S$.
Since $U_S$ is defined to be the closure of the set of elements that flow up to $\M_S$,
it includes all $[p,q]$ with $S$ straight and $p_j=0$ for all $j\in S^c$.
By dimension count, this containment is an equality, and
we have $\dim U_S = n-3$.

To see that $U_S$ is smooth, it is sufficient to show that $U_S$ is smooth
at $[p,q]$ for all $[p,q]\in \M^{\gm}$.
First suppose that $[p,q]\in \M_T$ for some $T\in\spa$ containing $S$.
Suppose, without loss of generality, that
$T=\{1,\ldots,l\}$ and $S=\{1,\ldots,m\}$ for some $l\leq m$.
Konno computes an explicit local complex chart for $\M$ at the point $[p,q]$,
with coordinates $\{z_i, w_i\mid 3\leq i\leq n-1\}$ \cite{K2}.
With respect to these coordinates, a point $[p',q']$ has the property that
$S$ is straight and $p'_j=0$ for all $j\in S^c$ if and only if 
$w_i=0$ for all $3\leq i\leq l$
and $z_j=0$ for all $l+1 \leq j\leq n-1$.
Hence $U_S$ is smooth at $[p,q]$.

It remains only to show that $U_S$ is smooth at $\polsa = U_S\cap\pola$.
Let $$E = \{(p,q)\mid S\text{ is straight, }
p_j = 0\text{ for all }j\in S^c,\text{ and }\mu_{\C}(p,q)=0\},$$ and let 
$N = \{(p,q)\in E\mid p=0\}$.  The natural projection from $E$ to $N$
exhibits $E$ as a vector bundle over $N$, because the equation $\mu_{\C}(p,q)=0$
is linear in $p$.  We have $U_S = E\mod G = E^{\ses}/G$, and 
$\polsa = N\mod G = N^{\ses}/G$.  
The set $E\vert_{N^{\ses}}/G\subs E^{\ses}/G$ is
an open neighborhood of $\polsa$ inside of $U_S$, and is isomorphic to a vector
bundle over $\polsa$.
Since $\polsa$ is a polygon space it is smooth, hence $U_S$ is smooth 
in a neighborhood of $\polsa$.
\end{componentproof}

\vspace{-\baselineskip}
\begin{corollary}\label{compactification}
$U_S$ is a compactification of the conormal
bundle to $\polsa$ in $\pola$.
\end{corollary}

\begin{proof}
Choose a point $[q,0]\in\polsa$, and a decomposition
$$T_{[q,0]}\M = \nu_1\oplus \nu_2\oplus T_{[q,0]}\polsa\oplus E,$$
where $\nu_1$ is the normal space to $\polsa$ inside of $U_S$,
and $\nu_2$ is the normal space to $\polsa$ inside of $\pola$.
Proposition \ref{coreprops} tells us that $U_S$ and $\pola$
are both $\omega_{\C}$-lagrangian submanifolds of $\M$,
hence $\omega_{\C}$ gives a perfect pairing between $T_{[q,0]}\polsa$
and $E$.  It follows that $\omega_{\C}$ also gives a perfect pairing between
$\nu_1$ and $\nu_2$, and therefore that the normal bundle to
$\polsa$ inside of $U_S$ is dual to the normal bundle of $\polsa$
inside of $\pola$.
\end{proof}

\vspace{-\baselineskip}
\begin{remark}
This argument generalizes to the smooth intersection of any two
lagrangian submanifolds of a symplectic manifold.
\end{remark}

We next describe $U_S$ in polygon-theoretic terms, as a certain moduli
space of pairs of polygons in $\R^3$.

\begin{theorem}\label{us}
The core component $U_S$ is 
homeomorphic to the moduli space
of $n+1$ vectors $$\{u_i, v_j, w\in\R^3\mid i\in S, j\in S^c\},$$
taken up to rotation, satisfying the following conditions:
\begin{eqnarray*}
&1)&\hs\hs w + \sum_{j\in S^c}v_j = 0\\
&2)&\hs\hs \sum_{i\in S}u_i = 0\\
&3)&\hs\hs u_i\cdot w = 0\hs\hs\text{ for all }i\in S\\
&4)&\hs\hs \|v_j\| = \a_j\hs\hs\text{ for all }j\in S^c\\
&5)&\hs\hs \|w\| = \sum_{i\in S}\sqrt{\a_i^2 + \|u_i\|^2}.
\end{eqnarray*}
\end{theorem}

\begin{remark}\label{interp}
In more descriptive terms, a point in $U_S$ specifies two polygons in $\R^3$,
as in Figure~\ref{fig:Uspoly}.
\begin{figure}[h]
\begin{center}
\psfrag{w}{$w$}
\psfrag{v1}{$v_1$}
\psfrag{v2}{$v_2$}
\psfrag{v3}{$v_3$}
\psfrag{vnS}{$v_{n-|S|}$}
\psfrag{u1}{$u_1$}
\psfrag{u2}{$u_2$}
\psfrag{uS}{$u_{|S|}$}
\includegraphics[height=60mm]{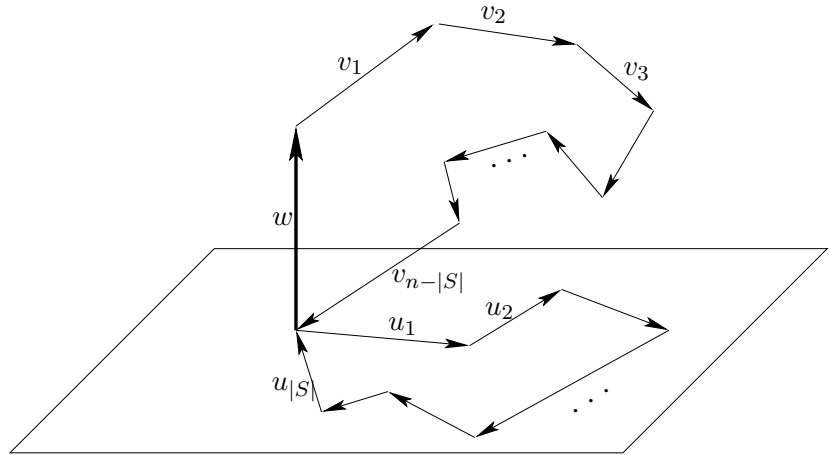}
\caption{
An element of $U_S$, represented by a spatial polygon with a distinguished
edge, and a planar polygon perpendicular to that edge.}
\label{fig:Uspoly}
\end{center}
\end{figure}
The first is the $n - |S| + 1$ sided polygon consisting of the vectors
$\{v_j\mid j\in S^c\}$ and $w$.  Each vector $v_j$ has length $\a_j$, and
$w$ has a variable length, always greater than or equal to $\sum_{i\in S}\a_i$.
This variable length is determined by the lengths of the edges in the second polygon, 
which consists of $|S|$ vectors $\{u_i\mid i\in S\}$, 
all contained in the plane perpendicular to $w$.
Note that this description also applies to the K\"ahler polygon space $\X$ 
by taking $S=\emptyset$.

By setting $u_i = 0$ for all $i$ we get $\polsa$,
the minimum of the Morse-Bott function $\Phi$ on $U_S$.
On the other hand, consider the submanifold of $U_S$ obtained by imposing
the extra condition that $\|w\| = \sum_{j\in S^c}\|v_j\|$.
Then the first of the two polygons is forced to be linear,
and we are left with $|S|$ vectors $\{u_i\}$ in the perpendicular plane
satisfying a certain norm condition
and adding to zero. 
Identifying this plane with $\C$ and dividing by the circle action rotating this plane, 
we obtain $\C P^{|S|-2}\cong \M_S$,
the maximum of $\Phi$ on $U_S$.
Other critical points of $\Phi$ occur whenever the first polygon is linear,
which is possible for finitely many values of $\|w\|$.
\end{remark}

\begin{usproof}
Suppose given a point $[p,q]_{\R}\in U_S$,
and let
$$u_i = q_ip_i + p_i^*q_i^*\hs\hs\text{ for all }i\in S,$$
$$v_j = (q_j q_j^*)_0\hs\hs\text{ for all }j\in S^c,$$
$$w = \sum_{i\in S}(q_iq_i^*)_0-(p_i^*p_i)_0.$$
These vectors live in $i\cdot\sut\cong\sut^*\cong\R^3$, which is endowed with the
metric
$u\cdot v = \frac{1}{2}\operatorname{tr}uv$,
invariant under the coadjoint action.
With respect to this metric, we have the equalities
$\|(qq^*)_0\| = \half |q|^2$ and $\|(p^*p)_0\| = \half |p|^2$,
hence conditions (1), (2), and (4) are immediate consequences 
of the moment map equations.

To verify condition (3), note that the vectors $\{q_i\mid i\in S\}$
are all proportional over $\C$, which implies that the vectors $(q_iq_i^*)_0$
are positive scalar multiples of each other.
Furthermore, the moment map equation $p_iq_i=0$ implies that $(p_i^*p_i)_0$
is a non-positive scalar multiple of $(q_iq_i^*)_0$, therefore
$w = \sum(q_iq_i^*)_0-(p_i^*p_i)_0$ is proportional over $\R_+$ to
$(q_iq_i^*)_0$ for any $i\in S$.
Then $u_i\cdot w = \half\tr u_i w$ is a multiple of
$$\tr u_i(q_iq_i^*)_0 = \tr u_iq_iq_i^* = \tr p_i^*q_i^*q_iq_i^*
= |q_i|^2\tr p_i^*q_i^* = 0,$$
where the first equality comes from the fact that
$q_iq_i^*-(q_iq_i^*)_0$ is a scalar multiple of the identity,
and $\tr u_i = 0$.

To check condition (5), we first compute the norm of $u_i$:
\begin{eqnarray*}
\|u_i\|^2 &=& \half\tr u_i^2\\
&=& \half\tr(q_ip_ip_i^*q_i^* + p_i^*q_i^*q_ip_i)\\
&=& |q_i|^2|p_i|^2\\
&=& |q_i|^2 (|q_i|^2 - 2\a_i).
\end{eqnarray*}
Since all of the vectors
$\{(q_iq_i^*)_0 -(p_i^*p_i)_0\mid i\in S\}$ point in the same direction,
we have $$\|w\| = \sum_{i\in S}\|(q_iq_i^*)_0\| + \|(p_i^*p_i)_0\|
= \sum_{i\in S} \half|q_i|^2 + \half|p_i|^2
= \sum_{i\in S} |q_i|^2 - \a_i
= \sum_{i\in S} \sqrt{\a_i^2 + \|u_i\|^2}.$$

We have defined a continuous map from $U_S$ to the moduli space of
vectors $\{u_i,v_j,w\}$ satisfying conditions (1)-(5),
and we claim that this map is a homeomorphism.
Since the source of this map is compact and the target is Hausdorff,
it is sufficient to show that the map is bijective.

Suppose given a collection of vectors $\{u_i,v_j,w\}\subs\sut$ 
satisfying conditions (1)-(5).  Using the adjoint action of $SU(2)$,
we may assume that $w$ is a positive scalar multiple of
$\footnotesize{\Big(\begin{array}{cc}
1&0\\
0& -1
\end{array}
\Big)}$.
By condition (3), this implies that for all $i\in S$, there exists $t_i\in\C$
with $u_i =
\footnotesize{\Big(\begin{array}{cc}
0& t_i\\
\bar{t_i}& 0
\end{array}
\Big)}$. 
For $j\in S^c$, we choose $q_j\in\C^2$ with $(q_jq_j^*)_0=v_j$,
and observe that $q_j$ is unique up to the action of $U(1)^n$.
We know that for all $i\in S$, $(q_iq_i^*)_0$ must be a positive multiple of $w$,
hence there exist $a_i,b_i\in\C$ such that
$$q_i=\binom{a_i}{0}\hs\text{ and }\hs p_i = (0\hs b_i)$$
for all $i\in S$.
In order to have $u_i = q_ip_i + p_i^*q_i^*$ and $\half|q_i|^2-\half|p_i|^2=\a_i$, 
we must have $$a_i b_i = t_i\hs\text{ and }\hs 
\half|a_i|^2-\half|b_i|^2=\a_i.$$
These equations uniquely define $a_i$ and $b_i$ up to the action of $U(1)^n$.
It follows from conditions (1)-(5) that 
$(p,q)\in\mur^{-1}(\a)\cap\muc^{-1}(0)$
and that $w=\sum_{i\in S}(q_iq_i^*)_0-(p_i^*p_i)_0$. 
This shows that our map is bijective, and thus completes the proof of Theorem \ref{us}.
\end{usproof}

\vspace{-\baselineskip}
\begin{remark}
Suppose that $S$ has only two elements; without loss of generality
we will assume that $S=\{1,2\}$.  Then forgetting $u_1$ and $u_2$
gives a diffeomorphism from $U_S$ to the ``vertical polygon space''
$V\!P(\a_3,\ldots,\a_n,\a_1+\a_2)$
defined in \cite{HK}, shown to be diffeomorphic to a toric variety.
More generally with $S=\{1,\ldots,k\}$, given any two-element subset $T\subs S$, the subvariety
of $U_S$ given by the equations $u_i=0$ for all $i\in S\smallsetminus T$
is diffeomorphic to $V\!P(\a_{k+1},\ldots,\a_n,\sum_T\a_i)$.
\end{remark}
\end{section}

\begin{section}{Cohomology rings}\label{s1}
In this section we use Theorem \ref{ordinary} to compute the circle-equivariant
cohomology of a hyperpolygon space $\M$, thus reproducing (by different means)
the results of \cite[\S 3]{HP2}.
Recall that we have $$\M = T^*\ctn\mmod G,$$ where $G$ is a quotient
of $U(1)^n\times SU(2)$ by $\Zt$.
We will simplify our computations by dividing first by
the torus $U(1)^n$.
We have
\begin{eqnarray*}
\M &=& \left(T^*\C^{2n}\right)\bigmmod G\\
&\cong&\left(\left(T^*\C^2\right)^n\bigmmod U(1)^n\right)\bigmmod SU(2)\\
&\cong& \prod_{i=1}^nT^*\C P^1\bigmmod SU(2),
\end{eqnarray*}
where the action of $SU(2)$ on each copy of $T^*\C P^1$
is induced by the rotation action on $\C P^1\cong S^2$.

\begin{proposition}\label{hpsurj}
The non-rationalized Kirwan map
$\K_{U(V)}:H^*_{\so\times U(V)}(T^*\C^{2n})\to\hso(\M)$ is surjective.
\end{proposition}

\begin{proof}
The map $\K_{U(V)}$ factors as a composition
$$H^*_{\so\times U(V)}(T^*\C^{2n})
\to H^*_{\so\times SU(2)}\left(\prod_{i=1}^nT^*\C P^1\right)
\overset{\K_{SU(2)}}\Longrightarrow \hso(\M),$$
where the first map is the Kirwan map for a toric hyperk\"ahler variety,
and therefore surjective by \cite{HP1}.  Hence it suffices to show that
$\K_{SU(2)}$ is surjective.

The level set $\muc^{-1}(0)$ for the action
of $SU(2)$ on $\prod_{i=1}^nT^*\C P^1$ is a subbundle of the cotangent bundle,
given by requiring the $n$ cotangent vectors to add to zero after being restricted
to the diagonal $\C P^1$.  In particular this set is smooth,
and its $\so\times SU(2)$-equivariant cohomology ring is canonically
isomorphic to that of $\prod_{i=1}^nT^*\C P^1$.
Hence $\K_{SU(2)}$ factors as
$$H^*_{\so\times SU(2)}\left(\prod_{i=1}^nT^*\C P^1\right)
\cong H^*_{\so\times SU(2)}\Big(\muc^{-1}(0)\Big)\to
\hso\Big(\muc^{-1}(0)\bigmod SU(2)\Big)
\cong \hso(\M),$$ where the map in the middle is the K\"ahler
Kirwan map.
Surjectivity of this map follows from the following more general lemma,
applied to the manifold $\muc^{-1}(0)$.

\begin{lemma}\label{cut}
Suppose given a hamiltonian action of $\so\times G$ on a symplectic
manifold $M$, such that the $\so$ component of the moment map is proper
and bounded below with finitely many critical values.
Then the K\"ahler Kirwan map
$\K:H^*_{\so\times G}(M)\to\hso(M\mod G)$ is surjective.
\end{lemma}

\begin{proof}
Extend the action of $S^1$ to an action on $M\times\C$ by letting
$S^1$ act only on the left-hand factor.  On the other hand,
consider a second copy of the circle, which we will call $\mathbb T$ to avoid confusion,
acting diagonally on $M\times\C$.
Choose $r\in \operatorname{Lie}(\mathbb T)^*\cong\R$ greater than
the largest critical value of the $\mathbb T$-moment map, and consider
the space $$Cut(M\mod G):=\left(M\times\C\right)\mod_{\! r} \mathbb T\times G
\cong\big((M\mod G)\times\C\big)\mod_{\! r} \mathbb T.$$
This space, which is called the {\em symplectic cut} of $M\mod G$ \cite{Le},
is an $S^1$-equivariant (orbifold) compactification of $M\mod G$.
We then have a commutative diagram
$$\begin{CD}
H^*_{\so\times G\times \mathbb T}(M\times\C)        @>>> H^*_{\so\times G}(M)\\
@VVV @VV\K V\\
\hso(Cut(M\mod G))        @>>> \hso(M\mod G).
\end{CD}$$
The vertical map on the left is surjective because the $G\times \mathbb T$ moment map is proper,
and the map on the bottom is surjective because the long
exact sequence in cohomology for $M\mod G\subs Cut(M\mod G)$
splits naturally, hence $\K$ is surjective as well.
\end{proof}

\noindent By applying Lemma~\ref{cut} to $M = \muc^{-1}(0)$,
this completes the proof of Proposition~\ref{hpsurj}.
\end{proof}

\vspace{-\baselineskip}
\begin{remark}\label{horn}
The argument in Proposition \ref{hpsurj} generalizes immediately
to show that the hyperk\"ahler Kirwan map for the quotient
$$\left(\prod_{i=1}^nT^*Flag(\C^k)\right)\bigmmod SU(k)$$
is surjective.  This is itself a quiver variety, and like the
hyperpolygon space, it has a moduli theoretic interpretation.
The K\"ahler quotient
$$\left(\prod_{i=1}^nFlag(\C^k)\right)\bigmod SU(k)$$
is isomorphic to the space of $n$-tuples of $k\times k$
hermitian matrices with fixed eigenvalues adding to zero, modulo
conjugation.
This space has been studied by many authors.  The classical
problem, due to Horn, of determining the values of
the moment map for which it is
nonempty, has only recently been solved \cite{KT}.  For a survey, see \cite{Fu}.
\end{remark}

To compute the kernel of the hyperk\"ahler Kirwan map for the
hyperpolygon space, we first need to study the abelian quotient
$$\NN := \prod_{i=1}^nT^*\C P^1\bigmmod T,$$ where $T\cong U(1)\subs SU(2)$ is a
maximal torus.\footnote{This is the hyperk\"ahler analogue of the {\em abelian
polygon space} from \cite{HK}.}
The space $\prod_{i=1}^nT^*\C P^1$ is a hypertoric variety
given by an arrangement of $2n$ hyperplanes in $\R^n$,
where the $(2i-1)^{\text{st}}$ and $(2i)^{\text{th}}$ hyperplanes
are given by the equations $x_i = \pm\xi_i$.
Taking the hyperk\"ahler quotient by $T$ corresponds on the level of arrangements
to restricting this arrangement to the hyperplane
$\{x\in\R^n\mid\sum x_i = 0\}$.
By Theorem \ref{hsm}, we have
$$\hso\left(\NN\right)\cong
\Q[a_1,b_1,\ldots,a_n,b_n,\delta,x]\Big/
\Big\la a_i-b_i-\delta, \hs a_ib_i\hs\Big{|}\hs i\leq n\Big\ra+
\Big\la A_S, B_S\hs\Big{|}\hs S\text{   short}\Big\ra,
$$
where
\begin{equation*}\label{asbs}
A_S= \prod_{i\in S}(x-a_i)\prod_{j\in S^c}b_j\hs\hs
\text{   and   }\hs\hs B_S = \prod_{i\in S}(x-b_i)\prod_{j\in S^c}a_j.
\end{equation*}
Here $\delta$ is the image in $\hso\left(\NN\right)$
of the unique positive root of $SU(2)$.
The Weyl group $W$ of $SU(2)$, isomorphic to $\Z/2\Z$, acts on this ring by fixing $x$
and switching $a_i$ and $b_i$ for all $i$.
Let $c_i = a_i+b_i$, and let $C_S=A_S+B_S$.
Let $$P=\Q[c_1,\ldots,c_n,\delta,x]\Big/
\Big\la c_i^2-\delta^2\hs\Big{|}\hs i\leq n\Big\ra$$ and
$$Q=P^W=\Q[c_1,\ldots,c_n,\delta^2,x]\Big/
\Big\la c_i^2-\delta^2\hs\Big{|}\hs i\leq n\Big\ra.$$
Let $$\I=\Big\la A_S, B_S\hs\Big{|}\hs S\text{   short}\Big\ra\subs P
\hspace{15pt}\text{and}\hspace{15pt}
\J=\Big\la C_S\hs\Big{|}\hs S\text{   short}\Big\ra\subs Q,$$
so that $$\hso(\NN)\cong P/\I
\hspace{15pt}\text{and}\hspace{15pt}
\hso(\NN)^W\cong Q/\J.$$
Note that all odd powers of $\delta$ in the expression for $C_S=A_S+B_S$ cancel out.
Then by Theorem~\ref{ordinary} and Remark~\ref{notsobad},
\begin{eqnarray*}
\hso(\M)&\cong& \frac{\hso\left(\NN\right)^W}{ann(e)}
\cong \frac{Q}{(e:\J)},
\end{eqnarray*}
where $e=\delta^2(x^2-\delta^2)$, and $(e:\J)$ is the ideal of elements
of $Q$ whose product with $e$ lies in $\J$.

If $S$ is a nonempty short subset, let
$m_S$ be the smallest element of $S$, $n_S$ the smallest
element of $S^c$, and $$D_S = \prod_{m_S\neq i\in S}(c_i-x)
\,\,\cdot\!\prod_{n_S\neq j\in S^c}(c_{n_S}+c_j)\in Q.$$

\begin{theorem}\label{hp}
The circle-equivariant cohomology ring of the hyperpolygon space
$\M$ is isomorphic to\footnote{The class
denoted by $c_i$ in \cite{HP2} differs from our $c_i$ by a sign, hence
to recover the presentation of \cite{HP2} we must replace $c_i-x$
with $c_i+x$ in the expression for $D_S$.}
$$Q\big/\big\la
D_S\mid \emptyset\neq S\text{   short}\big\ra.$$
\end{theorem}

\begin{proof}
We begin by proving that $e\cdot D_S\in\mathcal{J}$ for all nonempty
short subsets $S\subs\{1,\ldots,n\}$.  We will in fact prove the slightly
stronger statement
$$e\cdot D_S\in\Big\la C_T\hs\Big{|}\hs T\subs S\text{   short}\Big\ra
\subs\mathcal{J},$$
proceeding by induction on $|S|$.
We will assume, without loss of generality, that $n\in S$.
The base case occurs when $S=\{n\}$, and in this case
we observe that
$$e\cdot D_S = 2^{n-3}\cdot(x+c_n)\cdot\Big((2x-c_n)\cdot C_{\emptyset} - c_n\cdot C_S\Big).$$
We now proceed to the inductive step, assuming that the proposition
is proved for all short subsets of size less than $|S|$, and all values of $n$.
For all $T\subs S\smallsetminus\{n\}$, we have
$$\half\Big(C_T-C_{T\cup\{n\}}\Big)=(c_n-x)\cdot C'_T,$$
where $C'_T$ is the polynomial in the variables $\{c_1,\ldots,c_{n-1},\delta^2\}$
corresponding to the short subset $T\subs\{1,\ldots,n-1\}$.
Since $S\smallsetminus\{n\}$ is a short subset of $\{1,\ldots,n-1\}$
of size strictly smaller than $S$,
our inductive hypothesis tells us that $e\cdot D_S/(c_n-x)$
can be written as a linear combination of polynomials $C'_T$, where the coefficients
are quadratic polynomials in $\{c_1,\ldots,c_{n-1},\delta^2\}$.
Replacing $C'_T$ with $\half\left(C_T-C_{T\cup\{n\}}\right)=(c_n-x)\cdot C'_T$,
we have expressed $e\cdot D_S$ in terms of the appropriate polynomials.
This completes the induction.

Suppose that $F\in Q$ is an element of degree less than $n-2$
such that $e\cdot F\in\J$.  By the second isomorphism
of Theorem~\ref{ordinary}, this implies that $e'\cdot F\in\I\subs P$,
where $e'=\delta(x^2-\delta^2)$.
Consider the quotient ring $R$ of $P$ obtained by setting
$a_i^2=b_i^2=x=0$ for all $i$.  (Recall that $a_i=\half(c_i+\delta)$
and $b_i=\half(c_i-\delta)$.)  Then element $e'$ maps to zero in $R$,
while the generators $\{A_S, B_S\}$ of $\I$ descend to a basis
for the $n^{\text{th}}$ degree part of $R$.  This means that we
must have $e'\cdot F=0\in P$.
Using the additive basis for $P$ consisting of monomials
that are squarefree in the variables $c_1,\ldots,c_n$,
it is easy to check that $e'$ is not a zero divisor in $P$,
and therefore that $F=0$.

Finally, we must show that $\{D_S\mid \emptyset\neq S\text{  short}\}$
generates all elements of $(e:\J)$ of degree at least $n-2$.
We obtain this fact from the following technical lemma, the proof of which we defer
until the end of the section.

\begin{lemma}\label{tech}
The set $\{D_S\mid \emptyset\neq S\text{  short}\}$ 
descends to a basis for the degree $n-2$
part of the quotient ring $Q/\la x\ra$.
\end{lemma}

Let $F$  be an element of minimal degree $k\geq n-2$ that is in
$(e:\J)$ but not $\la D_S\mid \emptyset\neq S\text{  short}\ra$.
By Lemma \ref{tech}, $F$
differs from an element of
$\la D_S\mid \emptyset\neq S\text{  short}\ra$ by $x\cdot F'$ for some $F'$
of degree $k-1$.  By equivariant formality of $\hso(\M)$,
$$x\cdot F' = F\in (e:\J)\impl F'\in (e:\J),$$
which contradicts the minimality of $k=\operatorname{deg}F$.
Hence $\la D_S\mid \emptyset\neq S\text{  short}\ra=(e:\J)$, and the proposition
is proved.
\end{proof}

\vspace{-\baselineskip}
\begin{corollary}\label{konno}
The ordinary cohomology ring $H^*(\M)$ is isomorphic to
$$\Q[c_1,\ldots,c_n]\Big/\big\la c_i^2-c_j^2\mid i,j\leq n\big\ra
+\la\text{all monomials of degree $n-2$}\ra.$$
\end{corollary}

\begin{proof}
This follows from the fact that $H^*(\M)\cong \hso(\M)/\la x\ra$
for any equivariantly formal space $M$, and the observation in \cite{HP2}
that $\{D_S\mid \emptyset\neq S\text{  short}\}$ descends to a basis for the degree $n-2$
part of $Q/\la x\ra$.
\end{proof}

\begin{prooftech}
Let $\bk = \half(c_1+\ck)$ for all $k$, so that $c_k = 2\bk-d_1$.
The relations $c_k^2 = c_1^2$ translate to $d_k^2 = d_1d_k$ for all $k$,
and we have
$$Q/\la x\ra = \Q[d_1,\ldots,d_n]\Big/\big< d_k^2-d_1d_k\hs\big{|}\hs k\in\twn\big>.$$
For all short subsets $S$, put
$$\Sb = S\smallsetminus\{m_S\}\hspace{15pt}\text{and}\hspace{15pt}
\Lb = S^c\smallsetminus\{n_S\},$$
and let
$$v_S = (-1)^n\pjl (\bj+\bns-d_1)\times\pis(2\bi-d_1).$$
For all $A\subs\twn$, let
$$\ba = (-1)^{|A|} d_1^{n-2-|A|} \prod_{k\in A} d_k$$ for all $A \subsetneq \twn$.
Then $\{d_A\}$
is a basis for
the $(n-2)^{\text{nd}}$ graded piece of $Q/\la x\ra$,
and $v_S$ is equal to $(-1)^n\cdot 2^{-|\Lb|}$ times the image of $D_S$ in $Q/\la x\ra$.
Hence our the statement of Lemma \ref{tech} is that for all $A$, $\ba$ may be 
expressed as a linear combination of the elements $\{v_S\mid S\in\sa\}$.

\begin{convention}
The notation $S^c$ refers to the complement of $S$ inside of the set $\otn$,
while the notation $A^c$ refers to the complement of $A$ inside of the set $\twn$.
\end{convention}

\begin{claim}\label{vs}
We have the following expression for $v_S$ in terms of the basis $\{d_A\}$:
$$\vs =
\begin{cases}
\displaystyle{\sum_{\ss{\Lb\hs\subseteq A \\ m_S\notin A}}} \hs 2^{|\asb|} \hs\ba
& \text{if}\hspace{10pt}1 \in S^c;\vspace{.5cm}\\
\displaystyle{\sum_{S^c \nsubseteq A}} \hs 2^{|\asb|} 
\hs\ba & \text{if}\hspace{10pt}1 \in S.\\
\end{cases}
$$
\end{claim}

\begin{proof}
Any degree $n-2$ monomial
in $d_1,\ldots,d_n$ is equal to $(-1)^{|A|}d_A$, where $A$ is the set of $k>1$ such that
$d_k$ appears in the monomial.
Expanding $\vs$, we need to count (with sign) the occurrence of $d_A$ for each $A$.
In most cases we find that there is no cancellation, and the claim is straightforward.
The most difficult case occurs when
$1\in S$ (therefore $n_S=1$) and $\ms\in A$; 
in this case the number of times (with multiplicity) that $d_A$
occurs in $\vs$ is 
\begin{eqnarray*}
&& (-1)^n(-1)^{|A|}(-1)^{|\acsb|}\hs 2^{|\asb|}\sum_{E\subsetneq\Ac\cap\Lb}(-1)^{|E|}
\vspace{.3cm}\\
\vspace{.3cm}&=& (-1)^n(-1)^{|A|}(-1)^{|\acsb|}\hs 2^{|\asb|}
\left((1-1)^{|\Ac\cap\Lb|}-(-1)^{|\Ac\cap\Lb|}\right)\vspace{.3cm}\\
\vspace{.3cm}&=& (-1)^{n+|A|+|\acsb|+|\Ac\cap\Lb|+1}\hs 2^{|\asb|}\\
\vspace{.3cm}&=& (-1)^{2n}\hs 2^{|\asb|}\\
\vspace{.3cm}&=& \hs 2^{|\asb|}.
\end{eqnarray*}
We leave the remaining cases to be checked by the reader.
\end{proof}

\vspace{-\baselineskip}
\begin{claim}\label{ws}
Suppose that $1\in S$.
Let $S_0 = S$, and for $1\leq k\leq |S|$, let $S_k = S_{k-1}\setminus\{m_{S_{k-1}}\}$.
(In other words, let $S_k$ consist of the $|S|-k$ largest elements of $S$).
Then $$v_S + \displaystyle{\sum_{k=1}^{|S|-1}} \hs 2^{k-1}\hs v_{S_k} 
= \displaystyle{\sum_A}\hs 2^{|\asb|}\hs d_A.$$
\end{claim}

\begin{proof}
We proceed by induction to show that
$$v_S + \displaystyle{\sum_{k=1}^l} \hs 2^{k-1}\hs v_{S_k} =
\displaystyle{\sum_A} \hs 2^{|\asb|}\hs d_A - 
\hs 2^l\cdot\!\displaystyle{\sum_{\overline{S^c_{l+1}}\subs A}}
2^{|\asb_l|}\hs d_A.$$
The case $l=|S|-1$ is the statement of the claim.
The base case $l=0$ follows from Claim \ref{vs},
together with the observation that $\overline{S^c_{1}}=S^c$.
More generally, for all $l\geq 1$, we have 
$$\overline{S^c_{l+1}} = S^c\cup\{m_{S_1},\ldots,m_{S_l}\}.$$
Then 
\begin{eqnarray*}
v_S + \displaystyle{\sum_{k=1}^{l+1}} \hs 2^{k-1}\hs v_{S_k} &=&
v_S + \displaystyle{\sum_{k=1}^l}\hs 2^{k-1}\hs v_{S_k} + 2^l\hs v_{S_{l+1}}\\
&=& \displaystyle{\sum_A} \hs 2^{|\asb|}\hs d_A - 2^l\cdot\!\displaystyle{\sum_{\overline{S^c_{l+1}}\subs A}}
2^{|\asb_l|}\hs d_A + 2^l\cdot\!\sum_{\ss{\overline{S^c_{l+1}}\subs A \\ m_{S_{l+1}}\notin A}}
2^{|\asb_{l+1}|}\hs d_A
\end{eqnarray*}
by the inductive hypothesis and Claim \ref{vs}.
Using the fact that
$\asb_{l+1}=\asb_l$ when $m_{S_{l+1}}\notin A$,
this is equal to
$$\sum_A\hs 2^{|\asb|}\hs d_A - 2^l\cdot\!\sum_{\overline{S^c_{l+1}}\cup\{m_{S_{l+1}}\}\subs A}
2^{|\asb_l|}.$$
Finally, since $|\asb_{l+1}|=|\asb_l|-1$ when $m_{S_{l+1}}\in A$,
this reduces to
$$\sum_A \hs 2^{|\asb|}\hs d_A - 2^{l+1}\hs\cdot\!\displaystyle{\sum_{\overline{S^c_{l+2}}
\subs A}}
2^{|\asb_{l+1}|}\hs d_A,$$ thus proving our claim.
\end{proof}

For all short subsets $T$ containing $1$,
let $w_T=\displaystyle{\sum_A}2^{|A\cap \bar{T}|}\hs d_A$, 
which by Claim \ref{ws} is expressible as a linear combination of elements
of the set $\{v_S\mid \emptyset\neq S\in\mathcal{S}\}$.
Let $$x_S = 
\begin{cases}
\displaystyle{\sum_{1\in T\subs S}}(-1)^{|S|+|T|}w_T & \text{if $1\in S$,}\\
v_S & \text{if $1\in S^c$.}
\end{cases}$$
Our last task will be to prove that the transition matrix $\Upsilon$ 
taking the basis $\{d_A\}$ to 
the set $\{x_S\}$
is lower triangular with ones on the diagonal, and therefore invertible.
In order to make sense of ``the diagonal," 
we must first give an explicit bijection between the set of proper subsets of $\twn$
and the set of nonempty short subsets of $\otn$.
We do this as follows:  given $A\subsetneq\twn$, let 
$$S(A) = \begin{cases}
\Ac &\text{ if $\Ac$ is short,}\\
\otn\setminus\Ac = A\cup\{1\} &\text{ if $\Ac$ is long.}
\end{cases}$$
The rows of $\Upsilon$ will be indexed by $A$, 
and the sets will appear in lexicographic order within cardinality
class.  For example, when $n=4$, the order of the rows will be 
$\emptyset$, $\{2\}$, $\{3\}$, $\{4\}$,
$\{2,3\}$, $\{2,4\}$, $\{3,4\}$.  
The columns will be indexed by $S$ according to the bijection described above.

\begin{claim}\label{tri}
The matrix $\Upsilon$ is lower triangular with ones on the diagonal.
\end{claim}

\begin{proof}
First consider a column corresponding to a short subset $S$ that does {\it not} contain 1.
The entries in this column correspond to the coefficient of $d_A$ in $x_S = v_S$.
From Claim \ref{vs}, we see that $d_A$ appears in $v_S$ only if $\Lb\subs A\subs \Lb\cup\Sb$,
and if so it appears with a coefficient of $2^{|\asb|}$.
Since $1\notin S$, we have $\Lb = S^c\setminus\{1\} = \twn\setminus S$.
The diagonal entry corresponds to the set $A = \twn\setminus S = \Lb$, therefore
in this row we get the number $2^{|\asb|} = 2^{|\Lb\cap\Sb|} = 1$.
Since the set $A$ corresponding to a given row can never contain the set $B$ corresponding
to a 
lower row, the rows above the diagonal fail to satisfy the condition $\Lb\subs A$,
and we get all zeros.

Now consider a column corresponding to a short subset $S$ that {\it does} contain $1$.
In this case, the coefficient of $d_A$ in $x_S$ is
$$(-1)^{|S|}\displaystyle{\sum_{1\in T\subs S}}(-1)^{|T|}2^{|A\cap\bar{T}|}.$$
The diagonal entry corresponds to the set $A=\Sb$, and we get
\begin{eqnarray*}
(-1)^{|S|}\displaystyle{\sum_{1\in T\subs S}}(-1)^{|T|}2^{|\bar{T}|}
&=& (-1)^{|\Sb|}\displaystyle{\sum_{1\in T\subs S}}(-2)^{|\bar{T}|}\\
&=& (-1)^{|\Sb|}(1-2)^{|\Sb|} = 1.
\end{eqnarray*}
Any row above the diagonal corresponds to a set $A$ which does not contain $\Sb$.
Choose an element $l\in\Sb\setminus A$.  Then
\begin{eqnarray*}
(-1)^{|S|}\displaystyle{\sum_{1\in T\subs S}}(-1)^{|T|}2^{|A\cap\bar{T}|}
&=& (-1)^{|S|}\displaystyle{\sum_{l\in T}}(-1)^{|T|}2^{|A\cap\bar{T}|}
+ (-1)^{|S|}\displaystyle{\sum_{l\notin T}}(-1)^{|T|}2^{|A\cap\bar{T}|}\\
&=& (-1)^{|S|}\displaystyle{\sum_{l\notin T}}\left[(-1)^{|T|}2^{|A\cap\bar{T}|}
+ (-1)^{|T\cup\{l\}|}2^{|A\cap\bar{T}|}\right]\\
&=& 0.
\end{eqnarray*}
Hence $\Upsilon$ is lower triangular.
\end{proof}

Claim \ref{tri} tells us that each $d_A$ can be expressed as a linear
combination of elements of the form $x_S$, and therefore of elements
of the form $v_S$.  This completes the proof of Lemma \ref{tech}
\end{prooftech}
\end{section}

\begin{section}{Cohomology of the core components}\label{last}
In this section we compute the $S^1$-equivariant 
and ordinary cohomology rings of the core component
$U_S$ corresponding to a short subset $S\subs\otn$.  Since $U_S$ is the closure of a cell
in an even-dimensional equivariant cellular decomposition of $\M$, the restriction
map $\hso(\M)\to\hso(U_S)$ is surjective.  In particular, $\hso(U_S)$ 
is generated by restrictions
of the Kirwan classes $c_1,\ldots,c_n,x$.  
For our presentation, it will be convenient
to assume that $1\in S$, and to work with the classes 
$d_k = \half(c_1+c_k)$ introduced in the proof of Lemma \ref{tech}.  
With respect to these generators, we obtain the following result.

\begin{theorem}\label{eqcore}
The equivariant cohomology ring $\hso(U_S)$ is isomorphic to
$\Q[d_1,\ldots,d_n,x]/\mathcal{J}_S,$
where $\mathcal{J}_S$ is generated by the following four families:
\begin{eqnarray*}&1)&\hs\hs d_1 - d_i \hs\hs\text{ for all }\hs\hs i\in S\\
&2)&\hs\hs
d_j(d_1 - d_j)
\hs\hs\text{ for all }\hs\hs j\in S^c\\
&3)&\hs\hs \prod_{j\in R}d_j\hs\hs\text{ for all }\hs\hs R\subs S^c\text{ such that }R\cup S
\text{ is long}\\
&4)&\hs\hs (d_1 + x)^{|S|-1}\cdot\frac{1}{d_1}\left(\prod_{j\in L}(d_j-d_1)-\prod_{j\in L}d_j\right)
\hs\hs\text{ for all long subsets }L\subs S^c.
\end{eqnarray*}
\end{theorem}

\begin{corollary}\label{ordcore}
The ordinary cohomology ring $H^*(U_S)$ is isomorphic to
$\Q[d_1,\ldots,d_n]/\mathcal{I}_S,$
where $\mathcal{I}_S$ is generated by the following four families:
\begin{eqnarray*}&1)&\hs\hs d_1 - d_i \hs\hs\text{ for all }\hs\hs i\in S\\
&2)&\hs\hs
d_j(d_1 - d_j)
\hs\hs\text{ for all }\hs\hs j\in S^c\\
&3)&\hs\hs \prod_{j\in R}d_j\hs\hs\text{ for all }\hs\hs R\subs S^c\text{ such that }R\cup S
\text{ is long}\\
&4)&\hs\hs d_1^{|S|-2}\prod_{j\in L}(d_j-d_1)
\hs\hs\text{ for all long subsets }L\subs S^c.
\end{eqnarray*}
\end{corollary}

\begin{remark}\label{integer}
Each of these relations has a geometric interpretation.
For $i\in\{1,\ldots,n\}$, it is possible to construct a line bundle on $\M$ with equivariant
Euler class $d_i - d_1$
which has a section supported on the locus
where $q_1q_1^*$ and $q_iq_i^*\in\R^3$ point in opposite directions.  
Since this locus is disjoint
from $U_S$ when $i\in S$, we have $$1)\hs\hs d_i = d_1 \in \hso(\us)\hs\text{ for all }\hs i\in S.$$
Similarly,
$-d_j = -\half(c_1+c_j)$ is represented by the divisor $Z_{1j}\subs\M$ of points
on which $q_1q_1^*$ and $q_iq_i^*\in\R^3$ point in the same direction
\cite[\S 3]{HP2}.
Then by the previous reasoning, we obtain 
$$2)\hs\hs d_j(d_1 - d_j) = 0\in\hso(\us)\hs\text{ for all }\hs j\in S^c.$$
For any $R\subs S^c$, we may intersect the divisors $Z_{ij}\subs\M$ 
(defined in the analogous way) for all $j\in R$
to find that
the cohomology class $(-1)^{|R|}\prod_{j\in R}d_j$ is represented
by the subvariety $Z_R\subs \M$ of points with $q_j$ proportional to $q_1$
for all $j\in R$.
When restricted to $U_S$, this becomes $U_S\cap U_{R\cup S}$,
the closure of the unstable manifold for the critical locus $\M_{R\cup S}\cap U_S$
of the Morse-Bott function $\Phi|_{U_S}$.
In particular, we have $$3)\hs\hs \prod_{j\in R}d_j=0\in\hso(\us)\text{ if }R\cup S\hs\text{ is long}.$$
To understand the fourth family of relations, we note that the class
$$d_1+x = 2d_i-d_1+x = c_i+x \in \hso(U_S)$$ 
is represented by the divisor $W_i$ of points
with $p_i = 0$ for any $i\in S$ \cite[\S 3]{HP2}.  In particular,
$(d_1+x)^{|S|-1}$ 
is represented by the subvariety of points in $U_S$ on which $p_i = 0$
for all $i\in\Sb$, which is equal to $\X_S$ by the complex moment map condition.
Hence the fourth family of generators of $\mathcal{J}_S$ (or of $\mathcal{I}_S$)
can be interpreted geometrically as $(d_1+x)^{|S|-1}$ 
(respectively $d_1^{|S|-1}$
in the nonequivariant case) times
classes that vanish in $\hso(\X_S)$ (see Lemma \ref{pols}).
%
\end{remark}

\begin{eqcoreproof}
Let $\phi:\Q[d_1,\ldots,d_n,x]\to\hso(\us)$ denote the composition of the Kirwan
map with restriction to $\us$.  Our claim is that $\ker\phi = \mathcal{J}_S$.
For every short subset $T$ containing $S$,
let $$\phi_T:\Q[d_1,\ldots,d_n,x]\to\hso(\mtus)$$ denote
the composition of the Kirwan map with restriction to $\mtus$,
and let $$J_T = \ker\phi_T.$$
Similarly, let $$\phi_{\emptyset}:\Q[d_1,\ldots,d_n,x]\to\hso(\X_S)$$
be the natural map, and let
$$J_{\emptyset} = \ker\phi_{\emptyset}.$$
The kernel of the restriction map $\hso(U_S)\to\hso(U_S^{S^1})$
to the fixed point set of $U_S$ is a torsion module over $\hso(pt)$ \cite[3.5]{AB}, 
and Proposition \ref{formality}
tells us that $\hso(U_S)$ is a free $\hso(pt)$-module. 
Hence the restriction map is injective,
and we have $$\ker\phi = \ker\phi_{\emptyset}\cap\bigcap_{\ts}\ker\phi_T.$$
We know that $\mtus\cong\C P^{|S|-2}$ for all short $\ts$,
therefore $$\hso(\mtus)\cong\Q[h,x]/h^{\sz}.$$
Furthermore, we know that for all $i\in T$, the restriction of $d_i + x$ 
to $\hso(U_T)$ is represented by the divisor $W_i\cap U_T$ (see Remark \ref{integer}),
and therefore restricts to the class of a hyperplane on $\mtus$.
Hence $\phi_T(d_i+x)=h$ for all $i\in T$.
On the other hand, for $j\in T^c$, the class $d_j$ 
is represented by the divisor $Z_{1j}$ on $\M$, which is disjoint from $\mtus$,
hence $\phi(d_j) = 0$ for all $j\in T^c$.
Thus we conclude that
$$\ker\phi_T = \<d_1-d_i, \hs d_j, (d_1+x)^{\sz}\mid i\in T, j\in T^c\>.$$

\begin{lemma}\label{jt}
The intersection 
$\displaystyle{
\bigcap_{\ts}}\ker\phi_T$
is equal to
$$\<d_1-d_i, \hs d_j(d_1 - d_j), 
\displaystyle{\prod_{j\in R}d_j},
(d_1+x)^{\sz}
\hs\hs\bigg |\hs\hs i\in S, j\in S^c, R\cup S\text{ long}\>.$$
\end{lemma}

\begin{proof}
First, since the variable $x$ appears only in the generator $(d_1+x)^{\sz}$,
which is contained in every ideal on both sides of the equation, we may reduce
the problem to showing that
\begin{equation}\label{reduce}
\bigcap_{\ts}\big\langle d_1-d_i, \hs d_j\bigmid i\in T, j\in T^c\big\rangle = 
\<d_1-d_i, \hs d_j(d_1-d_j), \prod_{j\in R}d_j
\hs\hs\bigg |\hs\hs i\in S, j\in S^c, R\cup S\text{ long}\>
\end{equation}
in the ring $\Q[d_1,\ldots,d_n]$.
Both ideals cut out the (reducible) variety 
$$\bigcup_{\ts}Y_T\subs\operatorname{Spec}\Q[d_1,\ldots,d_n],$$
where $$Y_T = \big\{(z_1,\ldots,z_n\bigmid
z_i=z_1\hs\forall i\in S, z_j=0\hs\forall j\in S^c\big\}.$$
The left hand side of Equation \eqref{reduce} is an intersection of prime ideals,
and is therefore radical.
Thus by Hilbert's Nullstellensatz, it is sufficient to prove that
the right hand side of Equation \eqref{reduce} is radical.
This involves showing that the ideal is saturated, with
Hilbert polynomial equal to the constant $\#\{\text{short }\ts\}$.

The degree $k$ piece of the quotient
$$\Q[d_1,\ldots,d_n]/\<d_1-d_i, \hs d_j(d_1-d_j)
\mid i\in S, j\in S^c\>$$
has a basis of elements of the form
$$d_1^{e_1}\prod_{j\in S^c}d_j^{e_j},$$
where $e_j \in\{0,1\}$ for all $j>0$, and $e_1+\sum_{j\in S^c}e_j = k$.
The subset of these elements with the property that
$S\cup\{j\mid e_j = 1\}$ is short descends to a basis for the
degree $k$ part of the ring
$$\Q[d_1,\ldots,d_n]\bigg/\<d_1-d_i, \hs d_j(d_1-d_j), \prod_{j\in R}d_j
\hs\hs\bigg |\hs\hs i\in S, j\in S^c, R\cup S\text{ long}\>,$$
hence our ideal has the desired Hilbert polynomial.
It is also clear from this description that if an element
$a$ of the quotient ring is nonzero, so is $d_1^d\cdot a$ for any $d\geq 0$,
hence our ideal is saturated.
\end{proof}

It now remains to show that
$$\mathcal{J}_S = \<d_1-d_i, \hs d_j(d_1 - d_j), \prod_{j\in R}d_j,
\hs (d_1+x)^{\sz} \hs\hs\bigg |\hs\hs i\in S, j\in S^c, R\cup S\text{ long}\> \cap \ker\phi_{\emptyset}.$$
The fact that $\mathcal{J}_S$ is contained in the intersection is clear.
To show the opposite containment, consider an element
$$a + \eta\cdot (d_1+x)^{\sz}\in \<d_1-d_i, \hs d_j(d_1 - d_j), \prod_{j\in R}d_j,
\hs (d_1+x)^{\sz} \hs\hs\bigg |\hs\hs i\in S, j\in S^c, R\cup S\text{ long}\>,$$
with $$a\in \<d_1-d_i, \hs d_j(d_1 - d_j), \prod_{j\in R}d_j\hs\hs\bigg |\hs\hs
i\in S, j\in S^c, R\cup S\text{ long}\>,$$
and suppose that we also have
$$a + \eta\cdot (d_1+x)^{\sz}\in\ker\phi_{\emptyset}.$$

\begin{lemma}\label{pols}{\em \cite{HK}}
The kernel of $\phi_{\emptyset}$
is equal to 
$$\<d_1-d_i, \hs d_j(d_1-d_j), \prod_{j\in R}d_j, 
\hs (d_1+x)^{\sz}d_1^{-1}\left(\prod_{j\in L}(d_j-d_1)-\prod_{j\in L}d_j\right)\>,$$
where $i\in S, j\in S^c$, and $R,L\subs S^c$, with $R\cup S$ and $L$ both long.
\end{lemma}

Lemma \ref{pols} tells us that $a\in \ker\phi_{\emptyset}$,
therefore $$\eta\cdot (d_1+x)^{\sz}\in\ker\phi_{\emptyset}.$$
But $(d_1+x)^{\sz}$ is represented in $\hso(U_S)$ by the subvariety
$\X_S$ (see Remark \ref{integer}), hence 
$$0 = \phi_{\emptyset}(\eta\cdot (d_1+x)^{\sz})
= \phi_{\emptyset}(\eta)\cdot e(\X_S),$$
where $e(\X_S)$ is the equivariant Euler class of the normal bundle
to $\X_S$ inside of $U_S$.  Since the equivariant Euler class
of the normal bundle to a component of the fixed point set is never a zero-divisor,
we have $\eta\in\ker\phi_{\emptyset}$.
Then by Equation \ref{pols}, $$a+\eta\cdot (d_1+x)^{\sz}\in\mathcal{J}_S.$$
This completes the proof of Theorem \ref{eqcore}.
\end{eqcoreproof}

\vspace{-\baselineskip}
\begin{example}\label{projective}
For arbitrary $n$ and $\a$, suppose that $S$ is a maximal short subset.
Then Corollary \ref{ordcore} tells us that
$H^*(U_S)\cong\Q[d_1]/\langle d_1^{n-2}\rangle$.
We conjecture that in this case we in fact have $U_S\cong \C P^{n-3}$.
\end{example}

\begin{example}\label{again}
Consider the core component pictured in Example \ref{ooh}.
By Theorem \ref{ordcore} and Remark \ref{integer},
$$H^*(U_S) \cong \Q[d_1,d_3,d_4,d_5]
\Bigg/
\left<\begin{array}{c}
d_3(d_1-d_3),\hs d_4(d_1-d_4),\hs d_5(d_1-d_5),\hs d_3d_4,\hs d_3d_5,\hs d_4d_5,\\
d_1(d_1-d_3-d_4),\hs  d_1(d_1-d_3-d_5),\hs d_1(d_1-d_4-d_5)
\end{array}\right>,$$
where $d_1$ is the fundamental class of $\X_S$,
and $d_3$, $d_4$, and $d_5$ are the negatives of the fundamental classes of the curves
labeled $123$, $124$, and $125$, respectively.
Because the transverse intersection of two complex varieties is positive,
we know that $-d_1d_3[U_S] = 1$.
With respect to the basis $$\{d_1-d_3-d_4-d_5,\, d_3,\, d_4,\, d_5\},$$
the intersection form on $H^2(U_S)$ is represented by the matrix
$$\left(\begin{array}{cccc}
1&&&\\
&-1&&\\
&&-1&\\
&&&-1\end{array}\right).$$
It is likely that $U_S$ is isomorphic 
to the blow-up of $\C P^2$ at three points.
\end{example}

\begin{example}
Using the same $\a= (1,1,3,3,3)$,
consider the short subset $S=\{1,3\}$.
In this case, Theorem \ref{ordcore} tells us that
$$H^*(U_S)\cong\Q[d_1,d_2]/\left<d_1^2, d_2(d_1-d_2)\right>.$$
With respect to the basis $\{d_1-d_2,\, d_2\}$,
the intersection form on $H^2(U_S)$ is represented by the matrix
$\footnotesize{\Big(\begin{array}{cc}
-1&0\\
0& 1
\end{array}
\Big)}$,
hence $U_S$ is homeomorphic to the blow-up of $\C P^2$ at a single point.
\end{example}
\end{section}
\end{chapter}

\footnotesize{

}
\end{document}